\documentclass[12pt]{article}
\usepackage{graphicx}

\usepackage{amsmath}
\usepackage{amssymb}

\usepackage{comment}
\usepackage{verbatim}
\usepackage{afterpage} 

\usepackage{ifthen}
\newboolean{vmcode}  
\setboolean{vmcode}{false}
\newcommand{\cmdvmcode}[2]{\ifthenelse{\boolean{vmcode}}{\verbatiminput{#1}}{#2}}

\newcommand{\PP}{L}
\newcommand{\QQ}{Q}

\newcommand{\iK}{ Part II }

\newtheorem{theorem}{Theorem}

\begin{document}


\numberwithin{equation}{section}

\title{
Torsional rigidity for tangential polygons
%
}

\author{Grant Keady\\
\texttt{Grant.Keady@curtin.edu.au}}
\date{\today}

\maketitle
\vspace{-0.5cm}


\section*{Abstract}
{An inequality on torsional rigidity is established.
For tangential polygons this inequality is stronger than
an inequality of Polya and Szego for convex domains.
(A survey of related work, not in the journal submission, is presented.)
}


\section{Introduction}\label{sect:Intro}

In the most general situation $\Omega$ is a simply-connected plane domain.
However as we wish to compare our inequalities with earlier results proved for convex domains,
notably the Polya-Szego inequality (involving equation~(\ref{eqt:Bdef})), 
we have in mind convex domains.
Although Theorem~\ref{thm:thm1} concerns more general domains,
our main result derived from it, Theorem~\ref{thm:thm2},  concerns tangential polygons.
We will denote the area of $\Omega$, $|\Omega|$, by $A$, and 
its perimeter $|\partial\Omega|$ by $L$.

\subsection{The torsion pde, Problem (P(0))}\label{subsect:P0}

The elastic torsion problem is to find a function $u_0$ satisfying
$$ - \frac{\partial^2 u_0}{\partial x^2}- \frac{\partial^2 u_0}{\partial y^2}
 = 1\, {\rm in}\,\Omega ,\,\qquad
u_0
= 0 \, {\rm on}\,\partial\Omega \, .
\eqno({\rm P}(0))
$$
Define the {\it torsional rigidity} of $\Omega$ as,
\begin{equation}
\QQ_0
:= \int_\Omega u_0 .
\label{eqt:Qdef}
\end{equation}
(This differs by a factor of 4 from definitions elsewhere, e.g in~\cite{PoS51}.)
There are many identities and inequalities concerning $\QQ_0$.
For example, one identity is
\begin{equation}
\QQ_0
= \frac{1}{4}\, \int_{\partial\Omega} (x.n) \left(\frac{\partial u_0}{\partial n}\right)^2 .
\label{eqt:bndryQ}
\end{equation}
Define (for convex $\Omega$), in the notation of~\cite{PoS51},
\begin{equation}
B_\Omega
= \int_{\partial\Omega} \frac{1}{x\cdot{n}} .
\label{eqt:Bdef}
\end{equation}
An inequality, which we will call the Polya-Szego inequality, is
\begin{equation}
\QQ_0
\ge \frac{A^2}{4\ B_\Omega}  
( = {\rm{ \ for\ tangential\ polygons\ }} \frac{A^3}{2L^2} )
 \label{in:PS}
\end{equation}

Define 
$$ J(v) = \int_\Omega \left( 2v - |\nabla v|^2\right) .$$
The solution $u_0$ maximizes $J$ over functions vanishing on the boundary $\partial\Omega$.
(More precisely over functions $v\in{{\mathring W}^1_2(\Omega)}$.) 

A quantity occuring in our Theorem~\ref{thm:thm1} is
\begin{equation}
\QQ_1
=  \int_{\partial\Omega}  \left(\frac{\partial u_0}{\partial n}\right)^2 .
\label{eqt:Q1def}
\end{equation}
For tangential polygons from equation~(\ref{eqt:bndryQ}) we note that
there is a very simple equation relating $\QQ_1$ and $\QQ_0$: this is
needed for our Theorem~\ref{thm:thm2}.

\subsection{Problem (P($\infty$))}\label{subsect:Pinf}

Problem (P($\infty$) was defined in the statement of Theorem 2.2 of~\cite{KM93}.
Our notation here is as in that paper.
where $u_\infty$ solves
$$ - \frac{\partial^2 u_\infty}{\partial x^2}- \frac{\partial^2 u_\infty}{\partial y^2}
 = 1\, {\rm in}\,\Omega ,\,
\frac{\partial u_\infty}{\partial n}
= -\frac{|\Omega|}{|\partial\Omega|} \, {\rm on}\,\partial\Omega \,
{\rm and}\,  \int_{\partial\Omega} u_\infty = 0 .
\eqno({\rm P}(\infty))
$$
Define, as in~\cite{KM93} equations~(4.9) and (4.11),
\begin{equation}
\Sigma_\infty
= \int_\Omega u_\infty , \qquad  {\rm and\ \ } \Sigma_1= - \int_{\partial\Omega} u_\infty^2 ,
\label{eqt:SigDef}
\end{equation}
with $u_\infty$ satisfying Problem (P($\infty$)).

Once again there is a variational characterisation of the solutions.
This time $u_\infty$ is the maximizer of $J(v)$ as one varies over functions $v$ 
for which the integral of $v$ over the boundary $\partial\Omega$ is zero.
The variational approach can be used to establish the inequality
\begin{equation}
\Sigma_\infty \ge \QQ_0 .
\label{KM934p9}
\end{equation}
(See also~\cite{KM93} equation~(4.9) for a different approach.)

The solution $u_\infty$ when $\Omega$ is a tangential polygon is given in~\cite{KM93}
equations~(2.14) and (2.15).

\section{Relations between $u_0$ and $u_\infty$}

We already have~(\ref{KM934p9}) relating $\Sigma_\infty$ and $\QQ_0$.
The next result involves $\Sigma_1$ and $\QQ_1$.
(We haven't exp[ored to what extent the boundary smoothness might be relaxed.
We need to apply the Divergence Theorem.)

\begin{theorem}{\label{thm:thm1}
For (convex) domains $\Omega$ with Lipschitz boundary which is also piecewise $C^1$,
the following inequality is satisfied:
\begin{equation}
\frac{A^2}{L}-\frac{(\Sigma_\infty-\QQ_0)^2}{\Sigma_1}\le \QQ_1  .
\label{eqt:QQ1R1a}
\end{equation}
}
\end{theorem}
\par\noindent
Recall that $\Sigma_1<0$, so both terms on the left-hand side are positive.

\smallskip
\par\noindent{\it Proof.}
Considering ${\rm div}(u_0\nabla u_\infty)$ the Divergence Theorem gives
$$ \QQ_0 = \int_\Omega \nabla u_0\, \cdot\, \nabla_\infty .$$
Considering ${\rm div}(u_\infty\nabla u_0)$ the Divergence Theorem gives
$$\int_{\partial\Omega} u_\infty \frac{\partial u_0}{\partial n}
=    - \Sigma_\infty + \int_\Omega \nabla u_0\, \cdot\, \nabla u_\infty .
$$
These combine to give
$$\Sigma_\infty - \QQ_0 = \int_{\partial\Omega} u_\infty\left(- \frac{\partial u_0}{\partial n}\right) .$$
We now introduce an arbitrary constant $c$ and subtract $c\,A$ from both sides, giving
\begin{equation}
(\Sigma_\infty - \QQ_0) -cA
=  \int_{\partial\Omega} (u_\infty-c)\, \left(- \frac{\partial u_0}{\partial n}\right) .
\label{eqt:ceq}
\end{equation}
We now use the Cauchy-Schwarz inequality on the right-hand side to give
$$\left((\Sigma_\infty - \QQ_0) -cA\right)^2 \le
\QQ_1 \,  \int_{\partial\Omega} (u_\infty-c)^2 ,
$$
where we have used equation~(\ref{eqt:QQ0t}) for one of the integrals on the right-hand side.
Now, on using that the integral of $u_\infty$ around the boundary is zero, we have
$$  \int_{\partial\Omega} (u_\infty-c)^2 
= -\Sigma_1 + c^2 L .$$
Thus, for all real $c$
\begin{equation}
 \frac{\left((\Sigma_\infty - \QQ_0) -cA\right)^2}{-\Sigma_1 + c^2 L}
\le  \QQ_1 
\label{in:cin}
\end{equation}

The function of $c$ on the left is clearly nonnegative and bounded,
tending to $A^2/L$ as c tends to both plus and minus infinity.
It has two critical points: the one making the function 0 is clearly the minimum.
The other is at $c=c_*$ where
$$ c_* =  \frac{A\Sigma_1}{L (\Sigma_\infty - \QQ_0)} . $$
Subsitituting $c_*$ for $c$ in inequality~(\ref{in:cin}) yields the result of the theorem.
\hfill$\square$
\medskip

\section{Tangential polygons}\label{sect:Tang}

\subsection{Geometric results}\label{subsect:tangGeom}

A tangential polygon, also known as a circumscribed polygon, is a convex polygon that contains an inscribed circle (also called an incircle), a circle that is tangent to each of the polygon's sides.
For any (convex) tangential polygon, the area $A$, perimeter $L$ and inradius $\rho$
are related by
$$  A=\frac{1}{2} \rho L . $$
We always choose the origin of our coordinate system to be at the incentre.
There is some literature on tangential polygons, e.g.~\cite{AM04,RP01}.
Any triangle is a tangential polygon.
Quadrilaterals which are tangential include kites and hence rhombi.

Concerning $B$, defined in~\cite{PoS51} and here at equation~(\ref{eqt:Bdef}), another simple identity for (convex) tangential polygons is
$$ B = \frac{L}{\rho} ,$$
(and $B\ge{2\pi}$ with equality only for the disk,
and for any triangle, $B\ge{6\sqrt{3}}$):
see~\cite{Ai58}.

The quantities $\Sigma_\infty$ and $\Sigma_1$ can be expressed in terms of boundary moments
$i_{2k}$ -- moments about the incentre -- defined by
$$ i_{2k} = \int_{\partial\Omega} (x^2+y^2)^k . $$
We remark that the Cauchy-Schwarz inequality for the integrals implies that
$i_4\ge{i_2^2/L}$.

\medskip

{For methods to calculate $i_{2k}$  in terms of the tangential polygon's inradius
and angles, see~\cite{Ke20i}.}

\subsection{The new inequality for tangential polygons}\label{subsect:tangNew}

Since on the boundary $\partial\Omega$ of a tangential polygon we have 
$x\cdot{n}=\rho$, equation~(\ref{eqt:bndryQ}) yields,
for tangential polygons,
\begin{equation}
 \QQ_0 = \frac{\rho}{4}\int_{\partial\Omega} 
\left( \frac{\partial u_0}{\partial n}\right)^2 . 
\label{eqt:QQ0t}
\end{equation}
This enables us to rewrite the preceding theorem as follows.

\begin{theorem}{\label{thm:thm2}
For any tangential polygon $\QQ_0$ satisfies the  inequality, quadratic in $Q_0$,
\begin{equation}
\frac{A^2}{L}-\frac{(\Sigma_\infty-\QQ_0)^2}{\Sigma_1}\le \frac{4}{\rho} \QQ_0  .
\label{in:QQ1R1at}
\end{equation}
}
\end{theorem}

In our application we treat~inequality~(\ref{in:QQ1R1at}) as a quadratic inequality in $\QQ_0$
and it is satisfied if 
$$ \QQ_{0-} \le \QQ_0 \le \QQ_{0+} , $$
where, with
$$ \delta = -\frac{\Sigma_1}{L}\left(
 2A L^2\Sigma_\infty -L^3\Sigma_1 -A^4 \right) ,
$$
$$ Q_{0\pm}
= \frac{1}{A}\left(-L\Sigma_1+A\Sigma_\infty \pm\sqrt{\delta} \right) .
$$

\subsection{$\Sigma_\infty$ and $\Sigma_1$}\label{subsect:tangSigma}

For a tangential polygon
\begin{equation} 
u_\infty = c_0 -\frac{1}{4} r^2\qquad
{\rm where \ } r^2 = x^2+ y^2 , \ {\rm and\ \ }
c_0 = \frac{i_2}{4 L} .
\label{eqt:c0def}
\end{equation}
$c_0$ is such that the boundary integral is zero.
This was noted at equations~(2.14), (2.15) of~\cite{KM93}.
We find 
\begin{eqnarray}
\Sigma_\infty
&=& 
\frac{1}{16}\rho i_2 = \frac{A i_2}{8 L},
\label{eqt:tangPSigInf}\\
\Sigma_1
&=& -c_0^2\,  L + \frac{1}{2} c_0 i_2 -\frac{1}{16}  i_4
=\frac{1}{16}\left( \frac{i_2^2}{L}  -  i_4 \right) .
 \label{eqt:tangPSig1}
\end{eqnarray}

We can now readily rewrite our inequality~(\ref{in:QQ1R1at}) in terms of the geometric quantities $i_2$ and $i_4$.
Using  the  expressions~(\ref{eqt:tangPSigInf},\ref{eqt:tangPSig1})
in terms of $i_2$, $i_4$ we find that
the inequality of~(\ref{in:QQ1R1at}) is satisfied iff
the following inequality, quadratic in $\QQ_0$, is satisfied
\begin{equation}
f(\QQ_0)
:= 32 A\QQ_0^2 -
\frac{4\QQ_0}{\rho}\left(2A i_4 -\rho i_2^2 + A i_2 \rho^2\right) +
A\rho\, (A i_4 -\frac{3}{8} i_2^2\rho) 
\le{0} \ .
\label{eqt:fDef}
\end{equation}
Inequalities on the domain functionals $A$, $\rho$, $i_2$ and $i_4$
can be used to establish that both roots of $f(\QQ_0)=0$
are positive.
Denote the smaller root by $\QQ_{0-}$ and the larger by $\QQ_{0+}$.
The inequalities $\QQ_{0-}\le\QQ_0\le\QQ_{0+}$ improve,
for tangential polygons, some well known inequalities such as
the Polya-Szego inequality~(\ref{in:PS}),
$$ \QQ_0
\ge \frac{A^2}{4\ B_\Omega}  ( = {\rm{ \ for\ tangential\ polygons\ }} \frac{A^3}{2L^2}
=\frac{\rho^2 A}{8} ) .
$$

We now comment on the quadratic $f$.
Consider first a disk radius 1
$$ A_\odot=\pi,\ \rho_\odot=1, \ i_{2,\odot}=i_{4,\odot}= 2\pi , $$
$$\leqno{{\rm so}}\qquad\qquad\qquad
f_\odot(\QQ_0)
= 32\pi \QQ_0^2 - 4\QQ_0 (2\pi^2) +\frac{1}{2}\pi^3 = 32\pi(\QQ_{0-}\frac{\pi}{8})^2 .
$$
This agrees with that $\QQ_{0,\odot}=\pi/8$. 
Next consider an equilateral triangle,
$$ A_\Delta=\sqrt{3},\ \rho_\Delta=\frac{1}{\sqrt{3}}, \ i_{2,\Delta}=4,\ i_{4,\Delta}= \frac{16}{5} , $$
$$\leqno{{\rm so}}\qquad\qquad\qquad
f_\Delta(\QQ_0)
= 32\sqrt{3} \QQ_0^2 - 4\QQ_0\sqrt{3}(\frac{12\sqrt{3}}{5}) + \frac{6\sqrt{3}}{5}
= 32\sqrt{3} (\QQ_{0-}\frac{\sqrt{3}}{20})(\QQ_0 - \frac{\sqrt{3}}{4}) ,
$$
which is consistent with $\QQ_{0,\Delta}=\sqrt{3}/20$.

Consider next general tangential polygons.
Define $\QQ_B =\rho^2 A/8$, and recall the
Polya-Szego inequality $\QQ_0\ge\QQ_B$.
We have
$$ f(\QQ_B) = \frac{\rho^2 A}{2} \left(\rho A-\frac{1}{2}\rho i_2\right)^2 . $$
Thus the inequality $\QQ_0\ge\QQ_{0-}$ improves on $\QQ_0\ge\QQ_B$.
Consider next the upper bound
$\QQ_0\le{I_O/4}$ where $I_O$ is the polar moment of inertia about the incentre $O$,
which can also be written $\QQ_0\le\Sigma_\infty$.
Then $I_O=4\Sigma_\infty=\rho i_2/4$ and
$$f(\frac{1}{16}\rho i_2)
= -\frac{1}{2}\left(\frac{i_2}{2}-\rho A\right)\, \left(2A i_4-\rho i_2^2\right) .$$
Since both terms in parentheses are positive, one concludes that 
$\rho i_2/16\le\QQ_{0+}$ so the
bound $\QQ_0\le\QQ_{0+}$ is weaker than the earlier bound. 
Summarizing, we have
$$ Q_B=\frac{1}{8}\rho^2 A \le \QQ_{0-} \le \QQ_0 \le \frac{1}{16}\rho i_2 \le \QQ_{0+} .$$

Using the calculations of $\Sigma_\infty$ and $\Sigma_1$ for isosceles triangles with area $\sqrt{3}$
given in~\iK\ 
we show in Figure~\ref{fig:plIsos}  how the new inequality compares
with earlier inequalities.
As another check we note that perhaps the most studied isosceles triangle other than the equilateral is
the right isosceles triangle.
Let $\alpha$ be the apex angle of the isosceles triangle and
$\sigma=\tan(\alpha/4)$ so $0<\sigma<1$.
For this $\sigma=\sqrt{2}-1\approx{0.414214}$ and for area $\sqrt{3}$
its torsional rigidity is approximately 0.07827 (see~\cite{PoS51}) 
which is, as it must be, larger than $\QQ_{0-}$ which, at this $\sigma$, is 0.076511.
However, at just 2.5\% difference, it is too close to the curve to plot usefully.
\section{Conclusion, and open problems}

There are many questions.

We do not know if inequality~(\ref{eqt:QQ1R1a}) is implied by some other known inequality
for torsional rigidity.

Inequality~(\ref{in:QQ1R1at}) becomes an equality for circular disks and equilateral triangles.
It may be that these are the only shapes which achieve this.

For isosceles triangles with a given area,
as indicated in Figure~\ref{fig:plIsos}, $\QQ_{0-}$ is maximized
(and $\QQ_{0+}$ minimized) by the equilateral triangle.
We believe that this will be true for all triangles.
What happens with quadrilaterals, and more generally $n$-gons, is not known.

It would also be possible to check how well the inequality checks with computed torsional rigidities
for polygons, a few references being~\cite{Ha03,RSNC16,SC65}.
Perhaps checking for the regular polygons where the data is given in Table 1 of~\iK
would be the easiest place to start.

To date the work has been on domains for which $u_\infty$ is a quadratic polynomial.
(\cite{KW20} treats rectangles.)
Outside this class of domains,
there are other domains for which all the functionals 
in~inequality~(\ref{eqt:QQ1R1a}) can be found,
for example, the semi-circle has elementary function solutions for both $u_0$ and $u_\infty$ and
the various functionals can be found (using Maple to sum the series).

\newpage

\section*{Structure of the remainder of this document} 
The preceding part of this document, called Part I from now on, has been
accepted, subject to revision, for publication in {\it IMA Journal of Applied Mathematics}.\\
One referee was very positive and  concluded:\\
``This paper too is very well written, clear, and easy to follow. I strongly recommend publication in IMA Journal of Applied Mathematics"\\
The other referee wrote of the paper:\\
``the work was not put into context and background/other related studies were not discussed .... I think it’d be good to discuss previous work on this topic, so that it is clear what the new contribution of this work is. The derived inequalities will be interesting once a discussion is given (e.g of other inequalities which have appeared in the literature, why are these important etc.). "\\
The same referee also suggested:\\
``computing torsional rigidities for specific cases and comparing to other studies".\\
This supplement is intended to address this, while leaving the journal paper short and focussed on the new research. \\
Numerics for the torsional rigidities of regular polygons and of isosceles triangles 
(in the papers submitted to IMA) are repeated here in Parts I and IIa.
Additional numerics for rhombi are given near the end of Part III.

\newpage
\begin{itemize}
\item The earlier part of this document, Part I, is a preprint form of the IMA paper.
Some items addressing a referee's concerns with the original form or the paper
are in Part Ib.

\item Part II concerns geometric matters relating to tangential polygons.\\
Part IIa is adapted from material used in a different context in the paper~\cite{Ke20i}\\
 Part IIb contains geometric items not in the IMA paper(s).\\
 The first topics are related to $n$-gons including`isoperimetric inequalities',
 when, with the number of vertices $n$ fixed, with some fixed parameter (e.g. area)  
 regular polygons optimize over some other parameter (e.g. minimize the perimeter).\\
 A second topic is inequalities, sometimes not involving $n$ or at least allowing $n$ to range over
 positive integers, `circumgons' and `circum-$n$-gons', i.e. shapes in which part of the 
 boundary is the disk with radius the inradius.
 One of these is the `single-cap', the convex hull of the disk and a single point outside it:
 a circum-1-gon.
 Another is the `symmetric double-cap':  a circum-2-gon. 
 We will see these in connection with Blaschke-Santalo diagrams.
 
 \item In Part III we return to considerations of torsional rigidity.
 The emphasis is on bounds for convex domains,
 in particular convex polygons especially tangential polygons.
 Some sections are devoted to triangles, especially isosceles,
 and tangential quadrilaterals, especially rhombi.
 
 \item The remaining parts are only slightly connected to the torsion problem.\\
 Part IV concerns replacing the Dirichlet b.c. with a Robin b.c..\\
 Part V notes some other pde problems where tangential polygons are mentioned.
 
\end{itemize}

The treatment is very uneven.
I have not checked the more advanced recent pde papers.
Some of the geometry papers cited in Part II are very elementary.
The suggestion (by Buttazzo) that I look at Blaschke-Santalo diagrams has led to items
at present poorly integrated with the study of my bound $\QQ_{0-}$
(with just Part III~\S\ref{sec:BSQ} indicating one direction).
I hope a later version of this supplement will correct some of these defects.


\newpage
\begin{center}
{\large{{\textsc{ Part Ib
}}}}
\end{center}

\section*{Numerics for the isosceles triangle}

Numerics for the isosceles triangle were described earlier, but here is some amplification.

\begin{figure}[ht]
\centerline{\includegraphics[height=10cm,width=14cm]{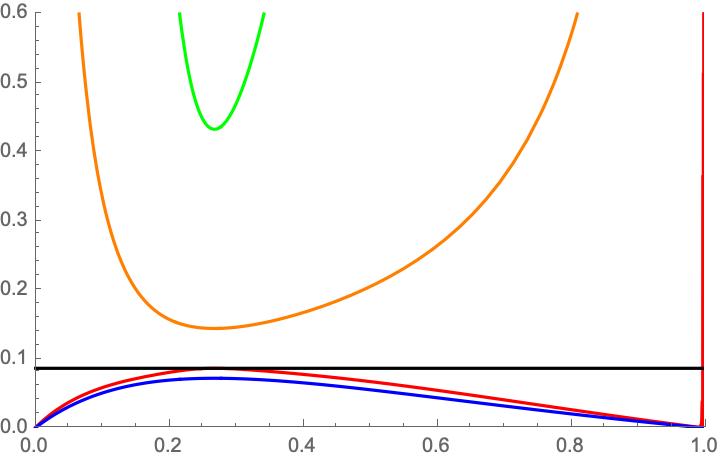}}
\caption{For an isosceles triangle with area$\sqrt{3}$.
$\sigma$ is tan of 1/4 of the apex angle.
Blue is $\QQ_B$, red is the new lower bound $\QQ_{0-}$,
black is $\QQ_\Delta$, orange is $\rho i_2/16$, green is $\QQ_{0+}$.
 }
\label{fig:plIsos}
\end{figure}

\begin{figure}[ht]
\centerline{\includegraphics[height=10cm,width=14cm]{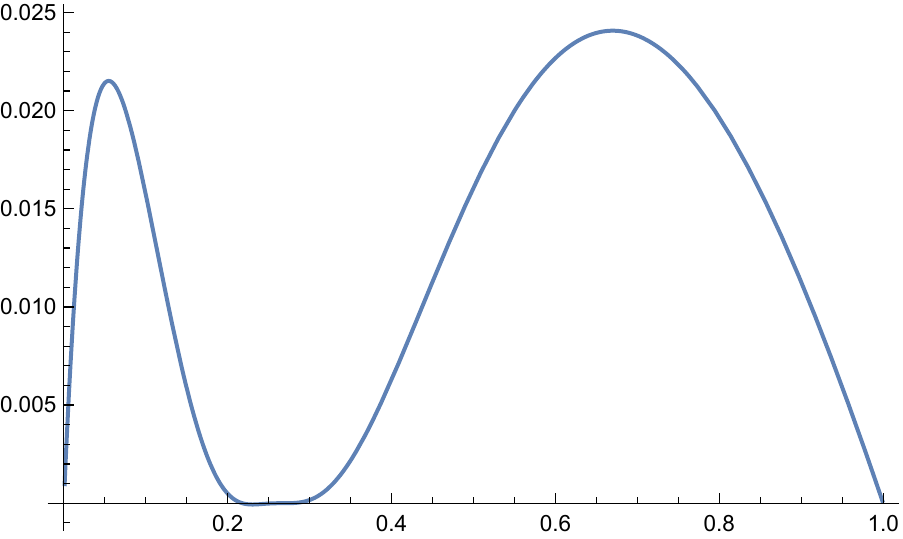}}
\caption{A plot, against $\sigma$, of the difference
$\QQ_{0-}-\QQ_{\rm Sol}$, the latter term defined in
equation(\ref{eqt:Sol92}).
 }
\label{fig:IIasolQ0m}
\end{figure}

Another lower bound on $\QQ_0$, as in~\cite{Sol92}, is that, amongst triangles with
a given inradius, the equilateral triangle has the smallest $\QQ_0$.
Thus
\begin{equation}
\QQ_0\ge \QQ_{\rm Sol} = \frac{9\sqrt{3}}{20} \rho^4 . 
\label{eqt:Sol92}
\end{equation}
For isosceles triangles this lower bound improves on our $\QQ_{0-}$
when the apex angle is slightly less than $\pi/3$.
See Figure~\ref{fig:IIasolQ0m}.

\clearpage

\section*{Speculation on equality in inequality~(\ref{in:QQ1R1at}).}

Inequality~(\ref{in:QQ1R1at}) becomes an equality for circular disks and equilateral triangles.
I don't have a proof, but it might be that one only gets equality
for these shapes.

There is just one inequality used in the proof of the theorems.
It is a Cauchy-Schwarz Inequality and it will be an equality
iff
$$ u_\infty-c_* = {\rm const}\, \frac{\partial u_0}{\partial n} , $$
all the way around the boundary.
The lhs is a quadratic function of $(x,y)$.
Consider next a genuine polygon, a tangential $n$-gon.
Without loss of generality consider a side parallel to the $x$-axis
as the interval $-1\le{x}\le{1}$, and $y=-h$,
with the $n$-gon being in the half-plane $\{ y>-h\}$.
The only possibility is
$$ \frac{\partial u_0}{\partial y} = -\frac{1}{4} (1-x^2) , 
\qquad{\rm on \ } y=-h,\ \ -1<x<1 ,$$
as the gradient of $u_0$ is zero at the corners.
This has implications for the values of $c_0$ in $u_\infty$
and $c_*$.
Further information might be obtained by considering sides adjoining
the $y=-h$ side,
with angles $\alpha_-$ at $x=-1$ and $\alpha_+$ at $x=1$.
One might get other restrictions on the $c_0$ and $c_*$ values,
and on the $\alpha\pm$.

Even if the disk and equilateral triangle are the only
bounded domains for which we get equality, it might be difficult
to show this.
If one considers an infinite wedge as `tangential polygon'
this might be a counterexample.
Consider a wedge, apex at the origin and symmetric about the
$x$-axis, $\theta=0$.
Suppose that the wedge has angle $\alpha$. 
Then
$$u_{0,{\rm wedge}}
= -\frac{r^2}{4}\left( 1- \frac{\cos(2\theta)}{\cos(\alpha)}\right) 
+ {\rm const} r^{\pi/\alpha}\cos(\pi\frac{\theta}{\alpha}) ,
$$
solves the torsion equation, has zero Dirichlet boundary data.
Here
$$ \frac{\partial u_0}{\partial n} = 
\frac{1}{r} \frac{\partial u_0}{\partial\theta} .
$$
It might be that we can set the const to 0 to see a counter example.
(The solution in a sector is given in~\cite{KCT97}.)

\newpage
\begin{center}
{\large{{\textsc{ Part IIa:
$\Sigma_\infty$ and $\Sigma_1$
for tangential polygons
}}}}
\end{center}

\section*{Abstract for Part IIa}
Items useful in calculation for tangential polygons in general,
and for particular cases, are here extracted from~\cite{Ke20i}.

\section{Outline of Part IIa}\label{sec:OutlineIIa}

In this part
the functionals $\Sigma_\infty$ and $\Sigma_1$ are calculated for various tangential polygons.
The starting point is
equation~(2.14), (2.15) from the 1993 paper~\cite{KM93}  --
repeated at appropriate times in this document -- which gives 
$u_\infty=c_0-\frac{1}{4}(x^2+y^2)$ for tangential polygons,
with $c_0$ such that the boundary integral of $u_\infty$ is 0.
\begin{itemize}
\item In \S\ref{sec:Tang}
we find that starting from $u_\infty$ 
the functionals $\Sigma_\infty$ and $\Sigma_1$ can be expressed in terms of various moments of inertia.

\item In \S\ref{sec:xtremeShapes} we treat 
the circular disk and the equilateral triangle.

\item
In \S\ref{sec:regn} we study regular $n$-gons.

\item Any triangle is a tangential polygon.\\
In \S\ref{sec:isos} we find the functionals for any isosceles triangle.

\item 
In~\S\ref{sec:tangQuad}.
we note papers relevant to work on tangential quadrilaterals.

\end{itemize}

We will use `tangential polygon' as defined before.
\begin{itemize}
\item
Genuine $n$-sided polygons for which every side is tangent to
the incircle
(and hence for which the intersection of the boundary with the incircle
is $n$ points, the points of tangency) will
be called {\it tangential $n$-gons}.
\item
When tangential polygon is the convex hull of 
$n$ points outside the incircle and the incircle we will
call it a {\it circum-$n$-gon}.
(This is a slightly different terminology than in~\cite{AM04}.
\end{itemize}

Tangential $n$-gons are particular cases of circum-$n$-gons.
The union over all $n$ from 0 to $\infty$ of  circum-$n$-gons
gives all tangential polygons (the $n=0$ case being considered
as the disk).
A circum-1-gon is also called a 1-cap,
a circum-2-gon is also called a 2-cap:
we will see these in Part IIb.

We use established terminology with a circum-4-gon being called
a tangential quadrilateral,
a circum-6-gon a tangential hexagon, etc.

\section{Tangential polygons}\label{sec:Tang}

\subsection{Geometric results}\label{subsec:tangGeom}

Here we continue from the basic geometric definitions and results 
given in Part IIa~\S\ref{sec:OutlineIIa}.

There are various well-known or elementary observations:
\begin{itemize}

\item Of the polygons with
fixed perimeter and angles, the tangential polygon has the greatest area.

\item Given two (convex) tangential polygons with the same incircle
their intersection is also a (convex) tangential polygon with the same incircle.

\end{itemize}
Papers concerning tangential polygons include~\cite{AM04,RP01}
(and many more particularly concerning tangential $n$-gons
for $3\le{n}\le{6}$ will be given at appropriate places later in
this document).
There is a literature on (convex) tangential polygons, one fact being that,
considering the polygons as linkages touching the incircle,
if the number of sides is odd
the polygon is rigid  but not if the number of sides is even.
Entertaining as such facts are, they do not appear to be relevant to our pde exercise.

We will need boundary moments
$i_{2k}$ -- moments about the incentre -- defined by
$$ i_{2k} = \int_{\partial\Omega} (x^2+y^2)^k ,
$$
and the polar area moment about the incentre
$$ I_{2k} = \int_\Omega (x^2+y^2)^2 . $$

(Caution. There are many results concerning moments about the centroid.
For example, as in~\cite{PoS51}, the moments about the centroid are minimized over $n$-gons with a given fixed area
by the regular $n$-gon.
We haven't checked in general if it is the case that the  moments about the incentre are minimized over $n$-gons 
by the regular $n$-gon, but when $n=3$ and one minimizes over isosceles triangles, the equilateral triangle is 
the minimizer.)

We remark that the Cauchy-Schwarz inequality for the integrals implies that
$i_4\ge{i_2^2/L}$ (and $I_4\ge{I_2^2/A}$).

\medskip

\subsubsection{$I_{2k}=\rho\,i_{2k}/(2k+2)$.}

A result which we first noticed in two special cases is the following.
\smallskip

{\par\noindent}{\bf Result.}
{\it For any tangential polygon
$A=\rho\,L/2$, $I_2=\rho\,i_2/4$ and, more generally,  $I_{2k}=\rho\,i_{2k}/(2k+2)$.}

{\par\noindent}{\it Proof.}
The Divergence Theorem can be used to give a boundary integral which equals $I_2$.
Consider polar coordinates, radial coordinate $r$, and the Laplacian of a function of $r$:
$$ \Delta u = \frac{1}{r}\frac{\partial}{\partial r} r \frac{\partial u }{\partial r} . $$
Also, for use in the following, for any tangential polygon $x.n=\rho$, here $x$ being the vector to a point on $\partial\Omega$.
To avoid using $x$ as a vector, equally $ r e_r.n=\rho$, with $e_r$ the unit vector in the $r$-direction.

With $u=r^2$ so $\Delta r^2 = 4$ the Divergence Theorem gives
$$ 2\rho L = \int_{\partial\Omega} 2 r e_r.n = \int_\Omega 4 = 4 A . $$
Similarly, with $u=r^4$ so $\Delta r^4 = 16 r^2$ the Divergence Theorem gives
$$ 4\rho i_2 = \int_{\partial\Omega}( \nabla r^4).n =\int_\Omega 16\, r^2 = 16 I_2 . $$

{\par\noindent}{\it Alternative Proof.}
Start with a tangential polygon $T(1)$ origin at the incentre, with radius 1.
Define $T(\rho)=\rho\, T(1)$ with second area moment $I_2(\rho)$ and
second boundary moment $i_2(\rho)$.
Along with $T(\rho)$ consider a similar scaled polygon $T(\rho+\Delta\rho)$
with $\Delta\rho$ small.
Then
$$I_2(\rho+\Delta\rho)-I_2(\rho) \sim \Delta\rho\, i_2(\rho)+O((\Delta\rho)^2) . $$
But, dimensionally,  $i_2(\rho)=\rho^3\, i_2(1)$ and $I_2(\rho)=\rho^3\, I_2(1)$.
Taking limits in the displayed equation gives the central part of
$$ 4\rho^3 I_2(1) = \frac{d I_2(\rho)}{d\rho}= i_2(\rho) = \rho^3 i_2(1). $$
This gives $i_2(1)=4I_2(1)$ and hence, more generally, $I_2(\rho)=\rho\,i_2(\rho)/4$
as asserted.

The same methods give $I_{2k}=\rho\,i_{2k}/(2k+2)$, and the case $k=0$ is $A=\rho\,L/2$.

\subsubsection{Methods to calculate $i_{2k}$.}

Denote with a bar just contributions from a vertical side, at $x=\rho$, extending from
$y=-\eta_-$ to $y=\eta_+$.
Then
$${\bar i}_{2k} = \int_{-\eta_-}^{\eta_+} (\rho^2+y^2)^k\, dy .$$
Hence
\begin{eqnarray*}
{\bar i}_0
&=& \eta_+ + \eta_- ,\\
{\bar i}_2
&=& \rho^2 {\bar i}_0 +\frac{1}{3} (\eta_+^3 + \eta_-^3) ,\\
{\bar i}_4
&=& -\rho^4 {\bar i}_0 + 2\, \rho^2 {\bar i}_2 +\frac{1}{5} (\eta_+^5 + \eta_-^5) .
\end{eqnarray*}
Assuming the tangential polygon has $n$ sides, this gives
\begin{eqnarray}
{  i}_0=L
&=& \sum_{k=1}^n (\eta _{k+} + \eta _{k-} ),
\label{eq:i0eta}\\
{  i}_2
&=& \rho^2 {  i}_0 +\frac{1}{3} \sum_{k=1}^n (\eta _{k+}^3 + \eta _{k-}^3) ,
\label{eq:i2eta}\\
{  i}_4
&=& -\rho^4 {  i}_0 + 2\, \rho^2 {  i}_2 +\frac{1}{5} \sum_{k=1}^n (\eta _{k+}^5 + \eta _{k-}^5)  .
\label{eq:i4eta} 
\end{eqnarray}
There is the geometrically evident fact that $\eta _{k+}=\eta _{k+1,-}$ as one traverses from one
of the polygon's sides to the next.

\medskip

The easiest case to consider is the regular $n$-gon, side $s_n=2\rho\tau_n$ where $\tau_n=\tan(\pi/n)$.
Then
\begin{eqnarray}
{  i}_0=L
&=& n\, s_n = 2n \rho\tau_n ,
\label{eq:i0reg} \\
{  i}_2
&=& \rho^2 {  i}_0 +\frac{n}{12} s_n^3  = n\, s_n \left(\rho^2 + \frac{1}{12} s_n^2 \right)
=\frac{2}{3} n \rho^3 \tau_n (3 + \tau_n^2), 
\label{eq:i2reg} \\
{  i}_4
&=& -\rho^4 {  i}_0 + 2\, \rho^2 {  i}_2 +\frac{n}{80} s_n^5
= n\, s_n \left( \rho^4 +\frac{\rho^2 s_n^2}{6} +\frac{s_n^4}{80}\right) , \\
&=& \frac{2}{15} n \rho^5 \tau_n\left( 15 + 10\, \tau_n^2 + 3\, \tau_n^4 \right)
 . \label{eq:i4reg}
\end{eqnarray}
These had been calculated independently from first principles,
as reported in~\S\ref{subsec:Generaln},
before our observations concerning general tangential polygons.

We now consider how $\eta_+$ and $\eta_-$ might be found in terms of vertex angles of the
tangential polygon.
Recall that in any tangential polygon the angle bisector at any vertex
passes through the incentre.
Suppose side $k$ lies between vertices $k$ and $k+1$.
Let the angle at vertex $k$ be $\alpha_k$.
(If there are $n$ sides the sum over all the $\alpha_k$ is $(n-2)\pi$.)
Then, with
\begin{equation}
T_k  = \frac{1}{\tan\frac{\alpha_k}{2}}=\tan(\frac{\pi-\alpha_k}{2}),
\label{eq:Tkdef}
\end{equation}
$$ \eta_{k-}= {\rho} T_k,\qquad
 \eta_{k+}= \rho T_{k+1}.
 $$
Summing over $k$
\begin{equation}
A= \rho^2 {\sum\atop{k}} T_k, \qquad
L= 2\rho  {\sum\atop{k}} T_k .
\label{eq:i0Tgen}
\end{equation}
 (This checks with our previous $A=\rho L/2$.)
 Another well-known identity is
 $$ \frac{L^2}{4A} =  {\sum\atop{k}} T_k. $$
 Similarly
\begin{eqnarray}
{  i}_2
&=& \rho^2 {  i}_0 +\frac{2\rho^3}{3} \sum_{k=1}^n T_k^3
= 2\rho^3\, \sum_{k=1}^n\left(T_k+\frac{ T_k^3}{3}\right),
\label{eq:i2Tgen}\\
{  i}_4
&=& -\rho^4 {  i}_0 + 2\, \rho^2 {  i}_2 +\frac{2\rho^5}{5} \sum_{k=1}^n  T_k^5 
\nonumber\\
&=&  2\rho^5\, \sum_{k=1}^n\left(T_k+\frac{ 2T_k^3}{3}+\frac{ T_k^5}{5}\right).
\label{eq:i4Tgen}
\end{eqnarray}

\subsection{$\Sigma_\infty$ and $\Sigma_1$}\label{subsec:tangSigma}

For a tangential polygon
$$ u_\infty = c_0 -\frac{1}{4} r^2\qquad
{\rm where \ } r^2 = x^2+ y^2 , $$
and $c_0$ is such that the boundary integral is zero:
\begin{equation}
c_0 = \frac{i_2}{4 L} .
\label{eq:c0def}
\end{equation}
where $i_{2k}$ are the boundary moments defined previously.
This was noted at equations~(2.14), (2.15) of~\cite{KM93}.
Hence
\begin{eqnarray}
\Sigma_\infty
&=& 
\frac{A\, i_2}{4 L}- \frac{1}{4} I_2
=\frac{1}{16}\rho i_2 = \frac{A i_2}{8 L},
\label{eq:tangPSigInf}\\
\Sigma_1
&=& -c_0^2\,  L + \frac{1}{2} c_0 i_2 -\frac{1}{16}  i_4
=\frac{1}{16}\left( \frac{i_2^2}{L}  -  i_4 \right),
 \label{eq:tangPSig1}
\end{eqnarray}
and the notation, as before, has $I$ for area moments, and $i$ for boundary moments.

\section{Equilateral triangle and disk}\label{sec:xtremeShapes}

The results for these domains are well known: see e.g.~\cite{MK94}.
For the unit disk
$$ A=\pi,\ L_\odot=2\pi,\ \QQ_{\odot,0}=\Sigma_{\odot,\infty}=\frac{\pi}{8},\
\Sigma_{\odot,1}=0, \
i_{\odot,2}=i_{\odot,4}=2\pi,\ I_{\odot,2}=\frac{\pi}{2}.$$

For an equilateral triangle with side $2a=s_3$,
$$A_{\Delta}= a^2\,\sqrt{3},\ L_{\Delta}=6a,\
\QQ_{\Delta,0} = \frac{\sqrt{3} a^4}{20} = \frac{A_{\Delta}^2}{20\sqrt{3}},\  $$
$$
\Sigma_{\Delta,\infty} = \frac{a^4}{4\sqrt{3}} =  \frac{A_{\Delta}^2}{12\sqrt{3}},\
 \Sigma_{\Delta,1}=  \frac{A_{\Delta}^{5/2}}{90\, \ 3^{1/4}} .$$
$$i_{\Delta,2}=  \frac{4}{3} \, 3^{1/4}\, A^{3/2},\
i_{\Delta,4}=  \frac{16}{45} \, 3^{3/4}\, A^{5/2} ,\
I_{\Delta,2}= \frac{\sqrt{3}}{9}\, A^2.$$
For triangles and disks, both with area $\pi$,
the St Venant isoperimetric inequality is consistent with
$$0.3927= \QQ_{\odot,0} =\frac{\pi}{8}>\QQ_{\Delta,0}=\frac{\pi^2}{20\sqrt{3}}=0.28491\ .$$
The inequality for the $\Sigma_\infty$ is
$$0.3927= \Sigma_{\odot,\infty} =\frac{\pi}{8}<
\Sigma_{\Delta,\infty} =\frac{\pi^2}{12\sqrt{3}}=0.47485\ . $$
\section{Regular $n$-gons}\label{sec:regn}

\subsection{General $n$}\label{subsec:Generaln}

We denote the inradius by $\rho$, the side by $s$, area by $A$,
perimeter by $L$ and the angle at the centre subtended by a single side by $\gamma$.
For the (regular) $n$-gon, simple geometry gives $\gamma_n=2\pi/n$ so
$$ \frac{s_n}{2\rho_n}=\tau_n
\qquad{\rm where\ } \tau_n=\tan(\frac{\gamma_n}{2}) =\tan(\frac{\pi}{n}), $$
and
$$A_n = \frac{n s_n\rho_n}{2}= \frac{n s_n^2}{4\tau_n}
= n\rho_n^2 \tau_n .$$
We will wish to specify the area (to be $\pi$), so we note
$$ s_n^2 =\frac{4 A_n \tau_n}{n} \qquad{\rm and\ }
L_n = n s_n  =2\sqrt{n A_n \tau_n}.$$
The inradius $\rho_n$ occurs in some formulae, so we note
$$\rho_n^2 = \frac{A_n}{n \tau_n}  . $$

We are not aware of any simple formula for the torsional rigidity $\QQ_{0}(n)$
for a regular $n$-gon, but there have been many numerical studies
(and theoretical studies starting from Schwarz-Christoffel conformal mapping).
Some numerical results will be given for particular instances later.

Formulae for $i_0=L$, $i_2$ and $i_4$ have been presented earlier at
equations~(\ref{eq:i0reg}), (\ref{eq:i2reg}) and~(\ref{eq:i4reg}).
In view of their central role and the detailed algebraic manipulations in their derivation
we record here an independent derivation.
Using polar coordinates, and considering the side with $x=\rho_n=r\cos(\theta)$ for which the polar angle at the centre lies between $\theta=-\pi/n$ and $\theta=\pi/n$,
$$ i_2(n) = n \int_{-\pi/n}^{\pi/n}
\left(\frac{\rho_n}{\cos(\theta)}\right)^2 \, \frac{dy}{d\theta}\, d\theta ,$$
in which $y=\rho_n\tan(\theta)$.
Thus
\begin{eqnarray*}
 i_2(n)
 &=& n\rho_n^3  \int_{-\pi/n}^{\pi/n} \frac{1}{\cos(\theta)^4}\, d\theta ,\\
 &=& \frac{2\,n\rho_n^3 \, \left(1+2\cos(\frac{\pi}{n})^2\right)\,
\sin(\frac{\pi}{n})}{3\,\cos(\frac{\pi}{n})^3} , \\
&=& \frac{2}{3} n \rho_n^3 \tau_n (3+\tau_n^2) ,\\
&=& \frac{2}{3} \sqrt{\frac{A^3}{n \tau_n}} (3+\tau_n^2) .
\end{eqnarray*}
 Similarly
 \begin{eqnarray*}
 i_4(n)
 &=& n\rho_n^5  \int_{-\pi/n}^{\pi/n} \frac{1}{\cos(\theta)^6}\, d\theta ,\\
 &=&  \frac{2\,n\rho_n^5 \, \left(3+4\cos(\frac{\pi}{n})^2+8\cos(\frac{\pi}{n})^4\right)\,
\sin(\frac{\pi}{n})}{15\,\cos(\frac{\pi}{n})^5} ,\\
&=& \frac{2}{15}  n \rho_n^5 \tau_n (15+10\tau_n^2+3\tau_n^4) .
\end{eqnarray*}

The area moment of inertia $I_2(n)$ is similarly calculated from that of
the isosceles triangle with apex at the origin.
Before noting the general relation at equation~(\ref{eq:tangPSigInf})
a short calculation -- just for regular polygons -- established that, for a regular polygon,
\begin{equation}
I_2(n) =\frac{\rho_n}{4}\, i_2(n) .
\label{eq:i2I2}
\end{equation}

Using equations~(\ref{eq:i0reg}), (\ref{eq:i2reg}) and~(\ref{eq:i4reg}),
equation~(\ref{eq:tangPSigInf}) becomes
\begin{equation}
\Sigma_\infty(n)
=  i_2(n)\, \left( \frac{A}{4 L_n}- \frac{\rho_n}{16} \right)
= \frac{ A^2\, (3+\tau_n^2)}{24 n\tau_n} .
\label{eq:tangPSigInf2}
\end{equation}
(From~(\ref{eq:tangPSigInf2}), $\Sigma_\infty>\Sigma_{\odot,\infty}=\QQ_{\odot,0}=\pi/8$
when $A=\pi$.)
Also
\begin{equation}
\Sigma_1(n)
= - \frac{1}{16 L}\left( Li_4 -i_2^2\right)
= -\frac{1}{90}\sqrt{\frac{ A^5\, \tau_n^5}{n^3}} .
\label{eq:tangPSig12}
\end{equation}

The original motivation for assembling this particular list of quantities
was that they were needed for a bound associated with slip flow down a
pipe with cross-section $\Omega$: see~\cite{Ke20i}.
Numeric values for the torsional rigidity are available in several references.
The entries for $J/A^2$, where $J=4Q_0$, in the following table are taken from~\cite{Ha03}.
See also~\cite{RSNC16} and, for $n=3,\ 4$ and $6$, also~\cite{PoS51}.
In Table~\ref{tbl:tbl1} we take our polygons to have area $\pi$.
In Table~\ref{tbl:tbl2} we take our polygons to have circumradius 1.

\begin{table}[h]
\begin{center}
\begin{tabular}{|| c | c | c | c| c | c ||}
\hline
$n$& $4\QQ_0/A^2$& $\QQ_0$ & $L$ & $\Sigma_\infty$& $-\Sigma_1$ \\
3& 0.11547&  0.28492& 8.0806&  0.47485&  0.14769 \\ 
4& 0.14058& 0.34687&  7.0898&  0.41123&  0.024296 \\
5& 0.14943& 0.36870&  6.7565&  0.39936&  0.007822 \\
6& 0.15340& 0.37850&  6.5978&  0.39571&  0.003349 \\
7& 0.15546& 0.38358&  6.5086&  0.39426&  0.001689\\
8& 0.15664& 0.38649&  6.4530&  0.39359&  0.0009485 \\
9& 0.15736& 0.38827&  6.4159&  0.39325&  0.0005754 \\
10&0.15783& 0.38943&  6.3899&  0.39306&  0.0003699 \\
11&0.15815& 0.39022&  6.3709&  0.39294&  0.0002489 \\
12&0.15837& 0.39076&  6.3566&  0.39287&  0.0001738  \\
$\infty$&  0.15915 &0.3927 & 6.2832& 0.3927&0 \\
\hline
\end{tabular}
\caption{Regular polygons with area $\pi$}
\label{tbl:tbl1}
\end{center}
\end{table}

\begin{table}[h]
\begin{center}
\resizebox{\textwidth}{!}{
\begin{tabular}{|| c | c | c | c| c | c | c  | c ||}
\hline
$n$& $4\QQ_0/A^2$&$A_n$& $\QQ_0$ & $L_n$ & $\Sigma_\infty$& $-\Sigma_1$& $\rho_n$ \\
\hline
3& $\sqrt{3}/15$&   $3\sqrt{3}/4$ &$9\sqrt{3}/320$& $3\sqrt{3}$ &  $3\sqrt{3}/64$&  $3\sqrt{3}/320$& $1/2$\\ 
 & 0.11547& 1.2990 & 0.0487 & 5.1962 & 0.0812&  0.0162 & 0.5\\
\hline
4& 0.14058&  $2$ &             0.14058&  $4\sqrt{2}$ &  $1/6$&  $\sqrt{2}/180$& $1/\sqrt{2}$\\
  & & 2 & & 5.6569 & 0.1667 & 0.007857& 0.707107 \\
\hline
6&  0.15340&   $3\sqrt{3}/2$ &0.2589 &$6$ & $5\sqrt{3}/32$& $1/480$&$\sqrt{3}/2$ \\
 & &2.5981& & 6 & 0.270633& 0.0020833& 0.866025 \\
\hline
$n$& & $\frac{n}{2}\sin(\frac{2\pi}{n})$ & &$2n\sin(\pi/n)$  & (\ref{eq:tangPSigInf2}) &(\ref{eq:tangPSig12}) &$\cos(\pi/n)$\\
\hline
$\infty$&  $1/(2\pi)$ & $\pi$ &$\pi/8$ & $2\pi$& $\pi/8$&0 & 1\\
 &  0.180043& & &  &  & 0& 1\\
\hline
\end{tabular}
}
\caption{Regular polygons with circumradius 1}
\label{tbl:tbl2}
\end{center}
\end{table}

\clearpage

 \subsection{The square, $n=4$}\label{subsec:square}

For a square with side $s_4$,
$$A_\square= s_4^2, \ L_\square= 4s_4,\ I_{\square,2}=\frac{A^2}{6},\ 
\rho_4=\frac{1}{2}s_4, \ i_{\square,2}=\frac{4}{3} s_4^3 .$$
Hence, consistent with equation~(\ref{eq:tangPSigInf}), we have
$$\Sigma_{\square,\infty} = \frac{s_4^4}{24} =\frac{A_\square^2}{24} .$$
For a square and a disk each of area $\pi$
$$0.3927= \Sigma_{\odot,\infty} =\frac{\pi}{8}<
\Sigma_{\square,\infty} =\frac{\pi^2}{24}= \ 0.41123. $$

Calculation for a rectangle reported in~\cite{KW20} gives (with $a=b$)
$$ \Sigma_{\square,1} = -\frac{s_4^5}{720}  = - \frac{A_\square^{5/2}}{720}. $$

\section{Triangles, especially isosceles}\label{sec:isos}

\subsection{Geometric preliminaries and checks}\label{subsec:isosGeom}

\subsubsection{$A$, $L$, $\Sigma_\infty$ and $\Sigma_1$ for isosceles triangles}

Consider the isosceles triangle whose incentre is at the origin, whose base is length $2a$
and whose apex angle is $\alpha$.
Denote by $\rho$ the inradius of the triangle.
The triangle's vertices are at
$(a,-\rho)$, $(0,h-\rho)$ and $(-a,-\rho)$.
The case where $h=a\sqrt{3}$, $\rho=h/3$ corresponds to an equilateral triangle.

The area of the triangle is denoted by $A$.
The vertex angle is $\alpha$.
Change the notation, what was formerly denoted $T_k$, we will denote
$T_A$, $T_B$ and $T_C$.
It is convenient to define
$$ \sigma = \tan(\frac{\alpha}{4}) . $$
We have
$$T_A=\tan(\frac{\pi-\alpha}{2})=\frac{1-\sigma^2}{2\sigma},\
T_B=T_C=\tan(\frac{\pi+\alpha}{4})= \frac{1+\sigma}{1-\sigma} . $$
We remark though that an alternative parameter, alternative to $\sigma$, is $a$,
$2a$ being the length of the base.
Then
\begin{eqnarray*}
\frac{h}{a}
&=&\cot(\frac{\alpha}{2}) = \frac{1-\sigma^2}{2\sigma} ,\\
\frac{\rho}{a}
&=&\tan(\frac{\pi-\alpha}{4}) = \frac{1-\sigma}{1+\sigma} ,
\qquad \rho=\frac{2A}{L} ,
\end{eqnarray*}
and, using equation~(\ref{eq:i0Tgen}) or directly,
\begin{eqnarray*}
A
&=& a\,h = a^2\,  \frac{1-\sigma^2}{2\sigma}  ,\\
L
&=& 2\left(a+\sqrt{a^2+h^2}\right) = a\frac{(1+\sigma)^2}{\sigma} ,\\
B
&=& \frac{L}{\rho}=\frac{L^2}{2A} =\frac{(1+\sigma)^3}{\sigma(1-\sigma)}.
\end{eqnarray*}
There are several identities relating these geometric quantities.
A relation between the area, perimeter, and base is
\begin{equation}
4 L a^3 - L^2 a^2 + 4 A^2 = 0 .
\label{eq:AaL}
\end{equation}
This is consistent with the geometrically obvious
$$ L = 2 \left(a+\sqrt{a^2+(\frac{A}{a})^2 } \right) ,
$$
which can be regarded as an example of using $a$ as a parameter.
($L(a)$ is a positive convex function with, when $A=\sqrt{3}$ a minimum of 6 when $a=1$.)

Less immediately useful in this paper is the circumradius $R_V$.
The radius of the circumcircle 
(circle through the 3 vertice) is:
$$ R_V=  {\frac {a\,(a^{2}+(\frac{A}{a})^2)}{2A}}. $$
The centre of the circle lies on the symmetry axis of the triangle, this distance below the apex.

An isosceles triangle with given area $A$ is uniquely determined if $\sigma$ is given between 0 and 1
or if $a$ is given between 0 and infinity, in the latter case with
$$\sigma = -\frac{A}{a^2}+\sqrt{\frac{A^2}{a^4} + 1} . $$
There are various formulae for the inradius $\rho$ including
\begin{eqnarray*}
\frac{\rho}{a}
&=& \frac{ \sqrt{a^2+h^2}- a}{h} ,\\
&=& \frac{L-4a}{2h}.
\end{eqnarray*}

{\bf Result.} {\it  At fixed area, $B$ and $L$ are positive convex functions of $\sigma$ for $0<\sigma<1$,
each having its minimum at
$\sigma=2-\sqrt{3}$.
Over isosceles triangles with given area, the equilateral triangle minimizes each of $B$ and $L$.\\
At fixed area, $\rho$ is a positive concave function of $\sigma$ for $0<\sigma<1$,
having its maximum at
$\sigma=2-\sqrt{3}$.
Over isosceles triangles with given area, the equilateral triangle maximizes $\rho$.}

Integration gives, for the area moment,
$$ I_2(A,\sigma)
= \frac{A^2}{12}\, \frac{(1-\sigma)^6+12\sigma^2(1-\sigma)^2+16\sigma^3}{\sigma(1-\sigma)(1+\sigma)^3} .
$$
{\bf Result.} {\it  $I_2(A,\sigma)$ is a convex function of $\sigma$ for $0<\sigma<1$ with its minimum at
$\sigma=2-\sqrt{3}$.
Over isosceles triangles with given area, the equilateral triangle minimizes the area moment of inertia
about the incentre.}

There are various identities for $I_2$, e.g.
\begin{equation}
4 I_2-\frac{1}{6} A L^2 + \frac{8 A^3}{L a} - \frac{16 A^3}{L^2} +2 A a^2 = 0 .
\label{eq:Aai2}
\end{equation}
Eliminating $L$ gives
$$  36 a^4 A I_2^2 -12 a^2 (12 a^8 +11 a^4 A^2 +A^4) I_2 +
A (24 a^{12}+33 A^2 a^8 + 6A^4 a^4 +A^6) = 0 .
$$
{\bf ToDo.} The discriminant of the quadratic for $I_2$ is nice, and the 2 solutions for $I_2$ are reasonably simple
(and possibly tidier than the expressions in $\sigma$).
However surely only one is relevant, which one?\\
Might $i_2=I_2/(4\rho)$ (with $\rho=a^2(L-4a)/(2A)$) be neater than $I_2$.
See if there is an equation like~(\ref{eq:AaL}) involving just $i_2$, $A$ and $a$.
Equation~(\ref{eq:AaL}) involves  $i_0$, $A$ and $a$.
If there is, it might be that $i_2$ considered as a function of $a$ might be tidier than it is
as a function of $\sigma$.

\bigskip

For the boundary moment we split the integral into the part over the base, and over another of the sides:
$$ i_{2k} = i_{2k}(base) + 2\,i_{2k}(side) .$$
Once again, as in equation~(\ref{eq:i2I2}), we find, at $k=1$,
$$
I_2 =\frac{\rho}{4}\, i_2 .
$$
{\bf Result.} {\it
At fixed area, $i_2$ is a  positive convex functions of $\sigma$ for $0<\sigma<1$,
with its minimum at
$\sigma=2-\sqrt{3}$.
Over isosceles triangles with given area, the equilateral triangle minimizes $i_2$
the the boundary moment about the incentre.}

Equation~(\ref{eq:tangPSigInf}) becomes,
with $i_2$ calculated either directly or from~(\ref{eq:i2Tgen}),
\begin{eqnarray}
\Sigma_\infty(A,\sigma)
&=&  i_2\, \left( \frac{A}{4 L}- \frac{\rho}{16} \right)=\frac{1}{16}\rho i_2 = \frac{1}{4} I_2 ,
\nonumber\\
&=& \frac{ A^2}{48} \, \frac{((1-\sigma)^6+12\sigma^2(1-\sigma)^2+16\sigma^3)}{\sigma(1-\sigma)(1+\sigma)^3} .
\label{eq:isosPSigInf2}
\end{eqnarray}

Restating the preceding Result for $I_2$:\\
{\bf Result.} {\it  $\Sigma_\infty(A,\sigma)$ is a convex function of $\sigma$ for $0<\sigma<1$ with its minimum at
$\sigma=2-\sqrt{3}$.
Over isosceles triangles with given area, the equilateral triangle minimizes $\Sigma_\infty$.}\\

The formulae for $i_4$ and $\Sigma_1$ are more elaborate.
See~(\ref{eq:i4Tgen}).
Define
$$ p_1(\sigma)
=(1-\sigma)^{12}+9 \sigma (1-\sigma)^{10}-40 \sigma^3 (1-\sigma)^6+144 \sigma^5 (1-\sigma)^2+256 \sigma^6 .$$
One can show that $p_1$ is positive on $0\le\sigma\le{1}$.\\
The negative quantity $\Sigma_1$ is found as in~(\ref{eq:tangPSig1}):\\
\begin{eqnarray}
\Sigma_1
&=&
\frac{1}{16}\left( \frac{i_2^2}{L}  -  i_4 \right),
\nonumber\\
&=& -\frac{1}{360}\,\frac{A^3}{L\,(\sigma (1-\sigma) (1+\sigma))^3}\ p_1(\sigma) .
\label{eq:isosS1}
\end{eqnarray}


Our main test cases are the equilateral triangle which has $\sigma=2-\sqrt{3}$ and
the right isosceles triangle which has $\sigma=\sqrt{2}-1$.
(Numerical values for the torsional rigidities of other isosceles triangles are 
available, for example in~\cite{PSH54}.)

It would, of course, be possible to produce tables, 
as done in \S\ref{subsec:Generaln} in the different context of regular polygons,
for a range of vertex angles for the isosceles triangles.
The starting point for this would be existing results for $\QQ_0$,
combined with our formulae for $L$, $\Sigma_\infty$~(\ref{eq:isosPSigInf2})
and $\Sigma_1$~(\ref{eq:isosS1}).

\newpage

 \subsection{The right isosceles triangle}\label{subsec:rightIsos}

\begin{table}[h]
\begin{center}
\resizebox{\textwidth}{!}{
\begin{tabular}{|| c | c | c | c| c | c | c ||}
\hline
$\alpha$& $4\QQ_0/A^2$&$A_n$& $\QQ_0$ & $L_n$ & $\Sigma_\infty$& $-\Sigma_1$\\ 
\hline
$\pi/3$& $\sqrt{3}/15$&   $3\sqrt{3}/4$ &$9\sqrt{3}/320$& $3\sqrt{3}$ &  $3\sqrt{3}/64$&  $3\sqrt{3}/320$ \\ 
$\rho=1/2$ & 0.11547& 1.2990 & 0.0487 & 5.1962 & 0.0812&  0.0162\\ 
\hline
$\pi/2$& 0.10436 &  $1$ &   0.02609&  $2+2\sqrt{2}$ &  $(3-2\sqrt{2})/3$&  
$(131-91\sqrt{2})/90$\\ 
$\rho=\sqrt{2}-1$  & & 1 & & 4.8284 & 0.0572 & 0.0256285 \\ 
\hline
\end{tabular}
}
\caption{Isosceles triangles with circumradius 1}
\label{tbl:tblisos}
\end{center}
\end{table}

\section{Tangential quadrilaterals, especially kites and rhombi}\label{sec:tangQuad}

A quadrilateral is tangential if and only if the sums of lengths of each pair of opposite sides are equal.
Examples of tangential quadrilaterals are the kites, which include the rhombi, which in turn include the squares.
(A bicentric kite is an orthogonal kite: 
a bicentric rhombus is a square.)
Torsional rigidities have been found numerically, for rhombi in~\cite{RI54,SC65}.
The other quantities occuring in the lower bound $R$ are easily found.

For example, the relevant quantities for a rhombus are as follows.
Consider a rhombus with inradius $\rho$, area $A$.
Let $\alpha$ be the angle at an acute vertex.
The points of tangency of the incircle with the sides of the rhombus
divide each side into a smaller part $\eta_-$ and a larger part $\eta_+$.
Denote $\tan(\alpha/2)$ by $\tau$.
Then $\eta_-=\rho\tau$ and $\eta_+=\rho/\tau$.
From equation~(\ref{eq:i0eta})
$$ L=i_0 =4\rho\left( \tau +\frac{1}{\tau}\right) ,\quad
A =\frac{1}{2}\rho L = 2\rho^2\left( \tau +\frac{1}{\tau}\right), \quad
\rho =\sqrt{\frac{A}{2\left( \tau +\frac{1}{\tau}\right)}} . $$
(At fixed $A$, $L$ is minimized for the square, $\tau=1$.)
We have, from~(\ref{eq:i2eta},\ref{eq:i4eta}),
\begin{eqnarray*}
i_2 &=& \rho^3\left(
4\left( \tau +\frac{1}{\tau}\right)+
\frac{4}{3} \left( \tau^3 +\frac{1}{\tau^3}\right)
\right) ,\\
i_4 &=& \rho^5\left(
4\left( \tau +\frac{1}{\tau}\right)+
\frac{8}{3} \left( \tau^3 +\frac{1}{\tau^3}\right)+
\frac{4}{5} \left( \tau^5 +\frac{1}{\tau^5}\right)
\right) .
\end{eqnarray*}
These check, in the case $\tau=1$ with the quantities given in~\S\ref{subsec:square}.
Equations~(\ref{eq:tangPSigInf},\ref{eq:tangPSig1}) give 
$\Sigma_\infty$ and $\Sigma_1$ in terms of $A$, $L$, $i_2$ and $i_4$.
The results for general rhombi are used in
Part~IIb~\S\ref{subsubsec:BStangQuad} and in
Part~III~\S\ref{sec:RhombiQ}.

There are many geometric results concerning tangential quadrilaterals.
A tangential quadrilateral is bicentric if and only if its inradius (hence area) is greater than that of any other tangential quadrilateral having the same sequence of side lengths.
There may be similar results for some other domain functionals.

Amongst all quadrilaterals with a given area that which\\
minimizes perimeter $L$,\\
 maximizes $\QQ_0$ (or similarly $\dot r$)\\
is square: see~\cite{PoS51} p159.
The proof in~\cite{PoS51} involves symmetrisation, with kites to rhombi to rectangles then kites, etc..
Perhaps because of curiosity on how the successive symmetrisations performed
there have been numerical studies of the torsion problem for rectangles, kites and rhombi.
For rhombi an early example is~\cite{RI54}.

\newpage
\begin{center}
{\large{{\textsc{ Part IIb:
More geometry
for tangential polygons
}}}}
\end{center}

\section*{Abstract for Part IIb}
Further items on tangential polygons,  
additional to those in Part IIa (which are taken from~\cite{Ke20i}),
are collected here.
In particular some Blaschke-Santalo diagrams for some geometric functionals
are presented.


\section{Outline of Part IIb}\label{sec:OutlineIIb}

In as much as the main focus of these notes should be the lower bound 
$\QQ_{0-}$ of Part I, we remark that Parts IIa and IIb only provide
results on the geometric quantities entering the formula for $\QQ_{0-}$,
namely $\rho$, $A$, $L$, $i_2$, $i_4$.
We defer further treatment of $\Sigma_1$ and $\QQ_{0-}$ to Parts III and IV.
\medskip

There are two main, and different, sorts of geometries.
\begin{itemize}

\item One, especially prominent in~\S\ref{sec:BlSa},
especially~\S\ref{subsec:cap} and~\S\ref{subsec:rLR},
involves general tangential polygons, convex circumgons.
(When we count the extreme points outside the incircle these are,
if $n$ points, called circum-$n$-gons.)

\item The other concerns genuine $n$-gons.
Their boundaries consist solely of straight line segments.
The hope is that one can show that regular $n$-gons optimize
appropriate domain functionals over subsets of, sometimes all,
tangential $n$-gons.
The hope is sometimes realized, an easy example (in~\S\ref{subsec:triGenz})
being as follows:\\
{\bf Result.}{\it Let $\Omega_0$ be a tangential $n$-gon with inradius $\rho$,
and vertex angles $\alpha_k$. 
Let $\Omega_1$ be a tangential $n$-gon with the same inradius $\rho$
and the same vertex angles except that vertices $\alpha_i$ and $\alpha_j$
are each replaced by their average $(\alpha_i+\alpha_j)/2$.
Then each of $L$, $A$, $i_2$, $i_4$ and $d_O$ are reduced in the
change from $\Omega_0$ to $\Omega_1$, strictly so if
$\alpha_i\ne{\alpha_j}$.}\\
The `easy' above refers to the proof.
It wasn't immediately obvious to this author before the proof.
By way of contrast, the following seems almost self-evident.\\
{\bf Corollary.}{\it Amongst all tangential $n$-gons with a fixed inradius,
the regular $n$-gon minimizes each of $L$, $A$, $i_2$, $i_4$ and $d_O$.}\\
We have yet to check the possibility that $\QQ_{0-}$ behaves,
in this respect, the same as~\cite{Sol92} Theorem 1 gives for $\QQ_0$.
Discussion of this is defered to Part III.

\end{itemize}

Here is an outline of this part.

\begin{itemize}

\item In \S\ref{sec:Transformations}
we consider operations involving
tangential polygons.

\item In \S\ref{sec:Construction} we indicate how we coded to construct
tangential polygons in order to later compute domain functionals for them.

\item In \S\ref{sec:Duality} considerations of `duality' direct the study.

\item In \S\ref{sec:propIIb} we collect a somewhat miscellaneous set of
inequalities and geometric facts.

\item In \S\ref{sec:IIbTri}
we consider triangles:
in \S\ref{sec:quad} tangential quadrilaterals.

\item In \S\ref{sec:Bicentric} we note a few facts concerning bicentric polygons.
All triangles are bicentric. All regular polygons are bicentric.

\item In \S\ref{sec:BlSa} we present some information about Blaschke-Santalo diagrams for some geometric functionals.

\end{itemize}

Blaschke-Santalo results can sometimes lead in to proving isoperimetric
results. Here is the style of an example with $\cal F$ some domain functional.
\begin{itemize}
\item Amongst triangles with fixed $\rho$ and $L$ that which
$<$ optimizes $\cal F$ $>$
is $<$ squat $|$ tall $>$ isosceles.

\item Amongst \ldots isosceles triangles at fixed $A=\rho L/2$
that which $<$ optimizes $\cal F$ $>$
is equilateral.

\end{itemize}
As an example consider $d_O$ (defined and treated extensively 
in~\S\ref{sec:BlSa}).
At fixed $\rho$ and $A=\rho L/2$ squat isosceles triangles minimize $d_O$.
With this preliminary, when considering triangles with given $A$,
minimizing $d_O$ we need only consider isosceles triangles with that 
area.
(For $Q_{0-}$ I have in Part I seen that, amongst isosceles triangles
with a given area that which maximizes $\QQ_0$ is equilateral.)
See~\S\ref{sec:IIbTri} for more formal treatment, but for now
the following very informal discussion might make the second step plausible.
Consider now squat isosceles triangles with given incentre and
its base vertices on the circle radius $d_O$.
It is eminently plausible that increasing the inradius of this family of triangles will increase the area, suggesting that at fixed $d_O$ one gets
maximum $A$. Conversely one expects at fixed $A$ to get minimum $d_O$
at the equilateral triangle.

\section{Transformations involving tangential polygons}\label{sec:Transformations}

Changing scale by some factor $t$ changes
a tangential polygon with inradius $\rho$
to one with inradius $t\rho$.
Mostly we fix the inradius, and always have the origin of coordinates at
the centre of the incircle.
In this situation, as
we have already noted, in Part~IIa~\S\ref{subsec:tangGeom} that
given two (convex) tangential polygons with the same incircle
their intersection is also a (convex) tangential polygon with the same incircle.
(For tangential $n$-gons the number of vertices could increase.)

\subsection{Convex $n$-gon to tangential $n$-gon}

Given a convex $n$-gon, and hence its sequence of angles, all less than $\pi$,
then one can define a tangential $n$-gon with the same sequence of angles and
the same area.
This transformation reduces the perimeter. (See~\cite{AM04}.)

\subsection{Tangential $n$-gons to tangential $m$-gons, $m\ge{n}$}

As a particular case of the intersection of tangential polygons being tangential polygon we mention ``corner cutting".
Given a tangential polygon $\Omega$ and a half-plane $H$ containing the incircle of $\Omega$, then $H\cap{\Omega}$ is  a tangential polygon.
Both $L$ and $A$ are decreased while $\rho=2A/L$ stays constant.
The number of vertices increases.

Let $\Omega$ be a tangential polygon.
Let $\sigma(l,.)$ be reflection through a line $l$ through the origin.
Then $\sigma(l,\Omega)$ is a tangential polygon and so is its intersection with
$\Omega$.\\

In the case of tangential $n$-gons, the numbers of vertices changes.
For example, if $\Omega$ is an equilateral triangle and $l$ is parallel to a side,
the intersection is a hexagon.

\subsection{Tangential $n$-gons to tangential $m$-gons, $m\le{n}$}

Begin with a tangential polygon with $n\ge{4}$ vertices.
Moving a tangency point to an adjacent tangency point will result in an edge disappearing.
While the inradius stays the same, the new tangential polygon might not be bounded with
an example of this having the starting point as a square.

\subsection{Permutations of angles of tangential $n$-gons}

If one permutes the entries of a sequence of angles or a tangential polygon,
keeps the inradius the same,
one has another tangential polygon 
with the same $\rho$, $A$, $L$, $i_2$, $i_4$, $d_O$.
As an example consider reflection about any `diagonal':
the incentre moves, but the reflected polygon is tangential.

\subsection{Tangential $2m$-gons to 2-special $2m$-gons}

Another transformation which at least takes tangential quadrilaterals to
tangential quadrilaterals is $m$-descendant mapping: see~\cite{Wo81}.
This paper also defines a $n$-gon with $n=2m$ even to be {\it 2-special} if
the sum of lengths even-indexed sides is equal to that of the odd-indexed sides.
Also the {\it 2-descendant map} of a $n$-gon with sides $s_{{\rm in},j}$ is
the $n$-gon with sides
$$ s_{{\rm out},j} =\frac{1}{2}(s_{{\rm in},j}+s_{{\rm in},j+1} ). $$
The 2-descendant map of any 2-special $2m$-gon is 2-special.
In particular the 2-descendant map of a tangential quadrilateral is a tangential quadrilateral.
Repeated application of 2-descendant maps to an initial tangential $2m$-gon
would take one ever closer to a regular $2m$-gon.
We remark that the doubly stochastic circulant matrix $\frac{1}{2}M(n)$,
defined in~\S\ref{subsec:Circulant} is here applied to
${\mathbf s}_{\rm in}$.
The effect of many successive operations with the map leading to the
regular $n$-gon corresponds to the fact that the matrix powers tend
to the $1/n$ times the matrix $E$ all of whose entries are 1:
$$  \left(\frac{1}{2}M(n)\right)^k \rightarrow \frac{1}{n} E(n)\qquad
{\rm as\ \ } k\rightarrow\infty . $$

{\bf ToDo.} Check out guess that applying a 2-descendant map to a tangential hexagon
may not lead to tangential hexagon.
(It is known that being 2-special is necessary but not sufficient for a hexagon to be
tangential.)

\medskip
One can also consider $2m$-gons as linkages.
Again a tangential $2m$-gon can move as a linkage to another 2-special configuration.
For tangential quadrilaterals the linkage remains, when convex, a tangential quadrilateral
(but this will not be the case for 6-gons).

\subsection{Tangential $n$-gons to Tangential $n$-gons, preserving $\rho$}

\begin{figure}[ht]
\centerline{\includegraphics[height=10cm,width=14cm]{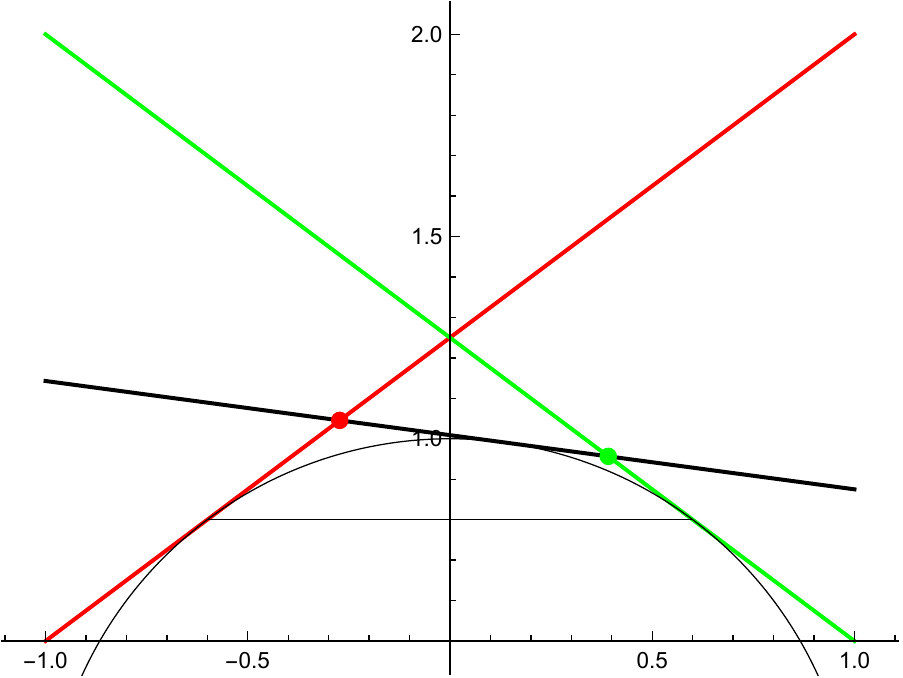}}
\caption{Diagram for `tilting transformation'}
\label{fig:IIbTilt}
\end{figure}

A transformation, called a `tilting transformation' is now defined.
In this context refer to Figure~\ref{fig:IIbTilt}.
A side is `tilted', its tangent contact point moved, so that
it becomes parallel to the line joining the contact points 
of adjacent sides.
All tangent points except one remain fixed.
Suppose the initial configuration is as shown in the figure.
The red and green tangent lines remain fixed, but suppose the
tangent point of the black side is considered movable.
The fixed tangent points either side of the movable one are at
$$ \zeta_- = (-X,h), \qquad \zeta_+ = (X,h) .$$
The coordinates of the movable tangent point might be taken as
$$ \zeta = \rho \left(\frac{1-t^2}{1+t^2},\frac{2t}{1+t^2}\right) \ \ 
{\rm\ and\ w.l.o.g.\ } \rho=1 .$$
The area $A(t)$ of the quadrilateral obtained from intersecting
the tangential $n$-gon with the half-space $\{y>h\}$
can be found, as can the length $L(t)$ of the three line segments
in the upper part of its perimeter.
Finding the formulae for $A(t)$ and $L(t)$ is not too onerous,
with one check being that $2A(t)-L(t)$ is independent of $t$.
Both $A(t)$ and $L(t)$ are minimized at $t=1$ which has the
topmost (black) tangent line parallel to $y=h$ and
the two consecutive angles of the tangential $n$-gon equal.

\medskip
Repeated application of the tilting transformation,
for varying sides, will,
in the limit, get one to a tangential $n$-gon with all angles equal.
If a tangential $n$-gon has all angles equal it is regular.
Thus, amongst all tangential $n$-gons with radius $\rho$,
that which has the smallest area (and perimeter $L=2A/\rho$)
is the regular $n$-gon.
\medskip

Similarly, amongst all tangential $n$-gons with radius $\rho$,
that which has the smallest $d_O$ is the regular $n$-gon.

\medskip

The `tilting transformation' is easy to visualize geometrically.
Less obvious is that at fixed $\rho$ one can reduce 
$A$, $L$, $i_2$, $i_4$, $d_O$ by averaging any two angles
(not necessarily adjacent angles as above).
The reduction of $L$ is a consequence of the convexity of the $\cot$
function on $(0,\pi/2)$, the quantities $T_k$ defined 
in Part IIa equation~\ref{eq:Tkdef}) and the representation of
$L$ as the sum of the $T_k$.
The other quantities $i_2$, etc. are  treated similarly.
More details are given in~\S\ref{subsec:triGenz}.


\subsection{Circum-$n$-gons to circum-$m$-gons, $m\le{n}$}

Let $\Omega$ be a tangential polygon with inradius $\rho_0$ and
with the maximum distance of the boundary to the origin $d_O$.
Then, for any ball $B$ centred at the origin, the convex hull
${\rm conv}(B\cup{\Omega})$ is a tangential polygon.\\
{\bf ToDo.} Prove this
and investigate the behaviour of $L/\rho$ as the radius of 
$B$ increases from $\rho_0$ to $d_O$.

\subsection{Minkowski sums}

The Minkowski sum of convex polygons is a convex polygon.

The Minkowski sum of two line segments is a parallelogram.\\
The Minkowski sum of two triangles is a hexagon (usually not tangential).\\
Any convex polygon is the Minkowski sum of triangles and line segments.
A reference for this (which I have yet to check) is page 177 of \\
I. M. Yaglom, V. G. Boltyanskii, (1961)
{\it Convex Figures}, New York: Holt, Rinehart and Winston.\\
This leads on to the following.

\noindent
{\bf Questions.} (i) Is {\it any} convex hexagon the Minkowski sum of 2 triangles?\\
(ii) Is, for $n\ge{6}$, any convex $n$-gon with origin at the centroid
the Minkowski sum of a small number (perhaps just 2) of tangential $m$-gons
(with $m\le{n/2}$),
now (unlike everywhere else in this document)  all with centroid at the origin?\\
A positive answer to a question like this might suggest, from
properties established for all tangential $m$-gons,
corresponding properties for convex $n$-gons.
And if this were to be the case, one can imagine establishing some 
domain-functional property 
for tangential $n$-gons and then getting something from the concavity of
the domain-functional under Minkowski sum.

\subsection{Rearrangements???}

For the purely geometric functionals $A$, $L$, $i_2$, $i_4$, $\QQ_{0-}$, $d_O$
studied in this document rearrangements might not be needed.
For functionals, like conformal inradius, transfinite diameter,
torsional rigidity, fundamental frequency -- functionals involving
integrals of gradient squared, etc. -- rearrangements are an appropriate tool.
\cite{PoS51} used Steiner symmetrization to establish, for triangles
and quadrilaterals, that, at given area, the regular $n$-gon optimizes.
These are equilateral triangle and square respectively.
However, Steiner symmtrizing polygons with  more vertices typically
increases the number of vertices.
\cite{SZ04,Sol92,SZ10} manage to use some sort of rearrangement preserving the
number of vertices of a convex $n$-gon. Dissymetrization?\\
{\bf ToDo.} Try to understand this. See also~\cite{Ba06}
\medskip

Polarizations seem to be building blocks for the rearrangements with which
I am more familiar, namely those used in~\cite{PoS51}.
Polarizations, alone, are not likely to be a tool for rearranging 
tangential polygons.
The following guesses and questions began with drawing sketches.

\begin{itemize}

\item The polarization of a triangle about an angle bisector just reflects the
triangle.

\item Can anything beyond the incircle staying fixed be said about the polarization of a triangle about any line through the incentre?
Similar question for any tangential polygon $\Omega$ about any line through the incentre?
Non-convex circumgons?

\end{itemize}

\cite{Sol92}\\
 (i) presents results of the kind that regular $n$-gons optimize over all
$n$-gons with the same area;\\
(ii) that tangential polygons are used (see Part III~\S\ref{sec:BoundsQIso}).\\
A process called `dissymetrization' gets used.
This, and polarization, get a mention in~\cite{SZ10}.

\subsection{Spaces of polygons?}

See~\cite{GPT17}.


\section{Construction of tangential polygons}\label{sec:Construction}

With just a few exceptions (on 1-cap, etc.)
to date our computations have been for tangential $n$-gons.

Our first method provided data for actually drawing up the polygon
and required an initial specification of the inradius.
After this prescribe $n$ (which in our computations so far just $n\le{6}$).
Then choose $n$ increasing values of $\theta_k$ in $(-\pi,\pi)$
with the maximum difference of consecutive $\theta_k$ less than $\pi$.
This yields $\exp(\theta_k)$ on the unit circle as tangent points.
(The restriction on the separation of the $\theta_k$ is
in order that a convex polygon is constructed.)
From each pair of consecutive tangent points, find the point of
intersection of the tangent lines.
These give the vertices of the tangential polygon and
there are standard formulae for perimeter $L$, area $A$,
in terms of the coordinates but it is easier to note that
with the tangent lengths one can find the $T_k$ and use
these to find $L$, $i_2$, $i_4$, etc.

\medskip

If one doesn't need to draw the polygon and is interested in
functionals that stay constant under change of scale,
e.g. $L/\rho$, one can use the $T_k$ as defined in
Part IIa equation~(\ref{eq:Tkdef}).
The tangent lengths are given by $\eta_k=\rho\,T_k$.
Simple formulae to find the other functionals $L$, $i_2$, $i_4$
are given Part IIa~\S\ref{subsec:tangGeom},
and another $d_O$ in~\S\ref{sec:BlSa}.

We start from result given at the beginning of~\S\ref{sec:Transformations}:\\
{\bf Existence Theorem.}{\it Given a convex $n$-gon, 
and hence its sequence of angles, all less than $\pi$, 
then one can define a tangential $n$-gon
whose incentre is at the origin, with the same sequence of angles.
}\\
Clearly once one has one, one can rescale by any factor preserving
the properties.
\medskip

\noindent{\bf Corollary. }{\it Given a set $S$ of $n$-numbers $0<\alpha_k<\pi$
summing to $(n-2)\pi$,
then the different tangential $n$-gons arising from the different
permutations of $S$ all have the same values for
$$\frac{L}{\rho},\ \
\frac{i_2}{\rho^3},\ \
\frac{i_4}{\rho^5},\ \
\frac{d_O}{\rho}.
$$
}
\medskip

In any event, for many calculations later in this part,
one can start directly with the tangent lengths,
or with the $\alpha_k$ 
or with the $T_k$: 
it is not always necessary to calculate the vertex coordinates.

In some contexts moving between tangential $n$-gons by permuting
angles may lose properties.
In particular, permuting angles of a bicentric $n$-gon will,
in general result in a tangential $n$-gon which is not bicentric.
(The simplest example would be any bicentric quadrilateral
with 4 different angles.
Permuting the angles must lose the property that the sum of
opposite angles is $\pi$.)

While $d_O$, the distance from the incentre origin to a vertex,
doesn't need the vertex coordinates and 
and is invariant, at fixed $\rho$ under permutations of the angles,
this may not be the case for the circumradius $R$.

\section{Duality} \label{sec:Duality}

Denote the inner product for plane vectors with a dot.
Define {\it polar-reciprocation} $\cal P$ of a point $x$ by
$$ {\cal P}(x) = \{ z\ |\ z\cdot{x} = 1 \} .$$
$\cal P$ takes a point to a line.
(This differs from the most common definition of polar in convex geometry in which
one has $z\cdot{x}\le{1}$ so points map to half-planes.)
Next continue the definition. Let $D$ be a set in the plane.  
Define ${\cal P}(D)$ by
$$ {\cal P}(D) = \{ z\ |\ z\cdot{x} = 1 \ \forall x\in{D} \} .$$
$\cal P$ takes the unit circle to itself.
$\cal P$ takes lines (not through the origin) to points:
in particular $\cal P$ takes a line tangent to the unit circle to its point of tangency with the unit circle.

See\\
\verb$https://en.wikipedia.org/wiki/Pole_and_polar$\\
\verb$https://en.wikipedia.org/wiki/Dual_polygon$

Thus the boundary lines of a tangential polygon map to the vertices of a cyclic polygon and vice-versa.

\bigskip

\noindent
{\bf `Vertex-side' duality, adapted from wikipedia}

As an example of the side-angle duality of polygons we compare properties of the cyclic and tangential polygons, especially quadrilaterals.

{\small
\begin{center}
\begin{tabular}{|c|c|}
\hline
Cyclic $n$-gon&	Tangential $n$-gon\\
\hline
\hline
Circumscribed circle&	Inscribed circle\\
\hline
Perpendicular bisectors of the sides are& 
Angle bisectors are\\
concurrent at the circumcentre& concurrent at the incentre\\
\hline
$n$ even: The sums of the two pairs& $n$ even: The sums of the two pairs\\
 of opposite/alternate angles&  of opposite/alternate sides\\
 are equal& are equal\\
\hline
\end{tabular}
\end{center}
}
For a cyclic $2m$-gon the sum of the alternate angles is $(m-1)\pi$.
We remark that for quadrilaterals, $m=4$ $n=2$ the converse is true
but it is false for $n\ge{3}$.
The same is true for tangential $2m$-gons:
for $m\ge{3}$ being `2-special' is necessary but not sufficient for a
$2m$-gon to be tangential.

$n=6$, $m=3$.  Brianchon's theorem states that the three main diagonals 
of a tangential hexagon are concurrent.\\
The polar reciprocal and projective dual of the conics version of
Brianchon's theorem give Pascal's theorem.

Duality is evident again when comparing an isosceles trapezoid to a kite.

{\small
\begin{center}
\begin{tabular}{|c|c|}
\hline
Isosceles trapezoid&	Kite\\
\hline
\hline
Two pairs of equal adjacent angles&	Two pairs of equal adjacent sides\\
\hline
One pair of equal opposite sides&	One pair of equal opposite angles\\
\hline
An axis of symmetry through& An axis of symmetry through\\
one pair of opposite sides& one pair of opposite angles\\
\hline
Circumscribed circle&	Inscribed circle\\
\hline
\end{tabular}
\end{center}
}

Let $P(n)$ be an ordered list of $n$ points $\zeta_k$ on a circle,
$P(n+1)=P(n)\cup\{\zeta_{n+1}\}$ with $\zeta_{n+1}$
after $\zeta_{n}$ and before $\zeta_{1}$.\\
Let $TP(P(n))$ be the tangential $n$-gon with the points of $P(n)$
its tangent points.\\
Let $CP(P(n))$ be the cyclic $n$-gon with the points of $P(n)$
its vertices.
The first entry in the table below indicates how areas change
on introducing the additional point on the circle.

\begin{center}
\begin{tabular}{|c|c|}
\hline
$|TP(P(n+1))|\le|TP(P(n))|$&
$|CP(P(n+1))|\ge|CP(P(n))|$\\
\hline
A tangential $2m$-gon has all sides equal iff&
A cyclic $2m$-gon has all angles equal iff \\
 the alternate angles are equal.&
the two sets of alternate sides are equal.\\
\hline
Equilateral tangential for $n=2m$  even&
Equiangular cyclic for $n=2m$ even\\
Opposite angles equal if $n/2$ is even&
Opposite sides equal if $n/2$ is even\\
\hline
\end{tabular}
\end{center}

See~\cite{deV11}.

\subsection{Tangential and cyclic polygons, continued}

The set of all convex sets is a lattice under operations of
intersection, $\Omega_1\cap\Omega_2$,
and convex-hull of union, ${\rm conv}(\Omega_1\cup\Omega_2)$.\\
The set of tangential polygons with incentre at the origin
and given inradius $\rho$ is closed under intersection.\\
The set of cyclic polygons with circumcentre at the origin
and given circumradius $R_V$ is closed under convex-hull-union.
(A cyclic polygon is the convex hull of its extreme points
all of which lie on the circle radius $R_V$.)\\
In the case of $n$-gons the numbers of vertices can increase.

\medskip

Let the coordinates of the vertices of a convex $n$-gon be $(x_k,y_k)$
with the vertices traversed in order 
(and vertex 1 identified with vertex $n+1$).\\
The polygon is cyclic with circumcentre O and circumradius 1 if 
the distance of every vertex from O is 1:
$$x_k^2 + y_k^2=1\qquad \forall k . $$
The polygon is tangential with incentre O and inradius 1 if 
the distance of every side from O is 1:
$$(x_{k+1}-x_k)^2 + (y_{k+1}-y_k)^2=  (x_{k+1} y_k - y_{k+1} x_k)^2
\qquad \forall k . $$
The tangency point on each line, the closest point to O, is
$$ x_t = \frac{y_{k+1}-y_k}{x_k\, y_{k+1}-y_k\, x_{k+1}}, \qquad
y_t = \frac{x_{k+1}-x_k}{x_k\, y_{k+1}-y_k\, x_{k+1}} . $$

\medskip

In establishing, by Steiner symmetrisation, isoperimetric properties
of $n$-gons, for $n=3$ and $n=4$, sequences of polygons which alternate
between tangential and cyclic occur.
Here is a quote from~\cite{PoS51}:
\begin{quote}
{\bf Of all quadrilaterals with a given $A$, the square has the smallest 
$L$, $I_c$ (polar moment of inertia about the centroid), 
$\overline r$, \ldots but the largest $\dot r$ and $\QQ_0$.}
\ldots\ldots it is sufficient to indicate a sequence of symmetrizations 
which transform, ultimately, a given quadrilateral into a square. 
Symmetrizing a given quadrilateral with respect to a perpendicular to one of ·its diagonals, we change it into a quadrilateral having a diagonal
as axis of symmetry. 
Symmetrizing this new quadrilateral with respect to a
perpendicular to its axis of symmetry, we change it into a rhombus. 
Symmetrizing the rhombus with respect to a perpendicular to one of its sides, we change it into a rectangle. 
Symmetrizing the rectangle with respect to a
perpendicular to one of its diagonals, we obtain another rhombus. 
Repeating the last two steps in succession, we obtain an infinite sequence in which rhombi alternate with rectangles.
\end{quote}
Rhombi are tangential polygons (with equal sides): \\
rectangles are cyclic (with equal angles).
\medskip

Steiner symmetrization is not applicable to showing that regular $n-gons$ 
optimize when $n\ge{5}$.
If one (initially at least) focuses on geometric quantities like
$$\frac{A}{L^2},\ \frac{i_2}{L^3},\ \frac{i_4}{L^4} ,
\qquad{\rm and\ }\ \frac{\QQ_{0-}}{A^2} ,$$
it may be possible to devise other transformations between $n$-gons\\
which alternate between tangential and  cyclic,\\
which change functionals like those immediately above monotonically and\\
which converge to the regular $n$-gon.
\medskip


\section{Miscellaneous properties of tangential polygons}\label{sec:propIIb}

\subsection{Circulant matrices and tangential $n$-gons}\label{subsec:Circulant}

This subsection treats questions like the following:\\
Given a set of $n$ of positive side lengths $(s_j)$
how can we recognize if there could be a tangential polygon
with these side lengths?

We begin with a connection between tangential polygons and
circulant matrices presented near the beginning of the wikipedia
page on tangential polygons.

Let $P$ and $M=I+P$ be the $n\times{n}$ circulant matrices as follows.
$P$ is the $n\times{n}$ cyclic permutation:
\begin{equation*}
P =
\begin{pmatrix}
0 & 1 & 0& \cdots & 0 \\
0& 0 & 1 & \cdots & 0 \\
\vdots&\vdots&\vdots& \ddots & \vdots  \\
1 & 0 & 0& \cdots & 0
\end{pmatrix}
\end{equation*}
The matrix $\frac{1}{2}M$ is doubly stochastic.

The wikipedia page states:
\begin{quote}
There exists a tangential polygon of $n$ sequential sides $s_1,\ldots, s_n$
if and only if the system of equations
\begin{equation}
M{\mathbf \eta} = {\mathbf s} ,\label{eq:Metas}
\end{equation}
has a solution $(\eta_1,\dots, \eta_n)$ in positive reals.
If such a solution exists, then $(\eta_1,\dots, \eta_n)$ are the tangent lengths
of the polygon
(the lengths from the vertices to the points where the incircle is tangent to the sides).
\end{quote}
Once one has the $\eta$ and $\rho$ the angles $\alpha_k$ are determined from
$$\eta_k (=\eta_{k-}) = \rho T_k \qquad{\rm where\ \ }
T_k = \frac{1}{\tan(\frac{\alpha_k}{2})} .$$
See PartIIa, equation~(\ref{eq:Tkdef}).
That the sum of the $\alpha_k$ is $(n-2)\pi$ leads to one further equation
which we record, as an aside, in the next subsubsection.

\subsubsection{Some relations between the $T_k$}

As before, consider $n$-gons and
denote the angle at vertex $k$ by $\alpha_k$,
with $k$ increasing as one goes around the convex $n$-gon in counterclockwise direction.
The sum over all the $\alpha_k$ is $(n-2)\pi$.

Fix the inradius $\rho$ as 1.
Suppose the points of tangency of the tangential $n$-gon are
$\zeta_j=\exp(i\phi_j)$.
Again $j$ increases as one goes around the convex $n$-gon in counterclockwise direction.
The angle at O formed by the lines $O\zeta_j$ and $O\zeta_{j+1}$ is
$\phi_{j+1}-\phi_j$.

\vspace{1cm} 

Then, with
\begin{equation}
T_k  = \frac{1}{\tan\frac{\alpha_k}{2}}=\tan(\frac{\pi-\alpha_k}{2}),
\label{eqBS:Tkdef}
\end{equation}
$T_k>0$ since the polygon is convex.
 Since we know the sum of $\alpha_k/2$:
\begin{equation}
\sum_{k=1}^n \frac{\alpha_k}{2}
 =\sum_{k=1}^n {\rm arccot}(T_k) = (n-2)\frac{\pi}{2} .
\label{eqBS:arccot}
\end{equation}
 Some, but not all, the information in this can be expressed in equations involving just
 rational functions of the $T_k$, i.e. without the transcendental arccot function.
 We denote the elementary symmetric polynomial of degree $k$ by
$$ {\rm SymmetricPolynomial}(k, \ldots ) ,$$
 and define
 $$ e_k =  {\rm SymmetricPolynomial}(k,[\frac{1}{T_1},\frac{1}{T_2},\ldots,\frac{1}{T_n}]) .$$
 For tangential $n$-gons, we first treat $n$ odd, then $n$ even.\\
 When $n$ is odd $\cos((n-2)\pi/2)=0$ so
 the cosine of the sum of all the $\alpha_k/2$ is $0$ we have
 \begin{equation}
 e_0 -e_2 + e_ 4- e_6 \dots  = 0 .
\label{eqBS:elPec}
\end{equation}
When $n$ is even $\sin((n-2)\pi/2)=0$ so
 the sine of the sum of all the $\alpha_k/2$ is $0$ we have
\begin{equation}
 e_1 -e_3 + e_ 5- e_7 \dots  = 0 .
\label{eqBS:elPes}
\end{equation}

We will need the formulae for perimeter, i.e. $L=i_0$, for $i_2$ and
for $i_4$ as given in Part IIa~\S\ref{subsec:tangGeom} namely
equations~(\ref{eq:i0Tgen}), (\ref{eq:i2Tgen}) and~(\ref{eq:i4Tgen}).

\subsubsection{Examples at $n=3$ or $4$}
For a triangle
$ A={\sqrt{\frac{L}{2}(\frac{L}{2}-a)(\frac{L}{2}-b)(\frac{L}{2}-c)}}$ and
since $L=2\sum\eta_k$
$$\leqno{\rm tang3gon:}\qquad
A= \sqrt{(\eta_1+\eta_2+ \eta_3)\eta_1\eta_2\eta_3 } .$$
For a tangential quadrilateral wikipedia gives
$$\leqno{\rm tang4gon:}
A= \sqrt{(\eta_1+\eta_2+ \eta_3+\eta_4)
(\eta_1\eta_2\eta_3 + \eta_2\eta_3\eta_4 + \eta_3\eta_4\eta_1 +\eta_4\eta_1\eta_2)} .
$$

\begin{figure}[ht]
\centerline{\includegraphics[height=11cm,width=9cm]{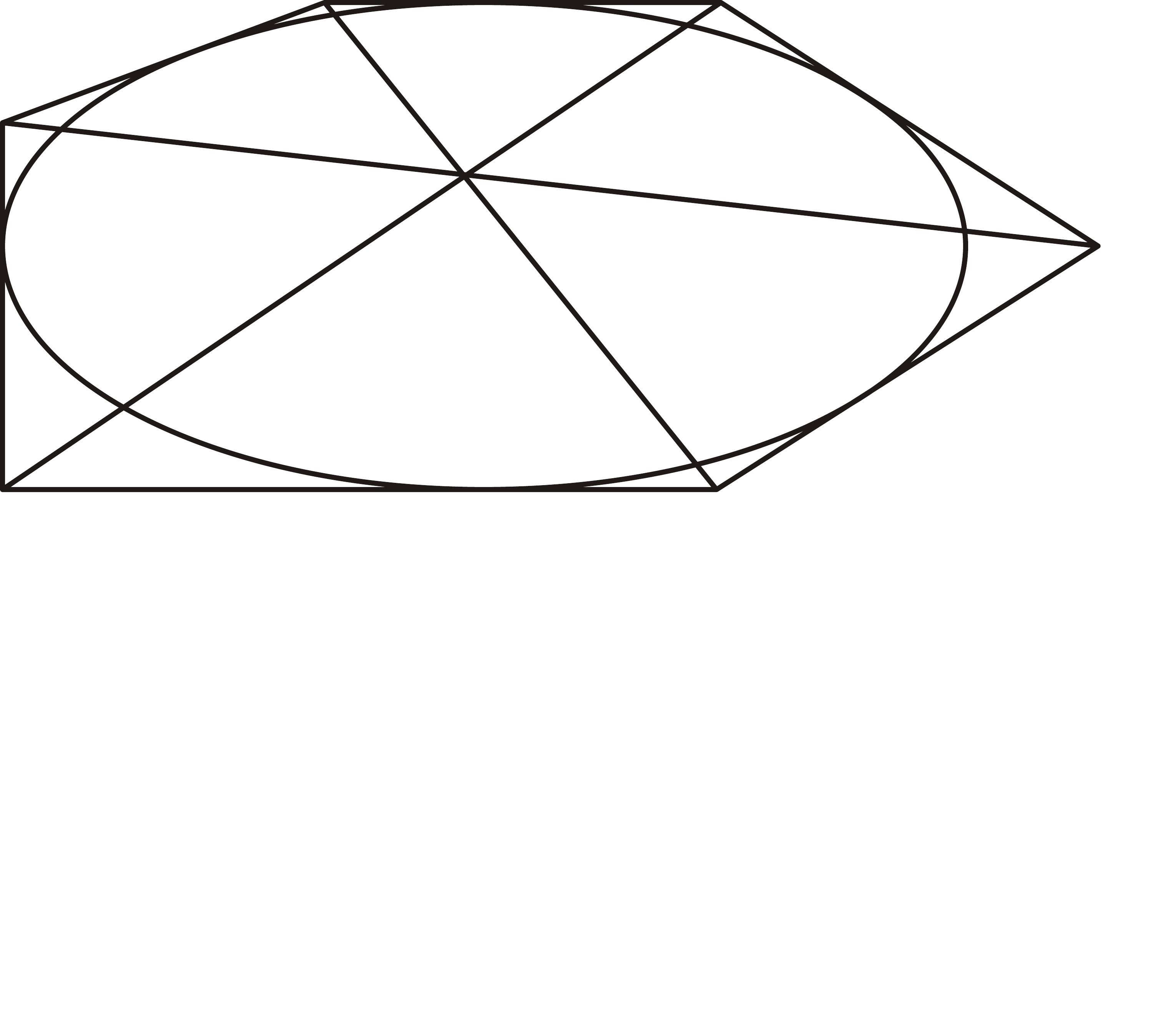}}
\vspace{-4cm}
\caption{From wikipedia. Brianchon's Theorem. Diagonals of a tangential hexagon are concurrent. }
\label{fig:Brianch}
\end{figure}

\subsubsection{Non-negative solution for $\mathbf\eta$ given $\mathbf s$?}

Return now to the linear equations~(\ref{eq:Metas}).
Denote the dependence of $M$ on $n$ by writing $M(n)$.
There is very different behaviour when $n$ is odd than when $n$ is even as,
for example,
$$ {\rm det}(M(n))= 1-(-1)^n .$$
Thus $M(n)$ is invertible when $n$ is odd, but not when $n$ is even.
There being many questions I have been, as yet, unable to answer, and
as properties of $M(n)$ might ultimately be useful, I have collected
several properties of $M(n)$ here, but not found uses for some of them yet.
Let ${\mathbf e}$ be the vector all of whose entries are 1.
Then, for all $n$,
$$ M(n)\, {\mathbf e} = 2 \, {\mathbf e} .$$
Geometrically this corresponds to a regular $n$-gon with side 2 and tangent lengths 1.
The eigenvalues of $P$ are the roots of unity, and the eigenvalues of $M$
require us just add 1 to these.
The eigenvectors are, of course, the same.
$${\rm CharacteristicPolynomial}(M(n),\lambda)
= -(-1)^n (1 - (\lambda - 1)^n) .
$$

The transpose $M^T(n)$ has exactly the same properties as described
for $M(n)$ in the preceding paragraph.
If one ignores the nonnegativity requirement, clearly there is,
for any rhs $\mathbf s$ a unique
`solution' to the linear equations when $n$ is odd.
When $n$ is even, there is only a `solution' when the rhs is
orthogonal to the nullspace of $M^T(n)$, i.e only when
\begin{equation}
\sum_{k\ {\rm odd}} s_k
=  \sum_{k\ {\rm even}} s_k ,
\label{eqIIb:evod}
\end{equation}
and, also, `solutions' are not unique.

$M(n)$ is normal, commutes with its transpose.
As $M\,M^T$ is symmetric and Toeplitz it is centrosymmetric.

A {\it persymmetric matrix} is
a square matrix which is symmetric with respect to the northeast-to-southwest diagonal.
$M(n)$ is persymmetric and, as such, satisfies
$$ M(n)\, J = J\, M^T(n),\qquad
{\rm where\ } J\ {\rm is\ the\ exchange \ matrix} .
$$
$J$ is the matrix with 1s on its northeast-to-southwest diagonal and 0s elsewhere.

Define next
$$M_i(n) = (I(n)+\sum_{k=1}^{n-1} (-1)^k P(n)^k )/2.$$
We have
$$ M(n)\, M_i(n) = \frac{1-(-1)^n}{2} I(n) . $$
For $n$ odd $M_i(n)$ is the inverse of $M(n)$.
Also (via Cayley-Hamilton Theorem)
$$ \sum_{k=1}^{n} (-1)^k {n\choose k} M(n)^k
= -(1-(-1)^n)\,  I(n) .$$

The obvious next question is `what conditions on the sides ensure there is
a {\em nonnegative} solution for $\mathbf \eta$'?
\smallskip

\noindent
{\bf Farkas Lemma.} {\it Exactly one of the following two assertions is true:\\
1.There exists an ${\mathbf \eta}\in{\mathbb{R}^n}$
such that $M{\mathbf \eta}={\mathbf s}$ and
${\mathbf {\eta}} \geq 0$.\\
2. There exists a  $\mathbf {y} \in \mathbb {R} ^{n}$ such that
$\mathbf {M} ^{\mathsf {T}}\mathbf {y} \geq 0$ and
$\mathbf {s} ^{\mathsf {T}}\mathbf {y} <0$.}\\
We have yet to devise a good use  for this lemma but suspect it will be relevant
to conditions for general $n\ge{3}$.
\medskip

Here we do not attempt to find the conditions for general $n$.
We will attempt to find necessary and sufficient conditions on $\mathbf s$ for the
existence of nonnegative $\mathbf\eta$, separately, for
each of $n=3$, 4, 5 and 6.

\medskip
\noindent{\bf Triangles.} When $n=3$,
$$ M(3)^{-1}= \frac{1}{2}
\left(
\begin{array}{ccc}
 1 & -1 & 1 \\
 1 & 1 & -1 \\
 -1 & 1 & 1 \\
\end{array}
\right) ,
$$
so one has nonnegative solutions for $\eta$ iff the nonnegative $s$
are such that the sum of any two is greater than (or equal to)
the remaining side's length.
This accords with the fact that any triangle is tangential
(indeed bicentric).

\medskip
\noindent{\bf Tangential pentagons.} When $n=5$,
$$ M(5)^{-1}= \frac{1}{2}
\left(
\begin{array}{ccccc}
 1 & -1 & 1 & -1 & 1 \\
 1 & 1 & -1 & 1 & -1 \\
 -1 & 1 & 1 & -1 & 1 \\
 1 & -1 & 1 & 1 & -1 \\
 -1 & 1 & -1 & 1 & 1 \\
\end{array}
\right) .
$$
Consider triples of sides in which just one pair of sides are adjacent.
There are five such triples.
There are nonnegative solutions for $\eta$ iff the nonnegative $\mathbf s$
are such that for any of these five triples
the sum of those in the triple is greater than (or equal to)
the sum of the remaining two sides' lengths.
\bigskip

Let $n=2m$ for $m\ge{2}$.
Define
$$ {\mathbf{nv}} = \left( \ (-1)^k \ \right) ,$$
which is a basis for the nullspace of $M(n)$
(and also of $M^T(n)$).
Suppose $\mathbf s$ satisfies equation~(\ref{eqIIb:evod}) so `solutions',
albeit without the nonnegativity condition imposed, exist.
These solutions are
\begin{equation}
{\mathbf \eta}_{\rm gen} =
{\rm PseudoInverse}(M(n))\,{\mathbf s} +c\, {\mathbf{nv}} ,
\label{eq:etaGen}
\end{equation}
but it remains to investigate the restrictions on $\mathbf s$ and $c$
needed so that
amongst the ${\mathbf \eta}_{\rm gen}$ there is at least one with ${\mathbf \eta}\ge{0}$.
In this document we will attempt this only for $n=4$ and $n=6$.
Before that, however, we consider $n=2m$ in general.
By standard properties
$$ M(n)\, {\rm PseudoInverse}(M(n))\, M(n) = M(n) .$$
Using the fact that when $n=2m$ the matrix $M(n)$ has the simple block structure
$$ M(n)
=\begin{pmatrix}
U& L\\
L& U
\end{pmatrix} ,
$$
it is easy to find $A$ and $B$ so that
$$ {\rm PseudoInverse}(M(n)
=\begin{pmatrix}
A& B\\
B& A
\end{pmatrix} .
$$
The matrix $L$ has just one 1 in the bottom left corner and
$U=M(m)-L$.
We have
\begin{eqnarray*}
A
&=&\frac{1}{2} ({\rm PseudoInverse}(U+L)+{\rm PseudoInverse}(U-L) ),\\
B
&=&\frac{1}{2} ({\rm PseudoInverse}(U+L)-{\rm PseudoInverse}(U-L) ) .
\end{eqnarray*}
We also have $L\,U\,L$ is the zero matrix, and
$L\,U^{-1}\,L=-L$.
\medskip

\noindent{\bf Tangential quadrilaterals.}
It is already known that no further restriction beyond
$$ s_1+s_3=\frac{L}{2}=s_2+s_4 $$
is needed to ensure that the quadrilateral is tangential.
So it remains just an exercise to show that the system of equations
$$ M(4){\mathbf \eta} =
\begin{pmatrix}
s_1\\
s_2\\
\frac{L}{2}-s_1\\
\frac{L}{2}-s_2
\end{pmatrix}
$$
has, for all $0<s_1<L/2$ and $0<s_2<L/2$ a positive $\eta$ solution.
When $n=4$,
$$ {\rm PseudoInverse}(M(4))
= \frac{1}{8}
\left(
\begin{array}{cccc}
 3 & -1 & -1 & 3 \\
 3 & 3 & -1 & -1 \\
 -1 & 3 & 3 & -1 \\
 -1 & -1 & 3 & 3 \\
\end{array}
\right) .
$$
There is no loss of generality in setting $L=2$ and
considering $1/2\le{s_1}<1$ and $1/2\le{s_2}<s_1$.
If we try the formula~(\ref{eq:etaGen})
we are led to consider the function $\phi$
$$\phi(s_1,s_2,c)
={\rm min}(1 + 2 s_1 - 2 s_2-c, -1 + 2 s_1 + 2 s_2+c,
1 - 2 s_1 + 2 s_2-c, 3 - 2 s_1 - 2 s_2+c),
$$
over the triangle in $(s_1,s_2)$ space defined in the preceding sentence.
Considering the final entry in the min defining $\phi$,
we have $\phi(s,s,0)<0$ for $3/4<s<1$, so we need to choose $c$
(which can depend on $\mathbf s$) appropriately.
We find, over the triangle in $(s_1,s_2$-space
$\phi(s_1,s_2,2s_1+2s_2-3)=0$, i.e. the last entry is zero but
the other 3 entries are nonnegative.
Except for a positive multiple, the other 3 entries
(first 3) are
$$ 1 - s_2, -1 + s_1 + s_2, 1 - s_1 . $$

\noindent{\bf Tangential hexagons.} Unlike the situation with $n=4$
extra conditions on the sides are needed.
The corresponding problem for cyclic hexagons is mentioned in
~\cite{deV16} (and the same author has other papers involving
tangential and cyclic hexagons,~\cite{deV02},~\cite{deV11}).

After writing the above, I found~\cite{BS05} gives the following.

\noindent
{\bf Theorem.}{\it There will be a tangential hexagon with given side lengths
$s_1,s_2,...,s_6$ if and only if the equality
$$ s_1 +s_3 +s_5 =s_2 +s_4 +s_6, $$
and the following nine inequalities are satisfied:
$$s_1 >0,s_2 >0,...,s_6 >0 , $$
\begin{eqnarray*}
s_1 - s_2 + s_3 &>& 0, \\
s_3 - s_4 + s_5 &>& 0, \\
s_5 - s_6 + s_1 &>& 0.
\end{eqnarray*}
}

\noindent
In words, the length of any side is less than the sum of the lengths
of the adjacent sides.
\medskip

\cite{BS05} also gives the conditions on $\mathbf s$ for an octagon to be
tangential.
When $n=6$,
$$ {\rm PseudoInverse}(M(6))
= \frac{1}{12}
\left(
\begin{array}{cccccc}
 5 & -3 & 1 & 1 & -3 & 5 \\
 5 & 5 & -3 & 1 & 1 & -3 \\
 -3 & 5 & 5 & -3 & 1 & 1 \\
 1 & -3 & 5 & 5 & -3 & 1 \\
 1 & 1 & -3 & 5 & 5 & -3 \\
 -3 & 1 & 1 & -3 & 5 & 5 \\
\end{array}
\right) .
$$
(On looking at the corresponding outputs for large $n=2m$ one finds that
$2n$ times ${\rm PseudoInverse}(M(n))$ has entries $\pm$ odd integers less than $n$.
And, as noted before, there is a block matrix structure too.)

\noindent{\bf ToDo.} Check out that the conditions on $\mathbf s$ 
given in~\cite{BS05}
are necessary and sufficient to ensure that $\eta$ is nonnegative.


\subsection{Geometric isoperimetric inequalities}\label{subsec:GeomIsoper}

Lets begin with the first, historic, instance of an isoperimetric
inequality in which regular $n$-gon optimizes over all $n$-gons.
The following argument is from geometers in ancient Greece, perhaps around 200BC.

{\bf L-A Isoperimetric Result.} Amongst convex $n$-gons with given perimeter,
that which has the largest area is the regular $n$-gon.

Start with a convex $n$-gon, $\Omega_0$.
Suppose $P_{i-1}$, $P_i$, $P_{i+1}$ are consecutive vertices.
(i) show that among all isoperimetric triangles with the same base the isosceles triangles has maximum area. 
Thus by changing $\Omega_0$ moving point $P_i$ to $P_i'$ with 
$P_{i-1}$, $P_i'$, $P_{i+1}$ isosceles, the area of the changed polygon is
increased.
Apply this process to all triples of consecutive vertices.
By iteration, one finds that the optimal polygon must be equilateral: 
call it $\Omega_1$.\\
(ii) show that if the polygon $\Omega_1$  is not equiangular, 
its area may be increased by redistributing perimeter 
from a pointy to a blunt angle until the two angles are the same. 

An alternative to (ii) is to use that
Among all $n$-gons with given side lengths, the cyclic $n$-gon has
the largest area.
And a cyclic $n$-gon with equal sides is regular.
\smallskip

(For this and related isoperimetric inequalities see~\cite{AHM09}.
See also~\cite{deV11} concerning equiangular and equilateral considerations.)

\medskip

By a $\cal D$ we shall mean some domain functionals, and we are interested in pairs of these for which one has a result of the form\\
{\it For tangential $n$-gons with fixed ${\cal D}_1$ the regular $n$-gon 
$<$maximizes$|$minimizes$>$ ${\cal D}_2$}\\
Table~\ref{tbl:tbl1ar} below presents some:

\begin{table}[h]
\begin{center}
\begin{tabular}{|| c | c | c | c ||}
\hline
${\cal D}_1$& & ${\cal D}_2$& Remark \\
\hline
area & min& perimeter& \\
area & max& inradius& $A=\rho L/2$  \\
inradius& min& perimeter& Jensen inquality\\
inradius& min& area& " \\
inradius& min& $i_2=\frac{16\Sigma_\infty}{\rho}$& " \\
inradius& min& $i_4$& " \\
\hline
\end{tabular}
\caption{Tangential $n$-gons}
\label{tbl:tbl1ar}
\end{center}
\end{table}


The Jensen inequality/convexity concerns convex functions $\phi$
on some interval $I$,  and that
$$ {\rm with }\ x_i\in{I} \ {\rm and }\ \sum a_i =1\ {\rm with }
a_i\ge{0}\ ,\  \phi(\sum a_i x_i)\le \sum a_i\phi(x_i) . $$
The inequality is reversed for concave $\phi$.
By using the formula for the length for a $n$-gon given at
equation~(\ref{eq:i0Tgen}) and the convexity of $\cot(\cdot/2)$ on
$(0,\pi)$ we have the first entry of the following.
\begin{itemize}
\item We use $\cot(\alpha/2)$ convex for $\alpha\in(0,\pi)$
and (\ref{eq:i0Tgen}).
We have
$$\frac{L}{2\rho n}
= \sum\frac{1}{n}\cot(\frac{\alpha_k}{2})
\ge \cot\left( \sum \frac{1}{n}\frac{\alpha_k}{2}\right)
=\cot(\frac{(n-2)\pi}{2 n})
\frac{L_n}{2\rho n} ,
$$
where $L_n$ is the perimeter of the regular $n$-gon with inradius $\rho$.
This establishes the entry in the table corresponding to fixed $\rho$,
mimimizing the perimeter (and since $A=\rho\,L/2$ also minimiing $A$).

\item Using that
$$ \cot(\frac{\alpha}{2}) +  \frac{\cot(\frac{\alpha}{2})^3}{3}\ \ 
{\rm convex\ for\ } \alpha\in(0,\pi), $$
and the formula (\ref{eq:i2Tgen}) the Jensen inequality approach
above gives that, at fixed $\rho$, $i_2$ is minimized by the
regular polygon with the same number of sides.

\item Starting from formula (\ref{eq:i4Tgen}), 
the result on $i_4$ follows in the same way.

\end{itemize}

{\it For cyclic $n$-gons with fixed ${\cal D}_1$ the regular $n$-gon 
$<$maximizes$|$minimizes$>$ ${\cal D}_2$}\\
Table~\ref{tbl:tbl2ar} below presents some:

\begin{table}[h]
\begin{center}
\begin{tabular}{|| c | c | c | c ||}
\hline
${\cal D}_1$& & ${\cal D}_2$& Remark \\
\hline
area & min& perimeter& \\
area & max& inradius&  \\
circumradius& max& perimeter& Jensen inquality\\
circumradius& max& area& " \\
\hline
\end{tabular}
\caption{Cyclic polygons}
\label{tbl:tbl2ar}
\end{center}
\end{table}

There is a huge literature even restricting to convex sets.\\
{\small
\verb$https://math.stackexchange.com/questions/749528/isoperimetric-inequality-isodiametric-inequality-hyperplane-conjecture-what$
}

Also the `stability' of the isoperimetric inequalities is studied 
in many different ways sometimes involving `isoperimetric deficit',
sometimes `Fraenkel asymmetry'.

\medskip
If $L$ is the perimeter of a convex polygon $\Omega$, $A$ its area, 
$\rho$ the inradius,
and $s$ the length of any chord through
the centre of a largest inscribed circle, then
$$ L^2 -4\pi A \ge\frac{\pi^2}{4} (s-2\rho)^2 .$$
This is a sharpened isoperimetric inequality for convex polygons.\\
H. Hadwiger, 
{\it Comment. Math. Helv.} {\bf 16} (1944),305-309.\\
Tangential polygons can be regarded as limits of tangential $n$-gons,
so the inequality applies to them and is
$$ L (L-2\pi\rho) \ge \frac{\pi^2}{4} (s-2\rho)^2 .$$
In notation as in the list at the beginning of~\S\ref{sec:BlSa},
$$ s \ge d_O+\rho\qquad{\rm so}\ \ 
L (L-2\pi\rho) \ge \frac{\pi^2}{4} (d_O-\rho)^2 . $$
With, as in~\S\ref{subsec:rLdO}, $x=L/d_O$, $y=\rho/d_O$,
$$ x \, (x-2\pi y) \ge \frac{\pi^2}{4} (1-y)^2 ,$$
or
$$ x = \frac{1}{2} \pi  \left(\sqrt{5 y^2-2 y+1}+2 y\right) . $$
This gives a curve bit left of the upper left 
1-cap curve of Figure~\ref{fig:IIbtangQuadxy}.

\bigskip
There is a considerable literature concerning moments of inertia {\it about the centroid}.
(There are, of course, situations involving appropriate symmetries when
the centroid and incentre will coincide.)
\\
{\bf Results.} {\it The equilateral triangle minimizes the moment of inertia, 
among all convex curves with given perimeter.}\\
References include~\cite{KH87}, \cite{Ti63}, \cite{MG74}.
See also~\S\ref{sec:centroids}.

\subsection{Further geometric items}

For every convex domain
$$ \pi \rho+ \frac{|\Omega|}{\rho}
\le |\partial\Omega|
\le 2\,  \frac{|\Omega|}{\rho} .
$$
See~\cite{CFG02} equation (8).
For tangential polygons
the right hand side is an equality 
and the left-hand side is just $L\ge{2}\pi\rho$.\\
There may be some use for Bonnesen inequalities, see~\cite{Oss79},
and stronger forms for tangential polygons.

\medskip

Repeat here, for the third time(!),
 the existence statement of~\S\ref{sec:Transformations}
and of~\S\ref{sec:Construction}.
Let there be given a convex polygon, and hence its sequence of angles, all less than $\pi$.
Then one can define a tangential polygon with the same sequence of angles and
the same area.

\medskip

A tangential polygon has a larger area than any other convex polygon with the same perimeter and the same interior angles in the same sequence.
(See~\cite{AM04,Kn44,YJ19}.)\\
Amongst all convex polygons with the same area and with the same interior angles in the same sequence\\
(i) those which have the smallest perimeter are tangential polygons, and\\
(ii) those which have the largest inradius are tangential polygons.

\medskip

Amongst all tangential quadrilaterals with a given sequence of side lengths,
that which maximizes the area is bicentric.\\
One starting point is the more general result.\\
Amongst all quadrilaterals with given side lengths, that which has maximum area is cyclic.\\
(Proof: Use Bretschneider's formula.)\\
(A very easy special case is that the area of a kite with given sides is maximized by the right kite.)

There may be generalization to $2m$-gons, in particular hexagons.\\
\verb$https://mathworld.wolfram.com/CyclicHexagon.html$\\
gives area in terms of sides.

\noindent{\bf Theorem.}
{\it  For any quadrilateral with given edge lengths, 
there is a cyclic quadrilateral with the same edge lengths.}

\noindent{\bf Theorem.}
{\it The cyclic quadrilateral has the largest area of all quadrilaterals with sides of the same length.}

\subsection{Cheeger constant}

For tangential polygons $\Omega$, the Cheeger constant is
$$ h_\Omega
= \frac{|\partial\Omega|+\sqrt{4\pi|\Omega|}}{|2\Omega|}.
$$
Clearly, at fixed area, since perimeter is minimized by the regular $n$-gon,
the regular $n$-gon minimizes $h_\Omega$ over tangential $n$-gons with given area.
Much more has been established.

{\it Among all simple polygons with a given area and at most $n$ sides, the regular $n$-gon minimizes the Cheeger constant.} (See~\cite{BF}.)

{\it If $\Omega$ is a convex polygon, we denote $\Omega_*$ the (unique up to rigid motions) circumscribed polygon which has the same area as $\Omega$ and whose angles are the same as those of $\Omega$, then
$$ h(\Omega) \ge h(\Omega_*) ,$$ 
with equality if and only if $\Omega=\Omega_*$ (up to rigid motions).}

\section{Triangles}\label{sec:IIbTri}

Given the area $A$ and the angles $\alpha$, $\beta$ and $\gamma=\pi-\alpha-\beta$ of the triangle,
 one can determine the inradius and thence, if needed, the sides.
We have
$$ A = \rho^2\left(\cot\frac{\alpha}{2}+\cot\frac{\beta}{2}+\cot\frac{\gamma}{2}\right) .$$
The semiperimeter $s=(a+b+c)/2$ is found from
$$ s= \frac{A}{\rho} =  \rho\,\left(\cot\frac{\alpha}{2}+\cot\frac{\beta}{2}+\cot\frac{\gamma}{2}\right)
=  \sqrt{A\, \,\left(\cot\frac{\alpha}{2}+\cot\frac{\beta}{2}+\cot\frac{\gamma}{2}\right)}.$$
Now, defining $f$ from
$$\frac{\sin(\alpha)}{a}=\frac{\sin(\beta)}{b}=\frac{\sin(\gamma)}{c}=\frac{1}{f} ,$$
we have
$$s=\frac{a+b+c}{2}=\frac{f}{2}\, \left( \sin(\alpha)+\sin(\beta)+\sin(\gamma)\right) ,$$
which determines $f$ and hence all the sides.

\medskip
Using
$$ T_A=\frac{1}{\tan(\frac{\alpha}{2})},\quad
T_B=\frac{1}{\tan(\frac{\beta}{2})},\quad
T_C= \frac{T_A+T_B}{T_A T_B-1} ,
$$
$i_2$ and $i_4$ can be found from
equations~(\ref{eq:i2Tgen}) and~(\ref{eq:i4Tgen}).
From these $\Sigma_\infty$ and $\Sigma_1$ can be found:
see Part III~\S\ref{sec:GenTri}

\begin{figure}[ht]
\centerline{\includegraphics[height=10cm,width=14cm]{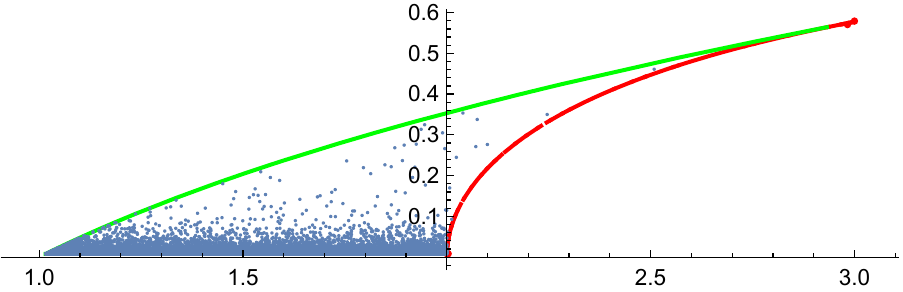}}
\caption{Triangles: $x=L/(2 d_O)$, $y=\rho/d_O$ with $\rho=1$.}
\label{fig:rLdOtri}
\end{figure}

With apologies for leaving the write-up in code, all the reasonable results are readily proved for triangles. 
The code also produces an example of a Blaschke-Santalo diagram,
shown in Figure~\ref{fig:rLdOtri},
 of the kind we will see later for other tangential polygons.

{\small
\begin{verbatim}
( * An isoperimetric result
Amongst triangles with a given inradius, taken as 1,
that which has the smallest perimeter is the equilateral triangle.
The allowed values of TA>0 and TB>0 are those for which TA*TB>1 *)
TCfn[TA_,TB_]:= (TA+TB)/(TA*TB-1);
Lfn[TA_,TB_]:= TA+TB+TCfn[TA,TB];
(* Lfn ia half the length *)
(* The hessian is positive definite, positive diagonal entries and Det>0 
on using TA*TB>1, etc. *)
hess = Map[Simplify, D[tmp, {{TA, TB}, 2}]];
Factor[Det[hess]]
(* Find minimum by checking where gradient is 0  *)
g = Map[Simplify, Grad[Lfn[TA, TB], {TA, TB}]]
Solve[{g[[1]] == 0, g[[2]] == 0}, {TA, TB}]
(* gives (TA,TB) = (Sqrt[3],Sqrt[3]) which is equilateral triangle *)

(* Define also *)
dOfn[TA_,TB_]:= Max[{TA,TB,TCfn[TA,TB]}];
(* One can also show that at given inradius, the triangle which
minimizes the distances from incentre to vertices is equilateral. *)

(* One can also imagine fixing not just the inradius but also dO
and with these TWO constraints finding 
(i) the shape with the maximum perimeter,
(ii) the shape with the minimum perimeter.
And find that they come out to be the obvious different sorts of isosceles triangles.
*)


(* Think of A as apex of triangle
Expect to see different behaviour for TA large - small apex angle
to what one gets when TA is small - big apex angle.
Let B the angle I will vary be <= Pi/2 so TB >=1.
Actually TB also gets restricted to be TB >= Max[1,1/TA]
This restiction causes TC>0 *)

(* fix TA, vary TB
 LAfn is a convex function and its first derivative is zero when TB=TC,
 i.e. triangle is isosceles 
 So, at fixed rho, the perimeter of a triangle with one angle fixed is
minimized when the triangle is isosceles having TB=TC= (1+Sqrt[1+TA^2])/TA
BCequal = Factor[TCfn[TA, TB] - TB]
Solve[BCequal == 0, TB]
Factor[D[Lfn[TA, TB], TB]/BCequal] (* clearly nonzero *)
Factor[D[Lfn[TA, TB], {TA, 2}]] (* clearly positive *)
Plot[Lfn[TA, (1 + Sqrt[1 + TA^2])/TA], {TA, 0.1, 8}]
(* convex function, minimum at TA=Sqrt[3] *)

tmp = Together[Simplify[D[Lfn[TA,(1 + Sqrt[1 + TA^2])/TA], TA]]] 
u = Sqrt[1 + TA^2];
Simplify[tmp - (-2 + u)*(1 + u)^2/(TA^2*u)] (* 0 *)
Simplify[tmp /. TA -> Sqrt[3]]. (* 0 *)

ppu =ParametricPlot[{Lfn[TA,(1+Sqrt[1+TA^2])/TA]/dOfn[TA,(1+Sqrt[1+TA^2])/TA],
 1/dOfn[TA,(1+Sqrt[1+TA^2])/TA]},{TA,0.01,Sqrt[3]},PlotStyle->Red];
ppl =ParametricPlot[{Lfn[TA,(1+Sqrt[1+TA^2])/TA]/dOfn[TA,(1+Sqrt[1+TA^2])/TA],
 1/dOfn[TA,(1+Sqrt[1+TA^2])/TA]},{TA,Sqrt[3],128},PlotStyle->Green];

pairs[AB_] := {Lfn[AB[[1]], AB[[2]]]/dOfn[AB[[1]], AB[[2]]], 
   1/dOfn[AB[[1]], AB[[2]]]};
rdm[Npts_] := 
  Module[{k} , 
   Map[pairs, 
    Table[{1 + RandomReal[{0, 32}], 1 + RandomReal[{0, 32}]}, {k, 1, 
      Npts}]]];
lp = ListPlot[rdm[20000], PlotRange -> All];
rLdOtri = Show[{ppu, ppl, lp}, PlotRange -> All]
\end{verbatim}
}

The dots in the figure are from random choices for two of the $T$ values
(corresponding to two angles).
The large number corresponding to large values of $d_O$ is from using
a uniform random distribution for the $T$ allowing large values
(so thin, or squat, triangles).
The red curve is from tall thin isosceles triangles.
The green curve is from tall thin isosceles triangles.

\clearpage

\subsection{Generalizing, angle-averaging}\label{subsec:triGenz}

We have also used Mathematica {\tt Minimize} on these tasks
and it worked well for $n=3$ and $n=4$ and {\tt NMinimize}
for larger $n$.
However, having already proved the result using Jensen's
inequality in \S\ref{subsec:GeomIsoper} one learnt (very) little from
the exercise.

However, one way to learn a little more is to use the convexity of
$\cot$ to `average two angles' as remarked upon in the
`tilting transformation' treated in~\S\ref{sec:Transformations}.
We keep $\rho$ fixed, say $\rho=1$.
Choose two vertices, with angles $\alpha_i$ and $\alpha_j$.
Let 
$$\sigma_i=\tan(\frac{\alpha_i}{4})\qquad{\rm so\ } \ \
T_i = \frac{1-\sigma_i^2}{2\sigma_i} .$$
Then, as
$$\tan(\frac{\alpha_i+\alpha_j}{4}) 
= \frac{\sigma_i + \sigma_j}{1-\sigma_i \sigma_j}) .$$
The convexity of the $\cot$ function on $(0,\pi/2)$ is
\begin{eqnarray*}
2\cot(\frac{\alpha_i+\alpha_j}{4})
&\le& \cot(\frac{\alpha_i}{2}) + \cot(\frac{\alpha_j}{2}) , \\
\frac{2(1-\sigma_i \sigma_j)}{\sigma_i + \sigma_j}
&\le& T_i + T_j .
\end{eqnarray*}
From this, and the representation of $L$ as a sum of $T_k$, we
see that $L$ is reduced by averaging two of the angles.

\medskip
The same argument works for $i_2$ and for $i_4$ on using the
convexity of $\cot^3$ and of $\cot^5$.
\medskip

It is also easy to work out the differences, by how much the quantities
$L$, $i_2$, etc. decrease.

{\small
\begin{verbatim}
(* Write T[1]=1/tan(alpha1/2) and T[2] in terms of u1 = tan(alpha1/4) and u2 resp. *)
T[1] = (1-u1^2)/(2*u1);
T[2] = (1-u2^2)/(2*u2);
(* After angle averaging/ optimal tilting when adjacent *)
tAve= (1-u1*u2)/(u1+u2);
Tout[1]= tAve;
Tout[2]= tAve;

Ldifference = Factor[2*(T[1]+T[2]-2*tAve)]
(* (((u1 - u2)^2*(1 - u1*u2))/(u1*u2*(u1 + u2))) *)

i2STfn[TT_]:= (TT+TT^3/3);
i2difference = Factor[2*(i2STfn[T[1]]+i2STfn[T[2]]-2*i2STfn[tAve])]
(* long expression - an obviously positive expression * Ldifference^3 *)
\end{verbatim}
}
\section{Tangential quadrilaterals}\label{sec:quad}

We begin with the context:

\begin{figure}[ht]
\centerline{\includegraphics[height=10cm,width=14cm]{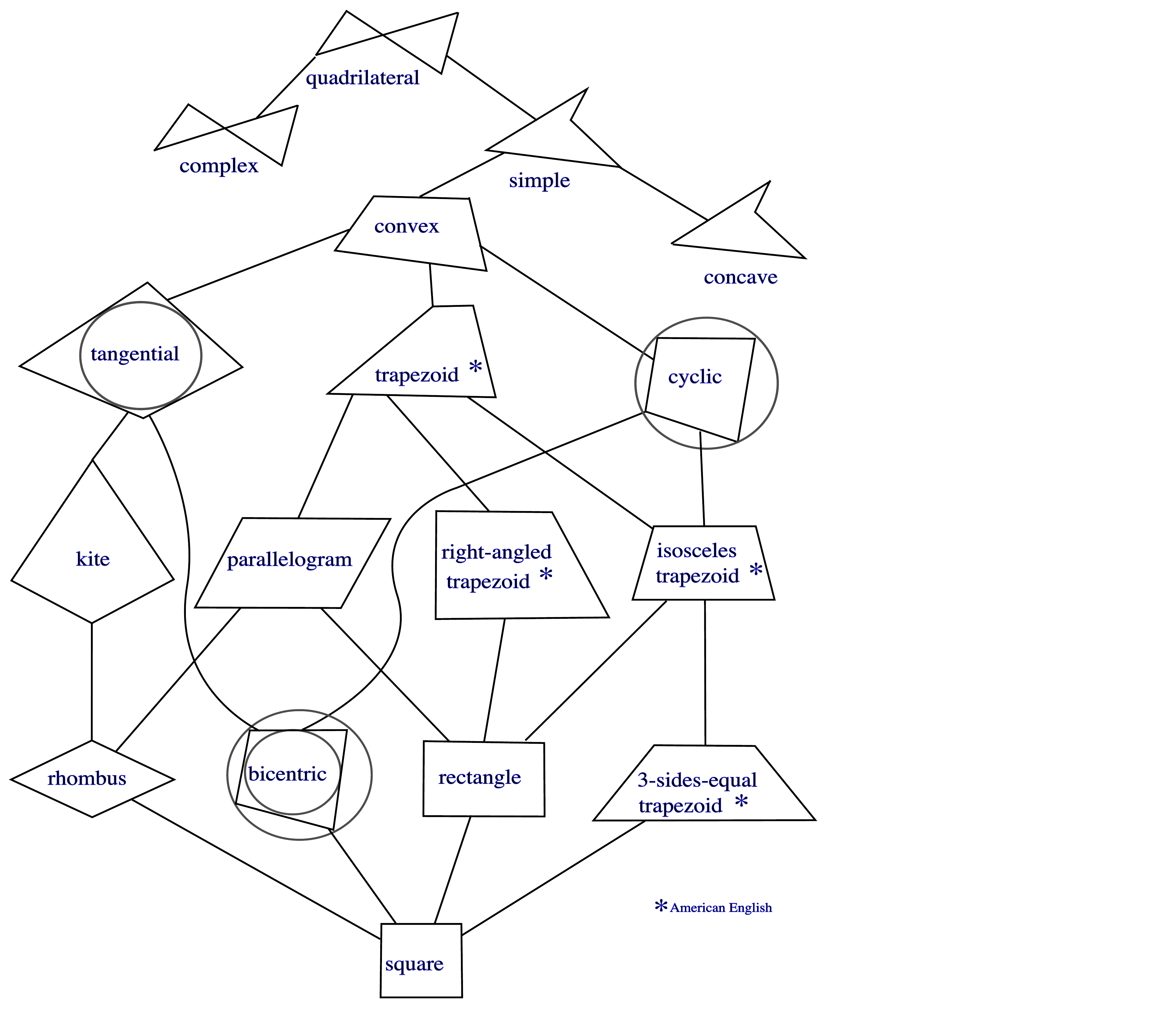}}
\caption{ wikipedia, attribution below }
\label{fig:wikQuad}
\end{figure}

Attribution for Figure~\ref{fig:wikQuad}:
By Alexgabi, jlipskoch\\
\verb$https://commons.wikimedia.org/wiki/File:Laukien_sailkapena.svg,$\\
CC BY-SA 3.0\\
\verb$https://commons.wikimedia.org/w/index.php?curid=34027107$
\clearpage

There is a huge literature on tangential quadrilaterals.
See~\cite{CS00,Gr08,Haj08,Ho11,JoBiMinArea,Jo10,Jo11,JoBiMaxArea,Mi09,Mi12}

In a tangential quadrilateral
the two diagonals and the two tangency chords are concurrent.

\noindent
{\bf Theorem.}{\it Let $ABCD$ be a tangential quadrilateral and 
$O_d$ be the point of intersection of its diagonals.
Invert, using $O_d$ as pole, each of the vertices, the inverse of $A$
denoted by $A_1$, etc.
 Let $A_1B_1C_1D_1$ be the quadrilateral obtained by these inversions.
Then$A_1B_1C_1D_1$ is a tangential quadrilateral.}

See~\cite{Mi12}.
Here is some additional  comment.
The intersection of the diagonals of the two quadrilaterals coincide.
Mobius maps map, in general, lines to circles, 
but lines through the pole are preserved.
Thus the angles between the diagonals at $O_d$ are the same (as the diagonals are).
The conformal map $z\rightarrow\frac{1}{z}$ {\it locally} preserves angles between lines,
and so does its conjugate.
I noticed some qualitative similarity between inputs and outputs.
In particular inputs like kites produced outputs like kites.
We already have that if the diagonals of the input cross at right angles
so will those of the output.

\begin{figure}[ht]
\centerline{\includegraphics[height=10cm,width=11cm]{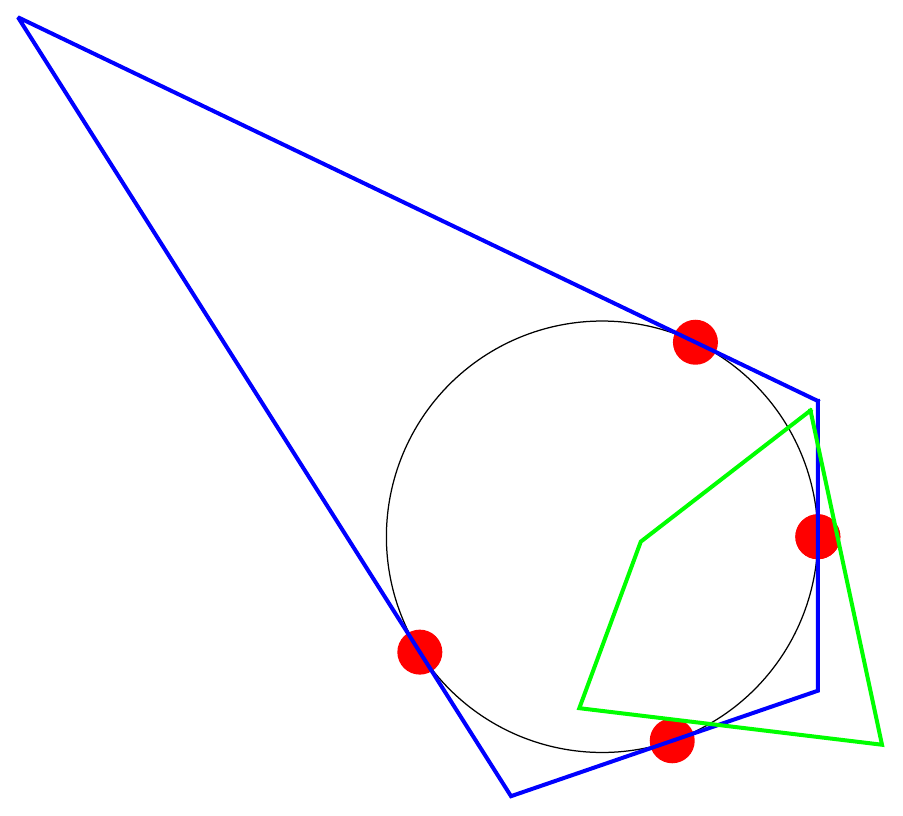}}
\caption{Inverse map. See~\cite{Mi12}.}
\label{fig:invQuad}
\end{figure}

The image under the inverse map with the incentre as pole seems 
often to take
tangential quadrilateral vertices to points
that look somewhat like vertices of a rhombus.

\clearpage

\section{Bicentric polygons}\label{sec:Bicentric}

In our calculations we often use $T_k$ as inout variables,
or equivalently $\eta_k=\rho\,T_k$.
There are additional identities for bicentrics beyond those that
one has for tangential $n$-gons.

wikipedia gives, for bicentric quadrilaterals,
$$\leqno{\rm bicentric4gon:}\qquad
\rho^2 = \eta_1\eta_3=\eta_2\eta_4 .$$
\cite{Ra05} studies bicentric hexagons and gives
$$\leqno{\rm bicentric6gon:}
\rho^2 = \eta_1\eta_3 +\eta_3\eta_5 + \eta_5\eta_1
=\eta_2\eta_4 +\eta_4\eta_6+\eta_6\eta_2.$$
There has been use of these in subsequent calculations.

The challenge put to me, admittedly in connection with $\QQ_{0-}$
rather than the geometric functionals treated in this Part IIb,
is: investigate Blaschke-Santalo diagrams for tangential polygons.
One triple that can be so investigated is $(\rho,L,R_V)$ for 
bicentric polygons.
For triangles, see the later section on Blundon's inequality,
e.g inequalities~(\ref{in:Blund}).

There is a huge literature on bicentric polygons dating back to the 19th century.
The topics include Fuss's Theorem, Poncelet’s porism and more.
Some of this is impressive so is presented below - with nothing original of mine -
but with the hope that some might be of use in future attempts to
produce Blaschke-Santalo diagrams.
\medskip

\subsection{Fuss's Theorem, Poncelet's Porism}

We would like to establish (polynomial) relations involving $(\rho,L,R_V)$.
Fuss's theorem(s) involves $(\rho,R_V,d)$ where
$d$ is the distance between the circumcentre and incentre of the bicentric polygon.
For triangles,
The quantity $d$ can be recognized in Blundon's inequality~(\ref{in:Blund})
and is, as in the first displayed equation below, $d=\sqrt{R_V(R_V-2\rho)}$.

The following is a (slightly adapted) quote from\\
\verb$https://mathworld.wolfram.com/PonceletsPorism.html$

\begin{quote}
The three numbers $(\rho,R_V,d)$ will not be arbitrary and along with $n$, 
they will have to satisfy certain relations. 
For the case of a triangle, one such relation is sometimes called the Euler triangle formula: 
$$ R_V^2 - 2 R_V\rho - d^2 = 0.  $$
One of popular notations for such relations 
(which is necessary and sufficient for existence of a bicentric polygon) 
can be given in terms of the quantities
$$ a=\frac{1}{R_V+d},\ \  b=\frac{1}{R_V-d},\ \  c=\frac{1}{\rho} .$$
For a triangle, the Euler formula has the form: 
$$a + b = c, $$
for a bicentric quadrilateral
$$a^2 + b^2 = c^2 .$$
The relationship for a bicentric pentagon is 
$$ 4(a^3 +b^3 +c^3) = (a+b+c)^3 .$$
Let
\begin{eqnarray*}
E_1 &=& -a^2+b^2+c^2 ,\\
E_2 &=& a^2-b^2+c^2 ,\\
E_3 &=& a^2+b^2-c^2 .
\end{eqnarray*}
The relationship for a bicentric hexagon is
$$ \frac{1}{E_1} + \frac{1}{E_2} = \frac{1}{E_3} 
$$
\end{quote}

\subsection{Triangles, again}

The sides $s_1$, $s_2$, $s_3$ of a triangle are the roots of the cubic
\begin{equation}
y^3 - L y^2+(\frac{L^2}{4}+\rho^2+4\rho\,R_V)y  - 2 L\rho R_V=0.
\label{eq:triCubic}
\end{equation}

\noindent{\bf ToDo.} Find conditions on the coefficients of this cubic
in order that it have 3 positive roots with the largest root less than
the sum of the other two.
Sturm sequences might be useful.
Do Blundon's inequalities come from this?

\subsection{Bicentric quadrilaterals}

The wikipedia article `Bicentric quadrilateral' contains many items.
There is a quartic equation with coefficients in terms of $L$, $\rho$ and $R_V$
whose solutions are the sides of a bicentric polygon:
\begin{equation}
y^{4} - L y^{3}+(\frac{L^{2}}{4}+2\rho^{2}+
2\rho {\sqrt{4R_V^{2}+\rho^{2}}})y^{2} - 
\rho L({\sqrt{4R_V^{2}+\rho ^{2}}}+\rho ) y+\rho ^{2}\frac{L^{2}}{4}=0 .
\label{eq:biQuadQuartic}
\end{equation}

\noindent{\bf ToDo.} Find conditions on the coefficients of this quartic
in order that it have 3 positive roots with the largest root less than
the sum of the other three.
Sturm sequences might be useful.
\bigskip

There are huge number of inequalities.

\noindent{\bf Theorem.}{\it
If a bicentric quadrilateral has an incircle and a 
vertex-circumcircle with radii $\rho$ and $R_v$ respectively, then its area $A$ satisfies
$$ \frac{1}{2}\rho\, L=A\ge 2\rho\sqrt{2\rho(\sqrt{4R_v^2+\rho^2}-\rho)} ,$$
where equality holds if, and only if, the quadrilateral is also an isosceles
trapezium.}\\
See~\cite{JoBiMinArea}.\\
For a square $\rho=1$, $R=R_v=\sqrt{2}$, $A=4$ gives equality in the
preceding inquality.
\medskip

For a bicentric quadrilateral
$$\frac{1}{32}(L^2-16A) \le R_v^2-\rho^2$$
with equality if and only if the hexagon is regular.\\

\subsection{Bicentric hexagons}

\noindent
For a bicentric hexagon
$$\frac{1}{36}(L^2-8\sqrt{3}\,A) \le R_V^2-\frac{4}{3}\rho^2$$
with equality if and only if the hexagon is regular.

For bicentric $n$-gons see~\cite{Mao96}.
\subsection{$R=R_v=d_O$ for regular polygons}

In a regular $n$-gon, the side $s_n=2\rho\tau_n$ where $\tau_n=\tan(\pi/n)$.
The circumradius $R$ and $d_O$ coincide.
We have
$$ 4(R^2-\rho^2) = s_n^2$$
so
$$ R=d_O = \sqrt{\rho^2+\frac{s_n^2}{4}}=\rho\sqrt{1+\tau_n^2} . $$
See Part IIa~\S\ref{sec:regn} for the formulae for $L=i_0$, $i_2$ and $i_4$.

\subsection{Bicentrics from regular}

Given any regular $n$-gon any choice of 3 of its vertices gives a
bicentric polygon as all triangles are bicentric.

Given a regular 8-gon, selection of its alternate vertices gives a 
bicentric 4-gon, a square.

Given a regular 9-gon, there is a selection of 5 of its vertices giving
a bicentric 5-gon. See~\cite{To07,To15}.
\section{Blaschke-Santalo diagrams for tangential polygons}\label{sec:BlSa}

\subsection{Definitions for Blaschke-Santalo diagrams}

Where we feel it helps the exposition
there will be some repetition of material already presented in PartIIa.
We always choose the origin of our coordinate system to be at the incentre.
The inradius $\rho$  and perimeter $\PP$ will occur in our diagrams,
and there will be various choices of a third geometric quantity
which, for the present, we denote by $q$.
In our Blaschke-Santalo diagrams the horizontal axis is usually $x=\PP/q$
and the vertical axis $y=\rho/q$.
(Further work in which $x=\PP/\rho$ and, as before, $y=\rho/q$
 might be undertaken, motivated by having the same $x$ for the
different $q$, e.g. $q$ being $i_2^{1/3}$, $Q_{0-}^{1/4}$, etc.)
The area of a tangential polygon is $A=\rho\,\PP/2$.
The quantity $x/y=\PP/\rho$ is denoted by $B$ in~\cite{PoS51} .
We have
 $B\ge{2\pi}$ with equality only for the disk,
and for any triangle, $B\ge{6\sqrt{3}}$):
see~\cite{Ai58}.
The perimeter and inradius of a regular $n$-gon are related by
$$\rho = \frac{\PP}{2n}\cot(\frac{\pi}{n}) , $$
and for a tangential $n$-gon 
$$B\ge{2n\tan(\pi/n)} . $$

The other domain functionals are as follows.\\
(1) There is a first set of geometric quantities:\\
$\bullet$ circumradius $R$ (radius of the smallest disk containing the region);\\
$\bullet$ the distance from the incentre to the boundary $d_0$;\\
(2) There are moments about the incentre:\\
$\bullet$ the second boundary moment $i_2$;\\
$\bullet$ the fourth boundary moment $i_4$;\\
(3) There are quantities leading to $\QQ_{0,-}$:\\
$\bullet$ $\Sigma_\infty=\rho\, i_2/16$;\\
$\bullet$ $0>\Sigma_1=(i_2^2/\PP -i_4)/16$;\\
$\bullet$ $\QQ_{0-}$ the lower bound on torsional rigidity treated in
Part I, i.e.~\cite{Ke20tIMA}.\\
We defer discussion of these last items, $\QQ_{0,-}$, until Part III.

\medskip

In connection with bicentric polygons, including triangles, the radius of the circle 
through the vertices is denoted $R_V$, and we always have $R\le{R_V}$.

\medskip

An outline of the remainder of this section follows.

\begin{itemize}
\item In~\S\ref{subsec:cap} we define cap domains which enter the account of
the $(\rho,\PP,R)$ diagram in the next subsection.

\item In~\S\ref{subsec:dOngon} we present formulae for $d_O$ for various
tangential $n$-gons.

\item In~\S\ref{subsec:rLR} we begin by noting that the account for general convex domains in~\cite{BCS03}
has on one of its boundaries the 2-cap.
We have done some computations, and have a belief that regular $n$-gons occur 
on one of the boundaries.
However, the circumradius $R$ doesn't seem to be a good lead-in for
computations related to ($i_2$, $i_4$ and) $\QQ_{0-}$.

\item \S\ref{subsec:rLdO} is the main subsection.
The quantity $d_O$ is a function easily defined in terms of the $T_k$
occuring in expressions for $L$, $i_2$ and $i_4$, and hence in $\QQ_{0-}$.
The expressions for $L$, $i_2$ and $i_4$ are given in Part IIa.
The only Blaschke-Santalo diagrams presented here are those for
$(\rho,L,d_O)$.
Future work may involve $(\rho,L,i_2)$ and $(\rho,L,i_2)$,
but for now there is just occasional comment on relevant inequalities
for particular tangential $n$-gons.

\item Though we believe it to be an aside to our main purpose of
investigating functionals related to $\QQ_{0-}$ over all
tangential polygons, in~\S\ref{subsec:rLRV} we 
propose to treat, for bicentric polygons the triple $(\rho,L,R_V)$.
So far, the main result is just that already in the literature for
triangles.
Some items in~\S\ref{sec:Bicentric} may be useful in future efforts.

\end{itemize}

\subsection{Some extremals amongst circumgons: the 1-cap and symmetric 2-cap}\label{subsec:cap}

Because our long-term motivation concerns $\QQ_{0-}$
and this involves $i_2$ (and $i_4$) se begin with the following.

\par\noindent
{\it Let $\rho>0$.
Among all plane convex sets $\Omega$ with fixed positive area $A$ (with $A>\pi\rho^2$) that contain a disk of radius $\rho$ around the origin,
$i_2$ is maximized if $\Omega$ is the single-cap, the convex hull or a disk of radius $\rho$ and a point
(the point being, up to rotation invariance, uniquely defined by the area constraint).}\\
The proof is elementary. See~\cite{HS20} where it is a Lemma needed prior to
establishing gradient bounds for the torsion function.
Exactly the same proof gives the corresponding result for $i_4$.
\medskip

Consider tangential polygons with $\rho$ and $d_O>\rho$ fixed.
Amongst these the 1-cap\\
$\bullet$ minimizes $A$, $L$, $i_2$ and $i_4$ and \\
$\bullet$ minimizes $\QQ_0$.\\
This is because of domain monotonicity of these functionals.\\
{\bf Question.} {\it  Given two tangential polygons $\Omega_0$
and $\Omega_1$ (with the same incentre? and)
 with $\Omega_0\subset\Omega_1$, is
$\QQ_{0-}(\Omega_0)\le\QQ_{0-}(\Omega_1)$?}

\begin{figure}[h]
$${\includegraphics[height=5cm,width=5cm]{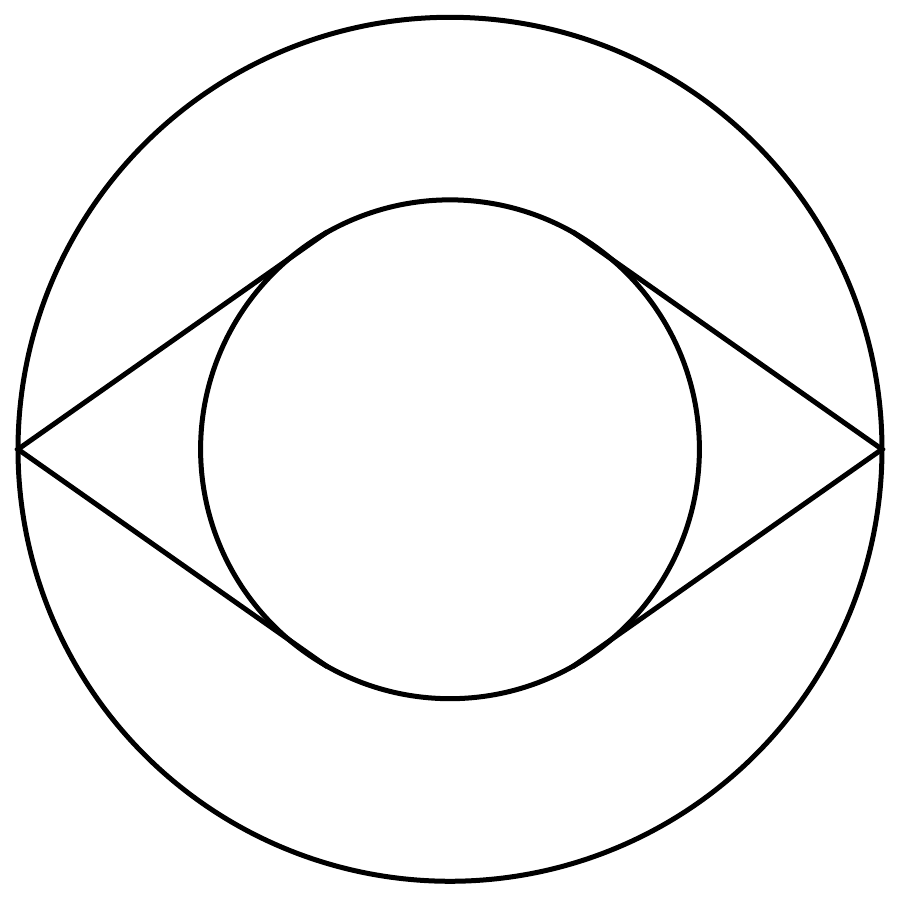}}
\qquad\qquad
{\includegraphics[height=5cm,width=5cm]{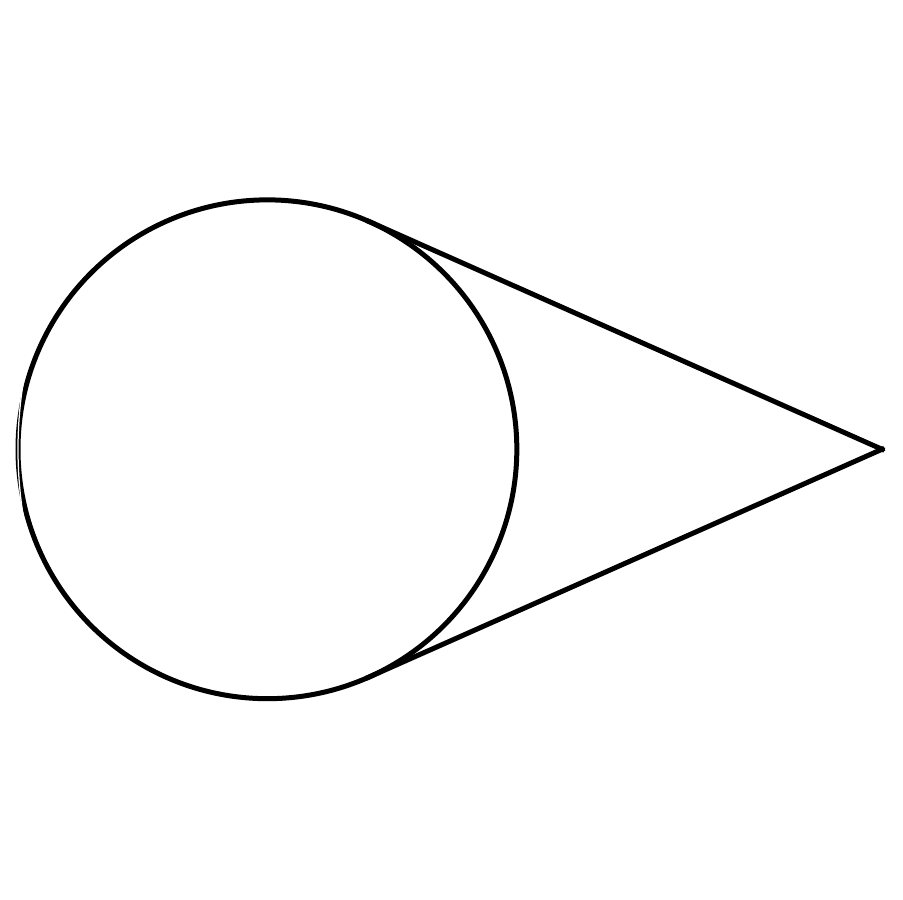}
}$$
\caption{Left: the symmetric 2--cap.
Right: the 1-cap. } 
\label{fig:tsquareplbNe0}
\end{figure}

We can easily calculate the perimeter of the 1-cap.
Denote by $\alpha$ the angle at its vertex on the circle radius $d_O$.
Then $\sin(\alpha/2)=\rho/d_O$ and
\begin{eqnarray*}
L
&=& (\pi+\alpha)\rho + 2\sqrt{d_O^2-\rho^2} ,\\
\frac{L}{d_O}
&=& \left( \pi +2\arcsin(\frac{\rho}{d_O})\right)\,\frac{\rho}{d_O} +
\sqrt{1-(\frac{\rho}{d_O})^2} .
\end{eqnarray*}
We will see this as occuring as the lower right boundary in
diagrams with $x=L/d_O$, $y=\rho/d_O$.
The limiting cases are
\begin{itemize}
\item $d_O\rightarrow\rho$, $y\rightarrow{1}$, $x\rightarrow{2\pi}$
corresponding to a disk;
\item $d_O\rightarrow\infty$, $y\rightarrow{0}$, $x\rightarrow{2}$
corresponding to a long flat shapes like the 1-cap.
\item $d_O\rightarrow\infty$, $y\rightarrow{0}$, $x\rightarrow{4}$
corresponding to a long flat shapes like the 2-cap.\\
\item The quadrilateral might have one side tending to zero,
say as a symmetric tangential trapezium tending to an isosceles triangle.
When the squat isosceles triangle becomes very thin,
$y\rightarrow{0}$, $x\rightarrow{4}$.
\end{itemize}
\bigskip

We will see the 2-cap in~\S\ref{subsec:rLR}.

\clearpage

\subsection{$d_O$ for various tangential $n$-gons}\label{subsec:dOngon}

The distances from the incentre of the vertices of a tangential
$n$-gon are given by 
$$\sqrt{\rho^2 +\eta_k^2} = \rho\sqrt{1+T_k^2} . $$

\subsubsection{$R$, $d_O$ for regular polygons}

In a regular $n$-gon, the side $s_n=2\rho\tau_n$ where $\tau_n=\tan(\pi/n)$.
The circumradius $R$ and $d_O$ coincide.
We have
$$ 4(R^2-\rho^2) = s_n^2$$ 
so
$$ R=d_O = \sqrt{\rho^2+\frac{s_n^2}{4}}=\rho\sqrt{1+\tau_n^2} . $$
See Part IIa~\S\ref{sec:regn} for the formulae for $L=i_0$, $i_2$ and $i_4$.

\subsubsection{$R$, $d_O$ for triangles}

For an equilateral triangle
$$ \PP= 6\sqrt{3}\rho, \qquad R = d_O = 2\rho . $$
In general $R$ and $d_O$ differ.
The circumradius of an acute angled triangle is the radius of the circle through the three vertices ($R=R_V$).
For a triangle with one interior angle measuring more than $\pi/2$, an obtuse triangle,
the circumradius is half the length of the longest side
(and $R<R_V$).

For triangles equation~(\ref{eqBS:elPec}) is
$$ 1-\frac{1}{T_1 T_2} -\frac{1}{T_2 T_3}  -\frac{1}{T_3 T_1} = 0.
$$

As in Part IIa, most of our calculations have, to date, concentrated on
isosceles triangles.
See Part IIa~\S\ref{sec:isos} for the formulae for $L=i_0$, $i_2$ and $i_4$.
\medskip
We return to triangles in~\S\ref{subsubsec:BStriangles}.
For now it is sufficient to record that
$$\leqno{\rm triangle:}\qquad
d_O = \rho\sqrt{1 +{\rm max}(T_1,T_2,T_3)^2} . $$

\subsubsection{$d_O$ for tangential quadrilaterals}

For tangential quadrilaterals equation~(\ref{eqBS:elPes}) is
 $$ \frac{1}{T_1} +\frac{1}{T_2}  +\frac{1}{T_3} +\frac{1}{T_4}-
 \frac{1}{T_1 T_2 T_3}-\frac{1}{T_1 T_2 T_4} - \frac{1}{ T_1 T_3 T_4 }-\frac{1}{T_2 T_3 T_4} = 0.
 $$

We return to tangential quadrilaterals in~\S\ref{subsubsec:BStangQuad}.
For now it is sufficient to record that
$$\leqno{\rm tang\ quad:}\qquad
d_O = \rho\sqrt{1 +{\rm max}(T_1,T_2,T_3,T_4)^2} . $$

\bigskip


\subsection{$(\rho,\PP,R)$}\label{subsec:rLR}

For any plane convex set
$$\rho\le{R},\ \PP\le{2\pi\,R},\ 2\pi\rho\le{\PP},\ 4R\le{\PP} . $$
There is some discussion of how to compute $R$ at\\
{\small
\verb$https://mathematica.stackexchange.com/questions/121987/how-to-find-the-incircle-and-circumcircle-for-an-irregular-polygon$
}

We can scale our shapes so there is no loss of generality in
setting $\rho=1$.
Define
$$ y = \frac{1}{R}\ \ {\rm and}\ \ x=\frac{L}{R} , $$
(with our $x$ a factor of $2\pi$ greater than than in~\cite{BCS03}).
$$ x=2\pi,\ y=1\qquad{\rm for\ a\ disk}. $$
\cite{BCS03}  establish, for convex sets,
$$4\left(\sqrt{1-y^2} +\ y\,{{\rm arcsin}(y)} \right)
\le x=\frac{\PP}{R} \le4\left(\sqrt{1-y^2} +{\rm arcsin}(y) \right) .
$$
The lower bound at the left corresponds to values for a 
symmetric 2-cap, which is a circumgon tangential polygon.
However, the right-hand, upper bound is for a convex shape which is
{\it not} a tangential polygon.
Our computations suggest that amongst the tangential polygon
shapes which will occur on the right-hand upper bound are the
regular polygons.

\begin{figure}[ht]
\centerline{\includegraphics[height=10cm,width=14cm]{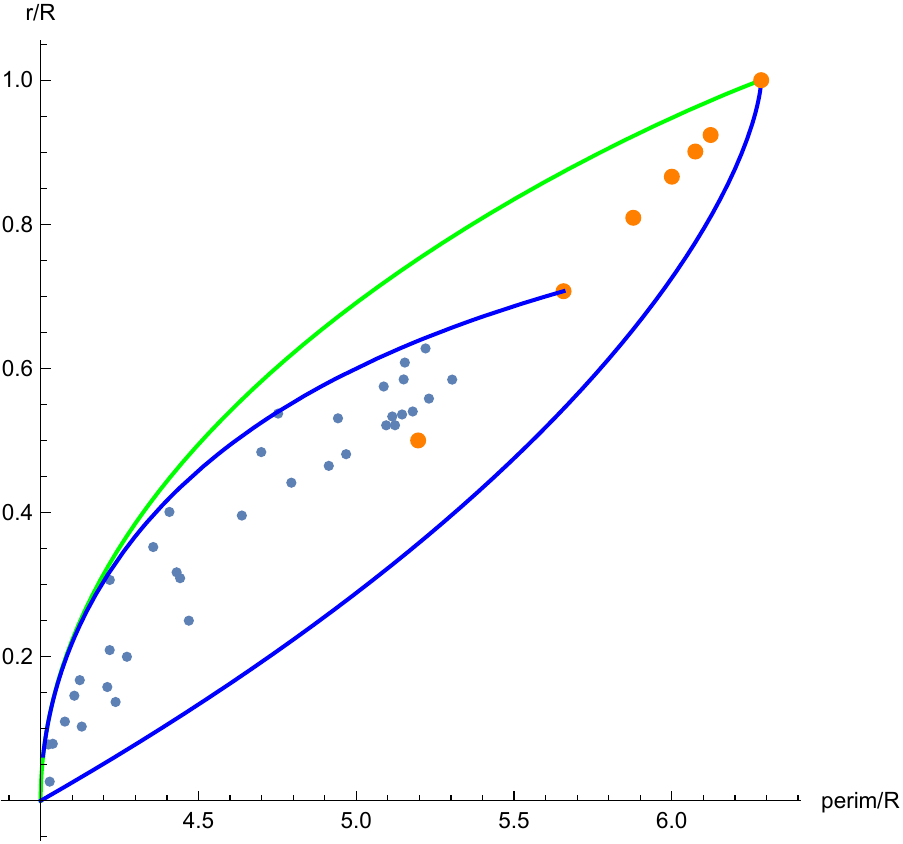}}
\caption{Except for the lower curve, some results on
$(L/R,\rho/R)$ pairs for tangential polygons.
See description in text.}
\label{fig:IIbrLR}
\end{figure}

In Figure~\ref{fig:IIbrLR} the right-most blue curve is the
shape from~\cite{BCS03} which is {\it not} a tangential polygon.
The upper-left green curve is form the symmetric 2-cap.
The red dots are some regular polygons, the rightmost upper dot the
circular disk.
The lowest is the equilateral triangle,
and the next up, joining the upper blue curve is the square.
The scatter of blue dots are from tangential quadrilaterals,
and the upper blue line from rhombi.

\clearpage

\subsection{$(\rho,\PP,d_O)$}\label{subsec:rLdO}


\subsubsection{Triangles}\label{subsubsec:BStriangles}

See also PartIIa \S\ref{sec:isos}.
\smallskip

\noindent
{\it Amongst all isosceles triangles with given inradius
\begin{itemize}
\item the equilateral triangle minimizes $d_O$,
\item the equilateral triangle minimizes $A$ and $\PP=2A/\rho$ and $\PP^2/A$, 
and
\item the equilateral triangle minimizes $i_2$ and $i_2/(A\,\PP)$.
\end{itemize}}

\medskip
\noindent
{\it At given $\rho$ and $A$ greater than the area of the equilateral
triangle of the same inradius,
amongst isosceles triangles
\begin{itemize}
\item
$d_O$ is maximised by tall isosceles triangles
(apex angle $\alpha<\pi/3$ small, $\sigma<2-\sqrt{3}$ small)

\item
$d_O$ is minimised by short squat isosceles triangles
(apex angle $\alpha>\pi/3$ near $\pi$, $\sigma>2-\sqrt{3}$ near 1).
\end{itemize}}
(Above may be true for all triangles.)

\medskip

The calculations for $i_2$ for the following are yet to be done.\\
{\bf Conjecture}
{\it At given $\rho$ and $A$ greater than the area of the equilateral
triangle of the same inradius,
amongst isosceles triangles
\begin{itemize}
\item
$i_2$ is maximised by tall isosceles triangles
(apex angle $\alpha<\pi/3$ small, $\sigma<2-\sqrt{3}$ small)

\item
$d_O$ is minimised by short squat isosceles triangles
(apex angle $\alpha>\pi/3$ near $\pi$, $\sigma>2-\sqrt{3}$ near 1).
\end{itemize}}

\subsubsection{Tangential quadrilaterals}\label{subsubsec:BStangQuad}

\noindent{\bf Rhombi}

\noindent
{\it Amongst all rhombi with given side,
\begin{itemize}
\item the square minimizes $d_O$,
\item the square maximizes $A$ and $\rho=2A/\PP$, and
\item the square minimizes $i_2/(A\,\PP)$.
\end{itemize}}

\noindent
These are easily established, as follows.
$$ A = 2\rho^2 (T+\frac{1}{T}) \qquad
\PP= 4\rho (T+\frac{1}{T}) . $$
Hence
$$\frac{\PP^2}{A} = 8 (T+\frac{1}{T}) , $$
which is minimized when $T=1$, the square.\\
Next
\begin{eqnarray*}
i_2 
&=& 2\rho^3\left( T+\frac{1}{T}+\frac{1}{3}(T^3+\frac{1}{T^3})\right) ,\\
&=& \frac{2\rho^3}{3} (T+\frac{1}{T})^3 , 
\end{eqnarray*}
so that
$$\frac{i_2}{A \PP}=\frac{1}{12}(T+\frac{1}{T}) , $$
which is minimized at $T=1$.
\medskip

We also have the following.

\noindent
{\it Amongst all rhombi with given inradius
\begin{itemize}
\item the square minimizes $d_O$,
\item the square minimizes $A$ and $\PP=2A/\rho$ and $\PP^2/A$, and
\item the square minimizes $i_2$ and $i_2/(A\,\PP)$.
\end{itemize}}
\medskip

\begin{figure}[ht]
\centerline{\includegraphics[height=10cm,width=14cm]{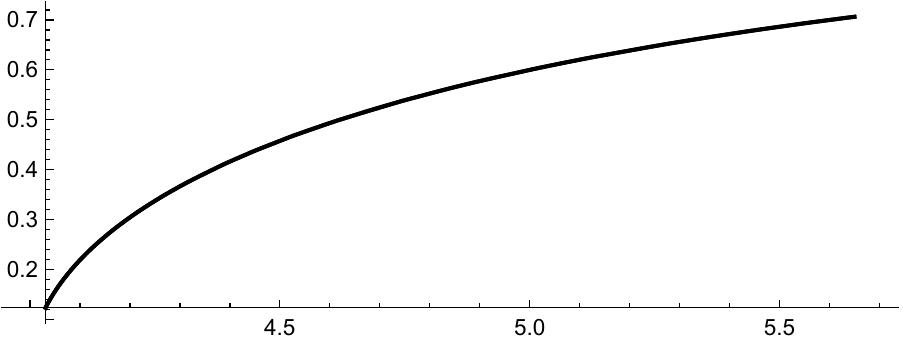}}
\caption{For rhombi, plot of $y=\rho/d_O$ against $x=L/d_O$.
The top right corner is $(1/\sqrt{2},4\sqrt{2})$ corresponding to a 
square.}
\label{fig:IIbRhombxy}
\end{figure}

A plot of $y=\rho/d_O$ against $x=\PP/d_O$ is shown in
Figure~\ref{fig:IIbRhombxy}

\clearpage

\noindent{\bf Kites}

\noindent
{\it Amongst kites with given area $A$ and distance $d$ between apexes,
the rhombus has
\begin{itemize}
\item the smallest perimeter $\PP$,
\item the largest inradius $\rho=2A/\PP$, and
\item the smallest $4I_2=\rho\,i_2$.
\end{itemize}}
A proof uses Steiner symmetrisation.
\medskip

One can easily recover results which, in a more general form, are
given in the next subsubsection.

\noindent
{\it Amongst kites with given sides $s_1$ and $s_2$ the right kite has the
largest area.}
\noindent
For ease of writing, suppose $s_1<s_2$.
The rhombus case can be treated separately.
Let $h$ be the distance between the vertices adjacent to unequal sides.
Let $\alpha_1$ be the angle at the vertex on both sides of length $s_1$.
Denote by $\beta$, the repeated angle, 
the angle at vertices adjacent to the unequal sides.
$$A =\frac{1}{2}d h\ 
{\rm which\ can\ be\ written\ }\ A = s_1 s_2 \sin(\beta). $$
The formula at the right can be deduced from that on the left
using the trigonometry:
$$  h = 2 s_1\sin(\alpha_1/2), \qquad 
d= s_2 \frac{\sin(\beta)}{\sin(\alpha_1/2)}, $$
the latter from the triangle sine rule.
Thus since $A= s_1 s_2 \sin(\beta)$, the area is maximized at $\beta=\pi/2$.

We remark 
\begin{eqnarray*}
\rho 
&=& \frac{ s_1 s_2\sin(\beta)}{s_1+s_2} ,\\
\frac{L}{\rho}
&=& \frac{2(s_1+s_2)^2}{s_1 s_2\sin(\beta)} .
\end{eqnarray*}
Also, on using $s_1<s_2$,
$$\frac{d_O}{\rho}
=\sqrt{1+{\rm max}(\cot(\frac{\beta}{2}),\cot(\frac{\alpha_1}{2}) )^2} ,
$$
and the sine rule for triangles gives $\alpha_1$ in terms
of $s_1$, $s_2$ and $\beta$.

\noindent{\bf Other tangential quadrilaterals}

For a quadrilateral with angles $\alpha_k$ and sequence of
sides $[s_1, s_2, L/2-s_1,L/2-s_2]$, the area is
$$ A = \sqrt{s_1\, s_2\, (L/2-s_1)\,(L/2-s_2)}
\sin(\frac{\alpha_1+\alpha_3}{2}) . $$
For any quadrilateral, since the sum of the angles is $2\pi$,
$$ \sin(\frac{\alpha_1+\alpha_3}{2}) 
= \sin(\frac{\alpha_2+\alpha_4}{2}) ) .$$
The formula for the area $A$ above gives the following.
\smallskip

\noindent
{\it Amongst all tangential quadrilaterals with given sequence of sides,
the bicentric quadrilateral maximizes $A$ and $\rho=2A/\PP$.\\
A kite is bicentric iff it is a right kite.}

\medskip
 
In general, amongst all kites with a given pair of sides it is {\it not} the case that
$i_2/(A\,\PP)$ is minimized by a right kite.
(Let $\beta$ be the repeated angle in the kite.
The quantity $i_2/(A\,\PP)$ plotted as a function of $1/\tan(\beta/2)$ appears to have a unique minimum
which for a rhombus is at $\beta=\pi/2$, but in general is not.)

\noindent{\bf More on quadrilaterals}

\cite{KKLY17} treat quadrilaterals with diagonals intersecting
at right angles.
If the side lengths are denoted $s_1$, $s_2$, $s_3$, $s_4$
these satisfy
$$ s_1^2 + s_3^2 = s_2^2 + s_4^2 .$$
If a quadrilateral is tangential and has its diagonals intersecting
at right angles it is  kite.
(See also~\cite{KK20}

\newpage

\noindent{\bf Bicentric quadrilaterals, $(\rho,\PP,d_O)$}

\begin{figure}[ht]
\centerline{\includegraphics[height=10cm,width=14cm]{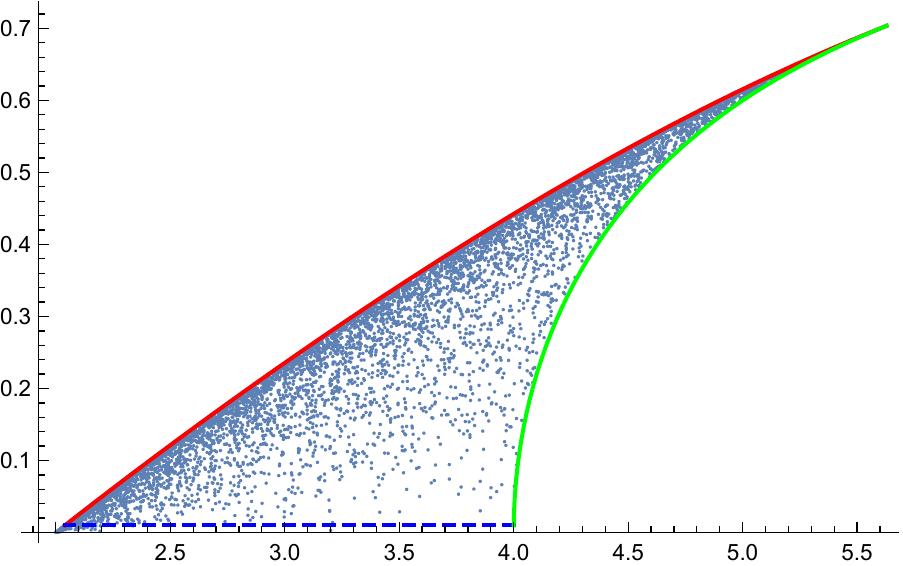}}
\caption{$x=L/d_O$, $y=\rho/d_O$ for bicentric quadrilaterals}
\label{fig:plBicenQuad}
\end{figure}

The upper left boundary curve in Figure~\ref{fig:plBicenQuad}
(red) corresponds to right kites: 
the lower right curve (green) 
corresponds to tangential isosceles trapeziums.
The dots are just from randomly generated tangential quadrilaterals.
The uppermost point $(4\sqrt{2},1/\sqrt{2})$ corresponds to a square.
The bottom boundary is approached by long thin shapes,
and the left-most point is $(2,0)$ corresponding to the limit
of thin shapes where $L\sim{2d_O}$.
\medskip

Another approach to the boundary curves involves working with
the tangent lengths $\eta_j$.
For a bicentric quadrilateral with $\rho=1$ we have
$$ \eta_1\eta_3=\rho^2=1= \eta_2\eta_4\ {\rm so}\ \ 
L= \eta_1+\frac{1}{\eta_1}+\eta_2+\frac{1}{\eta_2} .$$
Choose $\eta_1$ to be the largest of the $\eta_k$, so $\eta_1>1$
and $d_O=\sqrt{1+\eta_1^2}$.
Now $\eta_2$ must be less than or equal to $\eta_1$ and
greater than or equal to $1/\eta_1$, and consequently
So 
$$ L_-=2+\eta_1+\frac{1}{\eta_1}\le L 
\le 2(\eta_1+\frac{1}{\eta_1}) = L_+ .$$
Plotting $y=1/d_O$ against $x=L_-/d_O$ gives one of the boundary curves,
and the other is from $y=1/d_O$ against $x=L_+/d_O$.

As an aside we remark that for a bicentric quadrilateral
for which $\rho=1$ and $s_1$ ia the largest side (so $s_1\ge{2}$)
and $s_2$ is the larger of the other pair (so $s_2>1$),
$$\PP = \frac{s_1 s_2(s_1+s_2+\sqrt{4+(s1-s2)^2}}{s_1 s_2-1} .$$

\clearpage
\noindent{\bf General tangential quadrilaterals}

\begin{figure}[ht]
\centerline{\includegraphics[height=10cm,width=14cm]{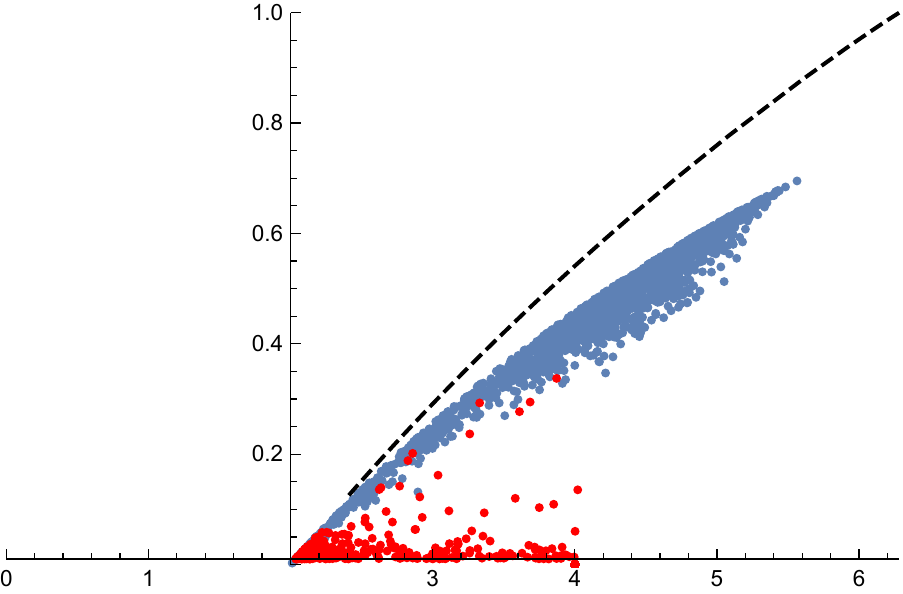}}
\caption{$x=L/d_O$, $y=\rho/d_O$ for tangential quadrilaterals}
\label{fig:IIbtangQuadxy}
\end{figure}

In Figure~\ref{fig:IIbtangQuadxy} the scatter of (blue) points
towards the left upper are from $T_k$ values with the first 3
not all that far from 1.
(Those blue dots that are near the bottom come from those with $T_4$ large,
such as occurs if the first 3 of the $T_k$ were near $1/\sqrt{3}$,
the $\alpha_k$ near $2\pi/3$.)
The red dots are from choosing the first 3 $T_k$ randomly in
the interval $(0,1000)$.
Large values of $T_k$ cause $d_O$ to be large, hence the cluster
close to $y=0$.
We haven't definitively defined the boundaries.
However we suspect the upper left boundary, at least for small enough $y$
($d_O$ large) have the minimum $L$ shapes somewhat like the 1-cap,
perhaps kites with a small apex angle at $d_O$.

The dashed curve is the bound from the 1-cap, discussed nearthe beginning
of this section.
\clearpage
 
\subsubsection{Tangential pentagons}
\vspace{1cm}
\noindent{\bf Bicentric pentagons}
\par\noindent
An elaborate formula for the area of a cyclic pentagon in terms of side lengths is given
in~\cite{Ro95}.
\vspace{0.5cm}

\subsubsection{Tangential hexagons}
\vspace{1cm}
\noindent{\bf Bicentric hexagons}
\par\noindent
A question that arises is, what conditions other than
$$s_1+s_3+s_5 = s_2+s_4+s_6 =\frac{L}{2} ,$$
are satisfied by the sides, or tangent lengths $\eta_j$, for
a bicentric hexagon. 
We have, for tangential hexagons,
$$ s_1= \eta_1+\eta_2,\ s_2= \eta_2+\eta_3,\ s_3 = \eta_3+\eta_4, $$
$$ s_4= \eta_4+\eta_5,\ s_5= \eta_5+\eta_6,\ s_6= \eta_6+\eta_1,$$
or, in the notation of~\S\ref{subsec:Circulant},
$$ M(6)\, {\mathbf \eta} = {\mathbf s} .$$
From~\cite{RK09} equations~(2.29),~(2.30)
$$\rho \sqrt{\frac{\eta_1 \eta_3 \eta_5}{\eta_1+\eta_3+\eta_5}}
=\eta_2 \eta_5 = \eta_1 \eta_4 = \eta_3 \eta_6 , $$
and, from~\cite{RK09} equation~(2.26)
$$ \eta_1 \eta_3 + \eta_3 \eta_5 + \eta_5 \eta_1
=\eta_2 \eta_4 + \eta_4 \eta_6 + \eta_6 \eta_2
=\rho^2 . $$
\vspace{0.5cm}

\noindent{\bf Bicentric hexagons}
An elaborate formula for the area of a cyclic hexagon in terms of side lengths is given
in~\cite{Ro95}.
\goodbreak
\subsection{$(\rho,R_V,L)$ and generalizing Blundon's inequality}\label{subsec:rLRV}

We hope to treat bicentric quadrilaterals in the future, but, for now,
here are results for triangles.

\subsubsection{Triangles}

\begin{figure}[ht]
\centerline{\includegraphics[height=10cm,width=14cm]{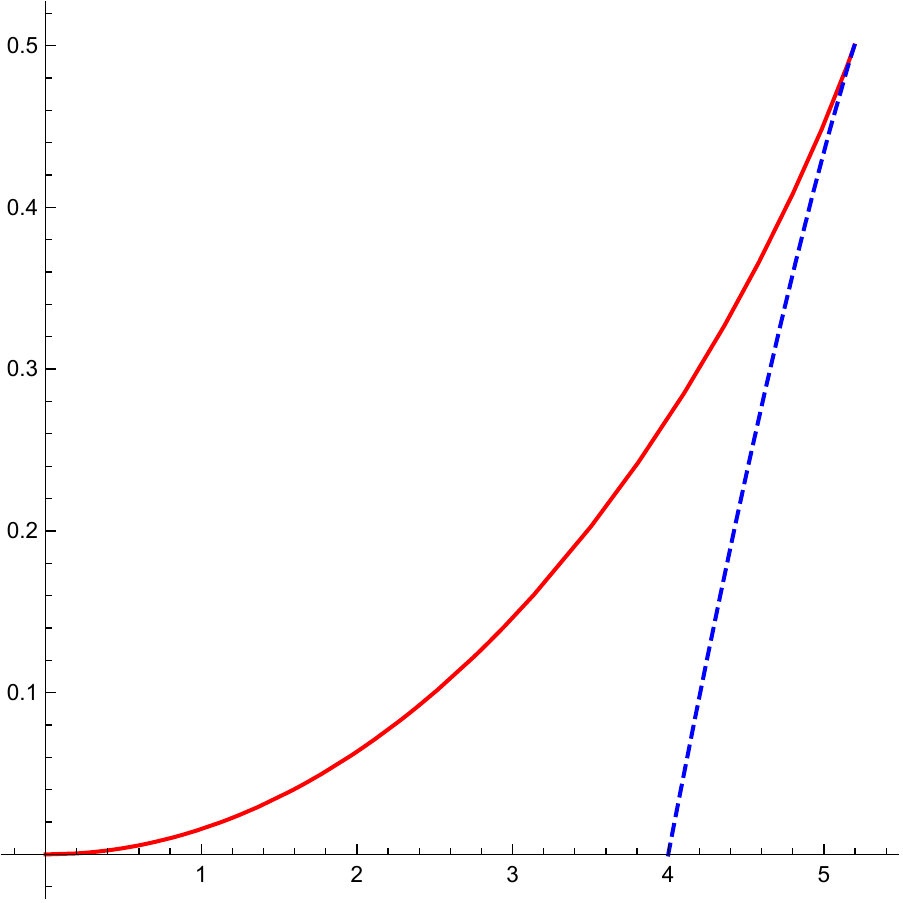}}
\caption{The lhs of the Blunden inequality is shown in red,
the rhs in blue dashed.
}
\label{fig:IIbBlund}
\end{figure}

One of the many entries in wikipedia's list of triangle inequalities is
the following, which in much of the literature is known as
Blundon's inequality:
{\small
\begin{equation}
2R_V^{2}+10R_V\rho-\rho^{2}-2\sqrt{R_V}(R_V-2\rho)^{3/2}
\leq \frac{L^{2}}{4}
\leq
2R_V^{2}+10R_V\rho-\rho^{2}+2\sqrt{R_V}(R_V-2\rho)^{3/2} .
\label{in:Blund}
\end{equation}
}
Quoting from wikipedia:
\begin{quote}
Here the expression
$$\sqrt {R_V^{2}-2R_V\rho}
={\rm dist(incentre,circumcentre)},$$
In the double inequality~(\ref{in:Blund}), the first part holds with equality if and only if the triangle is isosceles with an apex angle of at least $\pi/3$,
and the last part holds with equality if and only if the triangle is isosceles with an apex angle of at most $\pi/3$.
Thus both are equalities if and only if the triangle is equilateral.
\end{quote}
Let $x=L/R_V$, $y=\rho/R_V$. The allowed values of the $(x,y)$ pair is the region inside
the curves given by
$$4 (2+10 y - y^2 -2(1-2 y)^{3/2})
\le x^2
\le
4 (2+10 y - y^2 +2(1-2 y)^{3/2}) ,
$$
and shown in Figure~\ref{fig:IIbBlund}.

Some of the behaviour of isosceles triangles for the various
Blaschke-Santalo diagrams is indicated in the table below.

\begin{center}
\begin{tabular}{|c|c|c|}
\hline
 & squat &  tall\\
\hline
$L/R$ & $\sim{4}$&$\sim{4}$\\
$\rho/R$& $\rightarrow{0}$& $\rightarrow{0}$\\
\hline
$L/R$ & $\sim{4}$&$\sim{2}$\\
$\rho/R$& $\rightarrow{0}$& $\rightarrow{0}$\\
\hline
$L/R_V$ & $\sim{0}$&$\sim{4}$\\
$\rho/R_V$& $\rightarrow{0}$& $\rightarrow{0}$\\
\hline
\end{tabular}
\end{center}
\clearpage

%
\section{Moments about the centroid, $I_c$}\label{sec:centroids}

The moment of inertia about the centroid, $I_c$ is,
amongst $n$-gons with a given area, minimized by the regular $n$-gon:
\begin{eqnarray*}
I_c &\ge& I_c({\rm reg-n-gon}) ,\\
&=&\frac{1}{4} \rho({\rm reg-n-gon})\, i_2({\rm reg-n-gon}) , \\
&=&\frac{A^2 (3+\tau_n^2)}{6 n \tau_n}\ {\rm where\ } \ 
\tau_n = \tan(\frac{\pi}{n}) .
\end{eqnarray*}
(The calculations for the equations above are in Part IIa.)

The St Venant inequality is the left-hand side below
$$ \frac{\QQ_0}{A^2}\le \frac{1}{8\pi} \le \frac{I_c}{4 A^2}  \le \frac{I_2}{4 A^2} ,$$
where, as before, $I_2$ is the moment about the incentre,
and, for tangential polygons, $I_2=\rho i_2/4$.

We can add this inequality into that near the end of Part I~\S\ref{subsect:tangNew} that
for any tangential polygon
$$\frac{1}{8}\rho^2 A \le \QQ_{0-}\le\QQ_0\le\frac{1}{4}I_c\le\frac{1}{4}I_2\le\QQ_{0+} .$$

\subsection{Formulae for $I_2-I_c=A |z_c-z_I|^2 $}

The identity in the subsection heading is from the
Parallel Axis Theorem.\\
\verb$https://en.wikipedia.org/wiki/Parallel_axis_theorem$
Denote by $z_c$ the coordinates of the centroid,
and $z_I$ those of the incentre.
Take, as elsewhere, the incentre to be the origin.
Write
$$ d_{IG} = |z_c-z_I| . $$

\medskip

\noindent{\bf Triangles}

\noindent
For triangles $d_{IG}$ is given in terms of the side lengths in\\
\verb$https://mathworld.wolfram.com/Incenter.html$\\
Let the sides be denoted by $s_k$ and define
\begin{eqnarray*}
S_1 
&=& s_1+s_2+s_3 , \\
S_2 
&=& s_1 s_2+s_2 s_3+s_3 s_1 ,\\
S_3 
&=& s_1 s_2 s_3 .
\end{eqnarray*}
Then
\begin{equation}
d_{IG}^2 =\frac{ 5 S_1 S_2 -S_1^3 -18 S_3}{9 S_1} . 
\label{eq:dIG3}
\end{equation}

\medskip

\noindent{\bf Quadrilaterals}

For a rhombus the centroid and incentre coincide at the point of intersecion of
the diagonals.

We have yet to derive, for tangential quadrilaterals or special cases of these,
formulae analogous to~(\ref{eq:dIG3}).

See~\cite{KKLY17,KK20} for results concerning centroids, and characterizations
of kites, etc.
See also~\cite{My06,HS09}.

\newpage
\begin{center}
{\large{{\textsc{ Part III:
Other inequalities and properties for $\QQ_0$,\\
isoperimetric inequalities,   etc.
}}}}
\end{center}

\section*{Abstract for Part III}
A referee asked for a bit more on `context' and more numerics.
There is a huge literature on bounds for torsional rigidity and in this part I
focus on that most relevant to convex polygons.
\begin{itemize}
\item
I also return to more on triangles, especially isosceles triangles for which,
unsuprisingly, if one makes use of information specific to triangles, one can improve
on my $\QQ_{0-}$ bound of Part I.
\item
Numerics for rhombi again indicate that  $\QQ_{0-}$ is  quite close to the
actual torsional rigidity.
\end{itemize}


\section{Introduction to Part III}\label{sect:IntroIII} 

A famous open problem is as follows.\\
{\it Amongst all simple polygons with $n$ vertices with a given area, 
does the regular $n$-gon have the greatest torsional rigidity?}\\
\cite{PoS51} establishes this for $n = 3$ and $n = 4$ but the question is open for $n > 4$.
It remains open for tangential $n$-gons too.\\
It will be easier to investigate, for tangential $n$-gons,
 this for the lower bound $\QQ_{0,c}$ rather than the actual torsional
rigidity $\QQ_0$.
The results in this direction are, so far, very slight (just isoceles triangles).
We have yet to establish it for general triangles $n=3$.
\medskip

An outline for this part is as follows.

\begin{itemize}

\item In~\S\ref{sec:BoundsQIso} we present some previously published bounds.
We also introduce the style of isoperimetric inequalities
which will be studied in later sections.

\item  In~\S\ref{sec:genQ0} we review some bounds on $\QQ_0$.

\item In~\S\ref{sec:MoreIsos} we study isosceles triangles
There are many questions on
how the geometry affects the domain functionals.
We have, in~\S\ref{app:isosIsop}, some isoperimetric results: some functionals are optimized,
over isosceles triangles with given area, by the equilateral triangle.

\item In~\S\ref{sec:GenTri} (yet to be written!)
 we will, very briefly, consider extending the work
on $\QQ_{0-}$ to general triangles.

\item In~\S\ref{sec:RhombiQ} we show, numerically for rhombi, how close $\QQ_{0-}$
is to $\QQ_0$.\\
This parallels the work on isosceles triangles presented at the end of 
Part I.

\item In~\S\ref{sec:BSQ} we consider Blaschke-Santalo diagrams involving
$\QQ_0$, or $\QQ_{0-}$, for tangential polygons.

\item In~\S\ref{sec:QFurtherQs} I give further questions in connection
with $\QQ_{0-}$.
Some of these are of the form: if some property has been established
for $\QQ_0$, does $\QQ_{0-}$ have the same property.
Domain monotonicity is one such property, and one question is
whether it remains, for $\QQ_{0-}$
under corner cutting (see Part IIb~\S\ref{sec:Transformations}).
Another, not discussed there, is whether, amongst tangential $n$-gons
with a given inradius, the regular $n$-gon minimizes $\QQ_{0-}$.

\end{itemize}


\section{Bounds on $Q_0$, esp. isoperimetric}\label{sec:BoundsQIso}

Other papers involving the torsional rigidity of tangential polygons 
include~\cite{CFG02}\S{2.2} involving web functions
and~\cite{Sal18}.
Web functions are particularly simple for tangential polygons,
and align with the similar level curves of~\cite{PoS51},
level curves which are the same shape as the boundary.
In~\cite{Sal18} the bounds are in terms of
$$ I(q,\partial\Omega)
= \int_\Omega {\rm dist}(z,\partial\Omega)^q . $$
(The notation is that of~\cite{Ke07} as, in this Supplement,
there are already other uses for $I$ and $i$.
The capital $I$ is, as with our other use, an integral over the domain.)
For tangential polygons $I_0(\partial\Omega)=|\Omega|$ and
$$ I(q,\partial\Omega) = \frac{(p+1)(p+2)}{(q+1)(q+2)} I(p,\partial\Omega)\rho(\Omega)^{q-p}
= \frac{1}{(q+1)(q+2)}\, |\Omega|\, \rho(\Omega)^q .$$
We specialise a much more general theorem of \cite{Sal18} to the following.\\
{\bf Theorem~\cite{Sal18}}. 
{\it Let $\Omega$ be a tangential polygon.
Then
$$\QQ_0(\Omega) 
\ge \frac{1}{16} (p + 1)(p + 2)I(p,\partial\Omega)\,\rho(\Omega)^{2-p}
\qquad{\rm  where}\  -1 \le p < \infty.$$
Equality holds if $\Omega$ is a disk.}\\
In particular, with $p=2$
$$ \frac{1}{2} |\Omega|\rho^2 = 3 I(2,\partial\Omega) \le 4\QQ_0(\Omega) , $$
which recovers~\cite{PoS51}\S5.8 equation(7) on p100.
See inequality~(\ref{in:PS}) in Part I.
\cite{Sal18} is largely concerned with how this might generalize to domains
other than tangential polygons.

\medskip

The famous open problem stated at the beginning of this part  is as follows.\\
{\it Amongst all simple polygons with $n$ vertices with a given area, 
does the regular $n$-gon have the greatest torsional rigidity?}\\
Repeating from there:
\cite{PoS51} establishes this for $n=3$ and $n=4$ but the question is open
for $n>4$.

Of course the question can be asked with a smaller sets,
e.g. convex polygons or tangential polygons or bicentric polygons or ...
To date there are no answers for the torsional rigidity.
There are, however, for some other domain functionals: see Part IIb
for purely geometric ones. 
The conformal inradius and related radii are others that on fixing the area is extreme for the regular $n$-gon (see~\cite{SZ04}).
By a $\cal D$ we shall mean some domain functionals, and we are interested in pairs of these for which one has a result of the form\\
{\it For tangential $n$-gons with fixed ${\cal D}_1$ the regular $n$-gon 
$<$maximizes$|$minimizes$>$ ${\cal D}_2$}\\
Table~\ref{tbl:tbl1arIso} below presents some
(repeating some entries from Table~\ref{tbl:tbl1ar} of Part IIb):

\begin{table}[h]
\begin{center}
\begin{tabular}{|| c | c | c | c ||}
\hline
${\cal D}_1$& & ${\cal D}_2$& Remark \\
\hline
area & min& perimeter& \\
area & max& inradius& $A=\rho L/2$  \\
inradius& min& perimeter& Jensen inquality\\
inradius& min& area& " \\
inradius& min& $Q_0$& Solynin~\cite{Sol92,SZ10}\\
 & & &see below\\
\hline
\end{tabular}
\caption{Tangential polygons}
\label{tbl:tbl1arIso}
\end{center}
\end{table}

In~\cite{Sol92} the result, at fixed inradius $\rho$, that
$\QQ_0$ is minimized at the regular $n$-gon is first given,
in Theorem 1, for tangential $n$-gons,
and, after that, in Theorem 2, for more general $n$-gons.
(I have yet to check the proofs.)

\vspace{0.5cm}

\subsection{Bounds on $Q_0$ for convex $\Omega$}

We have, in Part I, reported the results, from~\cite{PoS51}
$$ \frac{A^2}{4B}\le \QQ_0 \le \frac{A^2}{8\pi} .$$
The left-hand inequality is in Part I, inequality~(\ref{in:PS}):
the right-hand inequality is the St Venant Inequality.
There is equality in both for disks.
\cite{PoS51} p99 gives, for convex domains,
$$ B\le \frac{2A}{\rho^2} ,$$
which is an equality for tangential polygons.
Combining these gives 
for convex domains, the left-hand inequality in
$$\frac{1}{8}\rho^2 A
\le \QQ_0 \le c \rho^2 A . $$
The left hand inequality is from~\cite{PoS51}\S5.8~equation~(7), and
equality occurs for $\Omega$ a disk.
The right hand inequality with $c=\frac{1}{3}$ is one of several due to Makai, and
equality is approached by rectangles as they become  long and slender,
leaving the question as to the best constant $c$ when restricted to 
tangential polygons.

From~\cite{PoS51} p254, for a thin rectangle
$$ \QQ_0 \sim\frac{1}{3} \rho^2 A
\ \ {\rm for}\ \ \rho\rightarrow{0} . $$

For thin isosceles triangles, with small vertex angle,
$$ \QQ_0 \sim \frac{1}{6} \rho^2 A
\ \ {\rm for}\ \ \rho\rightarrow{0} . $$
For the lower bound found in Part I 
$$\QQ_{0,-} \sim \frac{21}{128} \rho^2 A
\ \ {\rm for}\ \ \rho\rightarrow{0} .$$
This checks with $\QQ_{0-}\le\QQ_0$.



Another inequality for convex domains is
$$\frac{1}{3}\frac{A^3}{L^2}
\le \QQ_0
\le \frac{2}{3}\frac{A^3}{L^2} . $$
The left inequality is approached by thin rectangles.
The right inequality is approached by thin acute isosceles triangles.
(It is curious that thin rectangles occur as upper bounds in one
inequality in this section, and lower bounds in another.)
However, inequalities for convex sets,
$$ \frac{\rho}{2}\le \frac{A}{L}\le \rho(1-\frac{\pi\rho}{L})\le \rho ,$$
(in which the central inquality is equality for a disk,~\cite{SA00})
are consistent with
$$\frac{1}{3}\frac{A^3}{L^2}
\le \QQ_0
\le \frac{1}{3}\rho^2 A . $$

For tangential polygons the inequalities above give
$$\frac{1}{2}\frac{A^3}{L^2}=\frac{1}{8}\rho^2 A
\le Q_0 \le \frac{1}{6} \rho^2 A = \frac{2}{3}\, \frac{A^3}{L^2} . $$
The extreme domains are the disk (left) and
thin isosceles triangles (right).

See also~\cite{BBP20}.

\subsection{Cheeger constant and $\QQ_0$}

We have
$$ \QQ_0(\Omega) h(\Omega)^4
\ge \QQ_0(B) h(\Omega)^4 = 2\pi.
$$
Equality occurs only for disks.
See~\cite{LMR18}, and specialize to 2-dimensions.
Thus for our tangential polygons
$$\QQ_0(\Omega)\ge \frac{32\pi A^4}{(L+\sqrt{4\pi A})^4} .$$

\vspace{4cm}

Much is known about the elastic torsion function and,
in particular, its properties in convex domains.

It is known that in convex $\Omega$, the square root of the torsion function
$\sqrt{u_0}$ is concave.
In some proofs of this one uses that for solutions of the torsion equation
$\log(1-4H)$, with H the hessian determinant, is harmonic.
$u_0$ is, itself, not concave in any domain with corners.
However one can ask on what subset $\Omega_1$ of $\Omega$ is $u_0$ concave.
For an equilateral triangle (and for a disk), $\Omega_1$ contains the incircle.
(See~\S{11} of~\cite{KMtorsion}.)
One wonders what additional conditions, if any, might be needed for
it to happen in other tangential polygons that
$\Omega_1$ contains the incircle.

Improvements on the St Venant inequality involving Fraenkel asymmetry are given 
in~\cite{BP16}.
\medskip 

There is a Urysohn’s type inequality which we denote by (U).\\
(U): among convex sets with given mean width, the torsional rigidity is maximized by balls.\\
The mean width $w$ of any compact shape $\Omega$ in two dimensions is 
$L/\pi$, where $L$ is, as before, the perimeter of the convex hull of 
$\Omega$. 
So $w$ is the diameter of a circle with the same perimeter as the convex hull.
So in (U) it is just maximizing with same perimeter. 
Inequality (U) is weaker than St Venant.

There are inequalities involving polar moment of inertia about the centroid.\\
{\it For convex domains $\QQ_0 I_c A^{-4}$ is maximized by its value for an 
equilateral triangle.}
(Related results, if not exactly this, are in~\cite{Pol55}.)\\
We conjecture:\\
{\it For tangential polygons $\QQ_0 I_2 A^{-4}$ is maximized by its value for an 
equilateral triangle.}

Sperb $u_{\rm max}\le\rho^2$.
Strip $u_0 = (\rho^2-x^2)/2$

\section{Some bounds on $\QQ_0$}\label{sec:genQ0}

Recall the well-known estimates
\begin{eqnarray*}
\QQ_0
&\le& \frac{\pi}{8}\, \left(\frac{A}{\pi}\right)^2 = \QQ_{\odot,0} ,\\
\QQ_0
&\le& \Sigma_\infty .
\end{eqnarray*}
Define, in the notation of~\cite{PoS51},
\begin{equation}
B_\Omega
= \int_{\partial\Omega} \frac{1}{x\cdot{n}} .
\label{eq:Bdef}
\end{equation}
Some well known lower bounds are:
\begin{eqnarray*}
\QQ_0
&\ge& \frac{A^2}{4\ B_\Omega}  ( = {\rm{ \ for\ tangential\ polygons\ }} \frac{A^3}{2L^2} )
 ,\\
\QQ_0
&\ge& \frac{\pi}{8}\, {\dot{r}}^2 ,
\end{eqnarray*}
where ${\dot{r}}$ is the maximum interior mapping radius of $\Omega$.

\section{More on isosceles triangles}\label{sec:MoreIsos}

\subsection{Isoperimetric inequalities for isosceles triangles}{\label{app:isosIsop}}

Numerical data on functionals associated with
isosceles triangles was given in\S\ref{sec:isos}. 
Our main emphasis in the following will be isosceles triangles rather than
general triangles.
However, some general statements are available.
Recall that $\cal A$ is all (simply-connected) domains,
${\cal A}_n$ is all $n$-gons.
\begin{itemize}
\item 
$\QQ_0/A^2$ is maximized over ${\cal A}_3$ (triangles) for an equilateral triangle. 

\item 
$B=L/\rho$ is minimized over ${\cal A}_3$ (triangles) for an equilateral triangle. 
\item 
${\dot r}/\sqrt{A}$ is maximized over ${\cal A}_3$ (triangles) for an equilateral triangle. 
\end{itemize}

\subsubsection{$B$}\label{app:BIsosIsop}

$B=L/\rho$ for any tangential polygon, so for isosceles triangles
in particular.\\
The disk minimizes $B$ and $\QQ_0 B/A^2$ over tangential polygons.\\
The equilateral triangle minimizes $B$ and $\QQ_0 B/A^2$ over triangles.

\begin{tabular}{|| c | c | c | c | c | c  ||}
\hline
shape&$n$& $8\pi\QQ_0/A^2$& $\QQ_0/(A/\pi)^2$&$B=\frac{L}{\rho}$&$4\QQ_0 B/A^2$
  \\
\hline
disk&$\infty$& 1& $\pi/8$& $2\pi$& $1$
 \\
hexagon&6& 0.9643& 0.3786&  $4\sqrt{3}$&$1.063$
 \\
square& 4& 0.8834& 0.3469& $8$&$1.125$
 \\
\hline
equilateral $\Delta$&3& 0.7255&  0.2849&  $6\sqrt{3}$&$1.200$
 \\ 
right isosceles& &0.6557& 0.2575& $4(1+\sqrt{2})$&$1.217$
 \\
\hline
\end{tabular}

\subsubsection{$\dot r$}\label{app:rdotIsosIsop}

\bigskip
The calculation of $\dot r$ for polygons often involves the use of Schwarz-Christoffel conformal mappings.
Some exact values are given on p273 of~\cite{PoS51}. 

\begin{itemize}

\item (Aside, except for $n=3$.)
For a regular polygon with $n$ sides, and perimeter $L_n$,
$$ {\dot r}_n
=\frac{\Gamma(1-\frac{1}{n})}{2^{1+2/n}\, \Gamma(\frac{1}{2})\, \Gamma(\frac{1}{2}+\frac{1}{n})}\, L_n . $$

\item
Again citing~\cite{PoS51} p158, amongst all triangles with a given area that which maximizes $\dot r$ is
equilateral.

\item
For (regular polygons and) triangles we have
$$\pi {\dot r} {\bar r} = A, $$
$A$ being the area and $\bar r$ the transfinite diameter.

\item
Working from earlier more general results (Haegi 1951) in~\cite{Fi14rdot} it is given that
for an isosceles triangle, vertex angle $\alpha$ and base $2\sin(\alpha/2)$ the transfinite
diameter, denoted there by $\kappa$, is
$$\kappa(\alpha) = \frac{\sqrt{\pi+\alpha}}{8\pi^{5/2}}
\left(\frac{\pi+\alpha}{4\alpha}\right)^{\alpha/(2\pi)}
\, \frac{\sin(\alpha)^2}{\sin(\alpha/2)}\, \Gamma(\frac{\alpha}{\pi})\, \Gamma(\frac{\pi-\alpha}{2\pi})^2 .
$$
\end{itemize}

For isosceles triangles with area $A$ and vertex angle $\alpha$
$${\bar r} = \frac{2\sqrt{A\tan(\alpha/2)}}{\sin(\alpha/2)}\, \kappa \ \ {\rm so\ \ }
{\dot r} = \frac{A}{\pi{\bar r}} = \sqrt{A}\, \frac{\sin(\alpha/2)}{2\pi\sqrt{\tan(\alpha/2)}}\, \frac{1}{\kappa} .
$$

Some values of $\dot r$, as given in~\cite{PoS51}, are copied in the following table.

\begin{tabular}{|| c | c | c | c | c | c  | c ||}
\hline
shape&$ $& $\dot r$& $\dot r$&$\frac{\dot r}{\sqrt{A}}$&
$8\QQ_0/(\pi{\dot r}^4)$  \\
\hline
disk&radius $a$& $a$& $a$& $0.56419$& $1$
 \\
hexagon&side $s_6$& $\frac{2^{5/3}\, \sqrt{3}\pi}{\Gamma(1/3)^3}\, s_6$& $0.89850 s_6$&  
$0.55744$&$1.011$
 \\
square&side $s_4$& $\frac{4\sqrt{\pi}}{\Gamma(1/4)^2}\,s_4$&$0.53835 s_4$&
$0.53935 $& $1.058$
 \\
\hline
equilateral $\Delta$&side $s_3$ & $\frac{2\pi}{\Gamma(1/3)^3}\, s_3$&
 $0.3268\, s_3$ &  
$0.49665$&$1.209$
 \\ 
right isosceles& equal sides $a$&$\frac{4\sqrt{2\pi}}{3^{3/4}\, \Gamma(1/4)^2} \, a$& $0.33462 a$& 
$0.47320$&$1.325$
 \\
\hline
\end{tabular}

\subsubsection{Simple comments on torsion for isosceles triangles}\label{subsec:isosQ0}

As noted before, for an equilateral triangle  the torsion function is a cubic polynomial in $x$ and $y$ corresponding to forming the products of three linear terms, each linear term being 0 on one side of the triangle. 
There appear to be no other simple solutions for isosceles triangles, though there is a series formula for the torsional rigidity of the right isosceles triangle. 

In the table below, in the second, third and fourth columns
we take the base to be 2, the area to be $h$, $h$ being the height;
in the fifth and sixth columns the area is $\sqrt{3}$,
the height is $h$ and the base $2\sqrt{3}/h$.
\bigskip

\begin{tabular}{|l|c|c|c|c|c|}
\hline
infinitely acute & $h \rightarrow\infty $&$A\sim{h}$& $\QQ_0\sim \frac{1}{6} h $&
 $h \rightarrow\infty $& $\QQ_0\sim\frac{1}{2h}$\\ 
equilateral & $h=\sqrt{3}$& $A=\sqrt{3}$ &$\QQ_0=\frac{\sqrt{3}}{20}$&
$h=\sqrt{3}$& 
$\QQ_0=\frac{\sqrt{3}}{20}$\\
 right isosceles &$h=1$& $A=1$&$\QQ_0 = 0.026091$ & 
 $h=3^{1/4}$& 
  \\ 
 infinitely flat & $h\rightarrow{0}$&$A\sim{h}$&, $\QQ_0\sim \frac{1}{24} h^3$&
  $h\rightarrow{0}$& 
  $\QQ_0\sim\frac{h}{8}$ \\ 
 \hline
 \end{tabular}

\medskip
 

\goodbreak

\subsubsection{More variational bounds on $\QQ_0$ for isosceles triangles}

\noindent
The bounds on $\QQ_0$ we report here are often very old.

There are, of course, variational formulations of the Problem (P(0)). 
For positive, differentiable functions $v$, let
$$
E(v, 0, k) = \int_\Omega v^k \qquad{\rm  and\ \ } E(v,1,k) =  \int_\Omega |\nabla v|^k \  .
$$
With
\begin{equation}
\QQ_{0,\rm LB}(v)= E(v,0,1)^2 / E(v,1,2)                                 
\label{eq:e2}
\end{equation}    
it can be shown that for any smooth function $v$ vanishing on the boundary of $\Omega$, 
the expression $\QQ_{0,\rm LB}(v)$  provides a lower bound on the torsional rigidity 
$\QQ_0=\QQ_{0,\rm LB}(u)$. 
The theory for this is given in~\cite{PoS51}.

For isosceles triangles,  base $2a$, height $h$, 
a simple trial function $v$ to substitute into (\ref{eq:e2}) is
\begin{equation}
v_{\rm cub} = (y+\rho)  (1-(x/a)-((y+\rho)/h))  (1+(x/a)-((y+\rho)/h)) .   
\label{eq:e3}
\end{equation}          
Evaluating the integrals gives
\begin{equation}
\QQ_{0\rm{cub}}=\QQ_{0,\rm LB}(v_{\rm cub} )= \frac{a^3 h^3}{(30 a^2 + 10 h^2)}
= \frac{A^2}{30\tau +10/\tau},
\label{eq:e4}
\end{equation}   
with $\tau=a/h=\tan(\alpha/2)$ as before,
and the expression on the right of (\ref{eq:e4}) is, in fact, 
an approximation given in Roark's tables, valid 
for the the vertex angle $\alpha$ of the isosceles triangle lying between 40 and 80 degrees. 
It is exact for the equilateral case where the vertex angle is 60 degrees.
It would, we think, be an improvement to Roark's tables to have noted, 
after defining this rational function of $a$ and $h$, 
to have noted that $\QQ_0>\QQ_{0,\rm LB}(v_{\rm cub})$  for all positive values of $a$ and $h$, 
and after this to have noted the range of $h/a$ over which it provides a good approximation to $\QQ_0$.

We can compare this lower bound with our earlier $\QQ_{0-}$ as shown in
Figure~\ref{fig:plIsos}.
In the next figures the bound of~(\ref{eq:e4}) is shown brown, dashed.
Unsurprisingly it is good for triangles near equilateral improving on
$\QQ_{0-}$ there, but it is worse than $\QQ_{0-}$ when not near equilateral.
For example, for a right isosceles triangle, with area $\sqrt{3}$ we found
$$\QQ_0 =0.07827,\qquad
\QQ_{0-}=0.07651,\qquad
\QQ_{0\rm{cub}} = \frac{3}{40}= 0.075 .
$$

\begin{figure}[ht]
\centerline{\includegraphics[height=10cm,width=14cm]{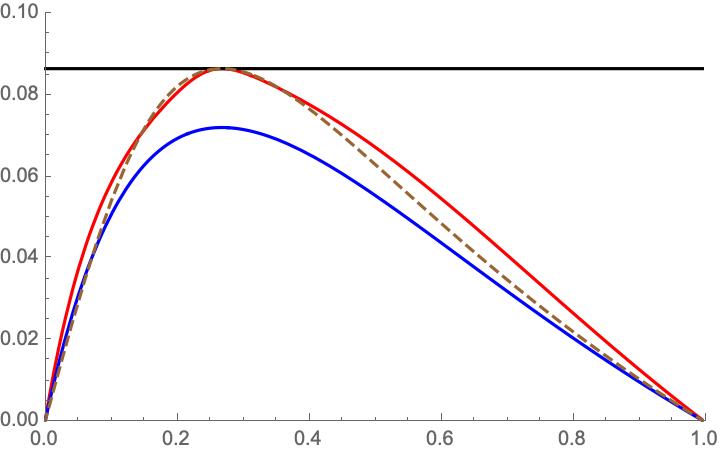}}
\caption{For an isosceles triangle with area$\sqrt{3}$.
$\sigma$ is tan of 1/4 of the apex angle.
Blue is $\QQ_B$, red is the new lower bound $\QQ_{0-}$,
black is $\QQ_\Delta$, brown dashed is the cubic one underdiscussion.
 }
\label{fig:plQ0LBcube1}
\end{figure}

\begin{figure}[ht]
\centerline{\includegraphics[height=10cm,width=14cm]{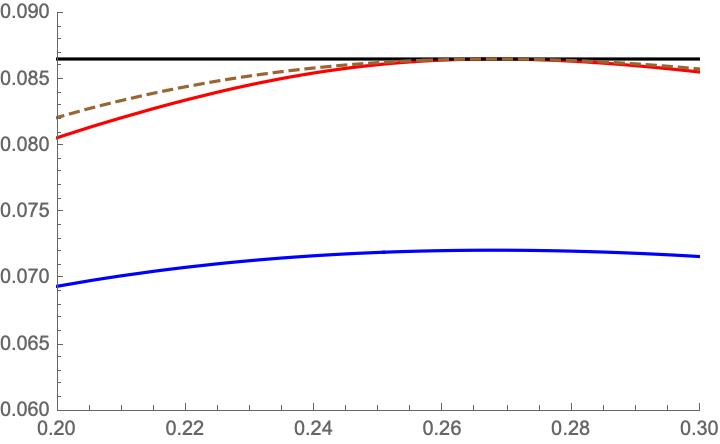}}
\caption{For an isosceles triangle with area$\sqrt{3}$.
$\sigma$ is tan of 1/4 of the apex angle.
Blue is $\QQ_B$, red is the new lower bound $\QQ_{0-}$,
black is $\QQ_\Delta$, brown dashed is the cubic one underdiscussion.
 }
\label{fig:plQ0LBcube2}
\end{figure}
\clearpage

\bigskip
Define $v_{\rm quad}$ to be even in $x$ and
$$v_{\rm quad} = y  (1+(x/a)-(y/h)) 
\ \ {\rm for}\ \ 0 < x < a. $$
Good lower bounds on $\QQ_0$ for when $h$ is small can be found using this as test function. 

See PartI, near Figure~\ref{fig:IIasolQ0m} for another lower bound on $\QQ_0$.
Many other bounds are available. See for example the Appendix, by Helfenstein of the paper~\cite{PSH54}. 

{\bf ToDo.} Check the  Helfenstein work.
Stretching in one direction
Define 
$$D(1,t)= \{ (tx,y) | (x,y) \in D \}. $$
In particular, $D(1,1)=D$. 
Domain monotonicity gives $QQ_0(D(1,t))$ is increasing in $t$.
In~\cite{PSH54} it is shown that  $t/\QQ_0(D(1,t))$ is increasing and concave in $t^2$. 
The appendix to their paper by Helfenstein makes good use of this in connection with isosceles triangles, finding upper and lower bounds on $\QQ_0$ which differ by no more than 12\%.
 These are, unfortunately, rather elaborate, and with cheap numerical computing, it is perhaps better to use the stretching result as yet another check on numerics.
 
 {\bf ToDo.} Main interest is in $\QQ_0(D(1,t))/|D(1,t)|^2$ where the maximum occurs for $t$ giving an equilateral triangle.
\bigskip

\subsection{Miscellaneous other bounds}

Consider an isosceles triangle with a base of length 2 and 
height of length $h$, and angles $\beta$, $\beta$, and $\pi-2\beta$ respectively.
Since $\tan(\beta)=h$, and the area $A$ of the triangle is $h$, it can be shown that
an inequality from~\cite{BFNT20} is
$$Q_0 \le \frac{1}{8} (1+A^2)^2\, (A-{\rm arctan}(A)) . $$
This inequality is, in~\cite{BFNT20}, used with triangles which are thin,
and, while satisfied for equilateral, it is weak there.

\section{General triangles}\label{sec:GenTri}

From the calculations of $L$, $i_2$ nad $i_4$
in~Part IIb~\S\ref{sec:IIbTri}, $\Sigma_\infty$ and $\Sigma_1$ can be found.
\section{Rhombi}\label{sec:RhombiQ}

The formulae for $A$, $L$, $i_2$, $i_4$,
$\Sigma_\infty$, $\Sigma_1$
are given in Part~IIa~\S\ref{sec:tangQuad}.
The numerical values of $\QQ_0$ are from~\cite{SC65}
(and checked with Mathematica NDSolve).

\begin{figure}[ht]
\centerline{\includegraphics[height=10cm,width=14cm]{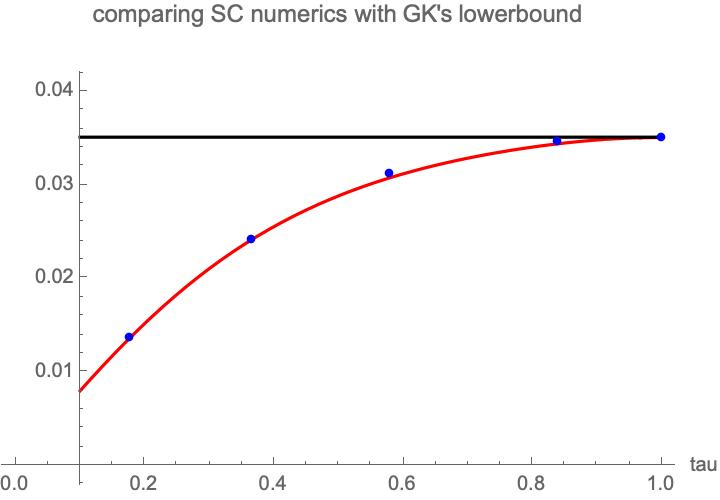}}
\caption{
Rhombus
 }
\label{fig:plRhomb}
\end{figure}
\clearpage



\section{Blaschke-Santalo for $\QQ_0$}\label{sec:BSQ}

If one were to assemble, for tangential polygons with given $\rho$ and $d_O$,
a Blaschke-Santalo diagram, I would expect that with $\QQ_0$ 
would be similar to that with the area (or equivalently the perimeter) given
earlier. 
In particular, at given $\rho$ and $d_O$,
\begin{itemize}
\item the 1-cap would minimize $\QQ_0$ thereby providing the left-upper boundary in the diagram, and
\item the regular polygons would be amongst the maximizing shapes thereby
providing a scatter of points on the right-lower boundary.
\end{itemize}

There may be independent interest in curves, or
other subsets, in the region, e.g. for
\begin{itemize}
\item triangles

\item kites/rhombi\\
other quadrilaterals.

\end{itemize}

A much easier task is to find the diagram for the lower bound $\QQ_{0-}$.
I expect everything stated above for $\QQ_0$ would occur in this
simpler situation.

\section{Further questions}\label{sec:QFurtherQs}

\noindent
{\bf Question.} Does `corner cutting' 
(defined in~\S\ref{sec:Transformations})
increase $\QQ_0/A^2$?\\
It would be (much) easier to begin with
$\QQ_{0-}$ and isosceles triangles 
each having its apex corner cut
to form a symmetric tangential trapezium.

Given the area of a regular $n$-gon, its inradius is,
as in Part IIa~\S\ref{sec:regn},
$$ \rho_n^2 = \frac{A}{n\tau_n} \qquad {\rm with\ \ }
\tau_n =\tan(\frac{\pi}{n}) . 
$$
Formulae for $i_{2,n}$, $i_{4,n}$, $\Sigma_{1,n}$ and $\Sigma_{\infty,n}$
are also given in Part IIa~\S\ref{sec:regn}.
Let $f_n(\QQ_0)$ be the quadratic in $\QQ_0$
defined at Part I, equation~(\ref{eqt:fDef})
using the values appropriate to the regular $n$-gon for $\rho$, $i_2$ and $i_4$.
The leading coefficient, i.e. coefficient of $\QQ_0^2$,
is the same for the two polynomials.

For an equilateral triangle with area $\sqrt{3}$, we found in Part I
$$ f_3(\QQ_0) = f_\Delta(\QQ_0)
= 32\sqrt{3} (\QQ_{0-}\frac{\sqrt{3}}{20})(\QQ_0 - \frac{\sqrt{3}}{4}) .
$$

We wish to show that, for any tangential $n$-gon with the same area
$$ \QQ_{0-}\le \QQ_{0-,n} .$$
This suggests that we try to show
$$ f(\QQ_{0,n})\le{0}, \qquad f_n(\QQ_0)\ge{0} .
$$

This leads to the question of how roots of a quadratic change
when its coefficients change.
Both polynomials $f$ are of the form
$$ 32 A(\QQ_0^2 +b*\QQ_0 +c) = 0 ,$$
with two positive roots
$$ 2\sqrt{c}<b<0,
\qquad 0<c<\frac{b^2}{4} .$$
If we could show that both coefficients $b$ and $c$ are largest for
the regular $n$-gon case, this would establish
that $\QQ_{0,n}$ was larger than any other tangential $n$-gon.

It would seem prudent
\begin{itemize}
\item to begin with $n=3$,
\item then $n=4$, initially with special cases (rhombi, kites, etc.)
\item then $n=5$,
\item then $n=6$,
\item then general $n$.
\end{itemize}
\newpage
\begin{center}
{\large{{\textsc{ Part IV:
Robin boundary conditions $\QQ(\beta)$,\\
isoperimetric inequalities,   etc.
}}}}
\end{center}

\section{Robin boundary conditions}

Part IV concerns a different problem than the earlier parts.

The function $R$ of~\cite{Ke20i} is a function of a nonnegative parameter
$\beta$ and geometric functionals $A$, $L$, $\Sigma_\infty$ and $\Sigma_1$
and $\QQ_0$.
Knowledge of these, for triangles, 
(as in Part II) is essential for the following questions.
\smallskip

Let $\beta\ge{0}$.
For solutions $u_\beta$ of Problem (P($\beta))$ of~\cite{Ke20i} in a triangle
define $\QQ_{\rm triangle}(\beta)$ as the integral of $u_\beta$ over the triangle.\\

\noindent{\bf Question.} {\it Consider triangles of fixed area, say $\sqrt{3}$.
Amongst these triangles is that which has the greatest  $\QQ(\beta)$ the
equilateral triangle?}\\
(For $\beta=0$ this is the case as proved in~\cite{PoS51}.)
(A similar question, but for the fundamental frequency is asked 
in~\cite{LS17}.)

\medskip
In~\cite{Ke20i} we presented a lower bound $R(\beta,\ldots)$ for $\QQ(\beta)$
and noted that it provided a good approximation to $\QQ(\beta)$.
For triangles the only argument of $R$ that isn't a relatively simple function
of the triangle's geometry is the torsional rigidity $\QQ_0$.
This leads on to simpler questions (for which some of Part III might be relevant).

\smallskip
\noindent{\bf Question.} {\it Consider triangles of fixed area, say $\sqrt{3}$.
What additional properties of the torsional rigidity $\QQ_{\rm triangle}(0)$
will ensure that for these triangles that which has the greatest  $R(\beta,\ldots)$ is the
equilateral triangle?}

\smallskip
\noindent{\bf Questions.} {\it Consider triangles of fixed area, say $\sqrt{3}$.
Replacing, in $R(\beta,\ldots)$, $\QQ_{\rm triangle}(0)$ by one of its bounds or approximations,
will ensure that for these triangles that which has the greatest  
$R(\beta,\ldots,\QQ_{0,{\rm approx}})$ is the
equilateral triangle?}\\
(The plural in `Questions' is because there are several possible bounds that might be used.)
\smallskip

\noindent There is some indication via numerics that, restricting to isosceles triangles,
equilateral triangles maximize.

\bigskip

More generally, similar questions can be asked for $n$-gons.
For these, does the regular $n$-gon optimize?
\newpage
\begin{center}
{\large{{\textsc{ Part V:
Miscellaneous PDE and tangential polygons
}}}}
\end{center}

Tangential polygons get a mention in other pde papers.
\bigskip

Papers by Solynin include treatment of the fundamental frequency
for tangential $n$-gons.
\bigskip

Consider evolution of temperature satisfying the heat equation,
zero Dirichlet data on boundary of 
a convex plane domain $C$, and unit initial value.
In the convex domain $C$ there will, at each fixed time, be a unique point 
at which the solution is a maximum, hot-spot.
\cite{MS08} characterises the tangential polygons which have stationary 
hot-spots, basically as the regular polygons.
We remark that the integral of the solution to the heat equation,
integrated w.r.t. time $t$ from 0 to infinity, solves the torsion problem.
So the unique maximum of the torsion function would have to be the
stationary hot-spot.
They must coincide with the unique critical point of the first eigenfunction.
What further restrictions are needed on a tangential polygon in order that
the incentre is a stationary hot spot? 

\newpage

\end{document}